\documentclass{emsprocart}  
\usepackage{amssymb, enumerate}
\usepackage{amscd}
\usepackage{xypic}
\usepackage{comment}

\newtheorem{thm}{Theorem}[section]
\newtheorem{cor}[thm]{Corollary}
\newtheorem{prop}[thm]{Proposition}
\newtheorem{lemma}[thm]{Lemma}

\newtheorem*{thm1}{Theorem 1}

\theoremstyle{remark}
\newtheorem{remark}[thm]{Remark}
\newtheorem{example}[thm]{Example}

\theoremstyle{definition}
\newtheorem{definition}[thm]{Definition}

\numberwithin{equation}{section}

\newcommand{\im}{\operatorname{Im}}

\newcommand{\Hom}{\operatorname{Hom}}
\newcommand{\End}{\operatorname{End}}
\newcommand{\Isom}{\operatorname{Isom}}

\newcommand{\id}{{\mathtt{Id}}}

\newcommand{\shHom}{\underline{\operatorname{Hom}}}
\newcommand{\shEnd}{\underline{\operatorname{End}}}

\newcommand{\shAut}{{\underline{\operatorname{Aut}}}}
\newcommand{\shIsom}{{\underline{\operatorname{Isom}}}}

\newcommand{\pr}{{\mathtt{pr}}}

\newcommand{\DR}{\mathtt{DR}}

\newcommand{\vac}{{\mathbf{1}}}

\newcommand{\Ad}{\operatorname{Ad}}
\newcommand{\one}{{1\!\!1}}

\newcommand{\Triv}{\operatorname{Triv}}
\newcommand{\MC}{\operatorname{MC}}

\newcommand{\Def}{\operatorname{Def}}

\newcommand{\Mat}{{\mathtt{Mat}}}

\newcommand{\SIGMA}{{\scriptstyle\Sigma}}
\newcommand{\Alg}{\operatorname{Alg}}

\newcommand{\Stack}{\operatorname{Stack}}
\newcommand{\Cat}{\mathbf{Cat}}
\newcommand{\Top}{\mathtt{Top}}
\newcommand{\Cov}{\mathtt{Cov}}

 \DeclareMathOperator{\cotr}{cotr}
\DeclareMathOperator{\Desc}{Desc} 
\DeclareMathOperator{\ad}{ad}
\DeclareMathOperator{\AlgStack}{AlgStack}
\DeclareMathOperator{\Algd}{Algd}

\newcommand{\stack}{\mathtt{St}}
\newcommand{\desc}{\mathtt{dd}}

\newcommand{\ca}{\mathcal{A}}
\newcommand{\uca}{\underline{\mathcal{A}}}

\newcommand{\cf}{\mathcal{F}}

\newcommand{\cj}{\mathcal{J}}

\newcommand{\cl}{\mathcal{L}}
\newcommand{\co}{\mathcal{O}}
\newcommand{\cu}{\mathcal{U}}
\newcommand{\cv}{\mathcal{V}}

\contact[bressler@math.ias.edu]{Paul Bressler, Institute for
Advanced Study, Princeton, NJ 08540, USA}

\contact[Alexander.Gorokhovsky@colorado.edu]{Alexander Gorokhovsky,
Department of Mathematics, University of Colorado, UCB 395, Boulder,
CO 80309-0395, USA}

\contact[rnest@math.ku.dk]{Ryszard Nest, Department of Mathematics,
Copenhagen University, Universitetsparken 5, 2100 Copenhagen,
Denmark}

\contact[tsygan@math.northwestern.edu]{Boris Tsygan, Department of
Mathematics, Northwestern University, Evanston, IL 60208-2730, USA}

\title[Deformations of gerbes]{Deformations of gerbes on smooth manifolds}

\author[P.~Bressler, A.~Gorokhovsky, R.~Nest, B.~Tsygan]{Paul
Bressler\thanks{Supported by the Ellentuck Fund}, Alexander
Gorokhovsky\thanks{Partially supported by NSF grant DMS-0400342},
Ryszard Nest , Boris Tsygan\thanks{Partially supported by NSF grant
DMS-0605030}}

\begin{document}

\begin{abstract}
We identify the $2$-groupoid of deformations of a gerbe on a
$C^\infty$ manifold with the Deligne $2$-groupoid of a corresponding
twist of the DGLA of local Hochschild cochains on $C^\infty$
functions.
\end{abstract}



\section{Introduction}

In \cite{bgnt1} we obtained a classification of formal deformations
of a gerbe on a manifold ($C^\infty$ or complex-analytic) in terms
of Maurer-Cartan elements of the DGLA of Hochschild cochains twisted
by the cohomology class of the gerbe. In the present paper we
develop a different approach to the derivation of this
classification in the setting of $C^\infty$ manifolds, based on the
differential-geometric approach of \cite{bgnt2}.

The main result of the present paper is the following theorem which
we prove in Section \ref{section: proof of main theorem}.

\begin{thm1}\label{main theorem}
Suppose that $X$ is a $C^\infty$ manifold and $\mathcal{S}$ is an
algebroid stack on $X$ which is a twisted form of $\mathcal{O}_X$.
Then, there is an equivalence of $2$-groupoid valued functors of
commutative Artin $\mathbb{C}$-algebras
\begin{equation*}
\Def_X(\mathcal{S})\cong
\MC^2(\mathfrak{g}_\DR(\mathcal{J}_X)_{[\mathcal{S}]}) \ .
\end{equation*}
\end{thm1}

Notations in the statement of Theorem \ref{main theorem} and the
rest of the paper are as follows. We consider a paracompact
$C^\infty$-manifold $X$ with the structure sheaf $\mathcal{O}_X$ of
\emph{complex valued} smooth functions.   Let $\mathcal{S}$ be a
twisted form of $\mathcal{O}_X$, as defined in Section \ref{tworms}.
Twisted forms of $\mathcal{O}_X$ are in bijective correspondence
with $\mathcal{O}_X^\times$-gerbes and are classified up to
equivalence by $H^2(X;\mathcal{O}^\times)\cong H^3(X;\mathbb{Z})$.

One can formulate the formal deformation theory of algebroid stacks
(\cite{lvdb2,lowen}) which leads to the $2$-groupoid valued functor
$\Def_X(\mathcal{S})$ of commutative Artin $\mathbb{C}$-algebras. We
discuss deformations of algebroid stacks  in the Section
\ref{defalgstack}. It is natural to expect that the deformation
theory of algebroid pre-stacks is ``controlled" by a suitably
constructed differential graded Lie algebra (DGLA) well-defined up
to isomorphism in the derived category of DGLA. The content of
Theorem \ref{main theorem} can be stated as the existence of such a
DGLA, namely $\mathfrak{g}(\mathcal{J}_X)_{[\mathcal{S}]}$, which
``controls" the formal deformation theory of the algebroid stack
$\mathcal{S}$ in the following sense.

To a nilpotent DGLA $\mathfrak{g}$ which satisfies
$\mathfrak{g}^i=0$ for $i<-1$ one can associate its Deligne
$2$-groupoid which we denote $\MC^2(\mathfrak{g})$, see \cite{G, G1}
and references therein. We review this construction in the Section
\ref{dgla}. Then the Theorem \ref{main theorem} asserts equivalence
of $2$-groupoids $\Def_X(\mathcal{S})$ and
$\MC^2(\mathfrak{g}_\DR(\mathcal{J}_X)_{[\mathcal{S}]})$.

The DGLA $\mathfrak{g}(\mathcal{J}_X)_{[\mathcal{S}]}$ is defined as
the $[\mathcal{S}]$-twist of the DGLA
\begin{equation*}
\mathfrak{g}_\DR(\mathcal{J}_X) := \Gamma(X;\DR(\overline
C^\bullet(\mathcal{J}_X))[1]) \ .
\end{equation*}
Here, $\mathcal{J}_X$ is the sheaf of infinite jets of functions on
$X$, considered as a sheaf of topological $\mathcal{O}_X$-algebras
with the canonical flat connection $\nabla^{can}$. The shifted
normalized Hochschild complex $\overline
C^\bullet(\mathcal{J}_X)[1]$ is understood to comprise locally
defined $\mathcal{O}_X$-linear continuous Hochschild cochains. It is
a sheaf of DGLA under the Gerstenhaber bracket and the Hochschild
differential $\delta$. The canonical flat connection on
$\mathcal{J}_X$ induces one, also denoted $\nabla^{can}$, on
$\overline C^\bullet(\mathcal{J}_X))[1]$. The flat connection
$\nabla^{can}$ commutes with the differential $\delta$ and acts by
derivations of the Gerstenhaber bracket. Therefore, the de Rham
complex $\DR(\overline C^\bullet(\mathcal{J}_X))[1]) :=
\Omega^\bullet_X\otimes\overline C^\bullet(\mathcal{J}_X))[1])$
equipped with the differential $\nabla^{can} + \delta$ and the Lie
bracket induced by the Gerstenhaber bracket is a sheaf of DGLA on X
giving rise to the DGLA $\mathfrak{g}(\mathcal{J}_X)$ of global
sections.

The sheaf of abelian Lie algebras $\mathcal{J}_X/\mathcal{O}_X$ acts
by derivations of degree $-1$ on the graded Lie algebra $\overline
C^\bullet(\mathcal{J}_X)[1]$ via the adjoint action. Moreover, this
action commutes with the Hochschild differential. Therefore, the
(abelian) graded Lie algebra
$\Omega^\bullet_X\otimes\mathcal{J}_X/\mathcal{O}_X$ acts by
derivations on the graded Lie algebra
$\Omega^\bullet_X\otimes\overline C^\bullet(\mathcal{J}_X))[1]$. We
denote the action of the form $\omega \in
\Omega^\bullet_X\otimes\mathcal{J}_X/\mathcal{O}_X$ by
$\iota_{\omega}$. Consider now the subsheaf of closed forms
$(\Omega^{\bullet}_X\otimes \mathcal{J}_X/\mathcal{O}_X)^{cl}$ which
is by definition the kernel of $\nabla^{can}$.  $(\Omega^k_X\otimes
\mathcal{J}_X/\mathcal{O}_X)^{cl}$ acts by derivations of degree
$k-1$ and this action commutes with the differential
$\nabla^{can}+\delta$. Therefore, for $\omega\in
\Gamma(X;(\Omega^2\otimes\mathcal{J}_X/\mathcal{O}_X)^{cl})$ one can
define the $\omega$-twist $\mathfrak{g}(\mathcal{J}_X)_\omega$ as
the DGLA with the same underlying graded Lie algebra structure as
$\mathfrak{g}(\mathcal{J}_X)$ and the differential given by
$\nabla^{can}+\delta +\iota_{\omega}$. The isomorphism class of this
DGLA depends only on the cohomology class of $\omega$ in
$H^2(\Gamma(X;
\Omega^\bullet_X\otimes\mathcal{J}_X/\mathcal{O}_X),\nabla^{can})$.

More precisely, for
$\beta\in\Gamma(X;\Omega^1_X\otimes\mathcal{J}_X/\mathcal{O}_X)$ the
DGLA $\mathfrak{g}_\DR(\mathcal{J}_X)_\omega$ and
$\mathfrak{g}_\DR(\mathcal{J}_X)_{\omega+\nabla^{can}\beta}$ are
canonically isomorphic with the isomorphism depending only on the
equivalence class $\beta + \im(\nabla^{can})$.

As we remarked before a twisted form $\mathcal{S}$ of
$\mathcal{O}_X$ is determined  up to equivalence by its class in $
H^2(X;\mathcal{O}^\times)$. The composition $\mathcal{O}^\times\to
\mathcal{O}^\times/\mathbb{C}^\times\xrightarrow{\log}
\mathcal{O}/\mathbb{C}\xrightarrow{j^\infty} \DR(\mathcal{J/O})$
induces the map $H^2(X;\mathcal{O}^\times)\to
H^2(X;\DR(\mathcal{J/O}))\cong H^2(\Gamma(X;
\Omega^\bullet_X\otimes\mathcal{J}_X/\mathcal{O}_X),\nabla^{can})$.
We denote by $[\mathcal{S}]\in H^2(\Gamma(X;
\Omega^\bullet_X\otimes\mathcal{J}_X/\mathcal{O}_X),\nabla^{can})$
the image of the class of $\mathcal{S}$. By the remarks above we
have the well-defined up to a canonical isomorphism DGLA
$\mathfrak{g}_\DR(\mathcal{J}_X)_{[\mathcal{S}]}$.

The rest of this  paper is organized as follows. In the Section
\ref{prel} we review some preliminary facts. In the Section
\ref{dgla} we review the construction of Deligne 2-groupoid, its
relation with the deformation theory and its cosimplicial analogues.
In the Section \ref{algstack} we review the notion of algebroid
stacks. Next we define matrix algebras associated with a descent
datum in the Section \ref{dglaloc}. In the Section \ref{defalgstack}
we define the deformations of algebroid stacks and relate them to
the cosimplicial DGLA of Hochschild cochains on matrix algebras. In
the Section \ref{jets} we establish quasiisomorphism of the DGLA
controlling the deformations of twisted forms of $\co_X$ with a
simpler cosimplicial DGLA. Finally, the proof the main result of
this paper, Theorem \ref{main theorem}, is given in the Section
\ref{section: proof of main theorem}
\section{Preliminaries}\label{prel}

\subsection{Simplicial notions}

\subsubsection{The category of simplices}
For $n=0,1,2,\ldots$ we denote by $[n]$ the category with objects
$0,\ldots,n$ freely generated by the graph
\begin{equation*}
0\to 1\to\cdots\to n \ .
\end{equation*}
For $0\leq p\leq q\leq n$ we denote by $(pq)$ the unique morphism
$p\to q$.

We denote by $\Delta$ the full subcategory of $\Cat$ with objects
the categories $[n]$ for $n = 0,1,2,\ldots$.

For $0\leq i\leq n+1$ we denote by $\partial_i = \partial^n_i :
[n]\to[n+1]$ the $i^{\text{th}}$ face map, i.e. the unique map whose
image does not contain the object $i\in[n+1]$.

For $0\leq i\leq n-1$ we denote by $s_i=s^n_i:[n]\to[n-1]$ the
$i^{\text{th}}$ degeneracy map, i.e. the unique surjective map such that $s_i(i)=s_i(i+1)$.

\subsubsection{Simplicial and cosimplicial objects}
Suppose that $\mathcal C$ is a category. By definition, \emph{a
simplicial object in $\mathcal C$} (respectively, \emph{a
cosimplicial object in $\mathcal C$}) is a functor
$\Delta^{op}\to\mathcal{C}$ (respectively, a functor
$\Delta\to\mathcal{C}$). Morphisms of (co)simplicial objects are
natural transformations of functors.

For a simplicial (respectively, cosimplicial) object $F$ we denote
the object $F([n])\in\mathcal{C}$ by $F_n$ (respectively, $F^n$).

\subsection{Cosimplicial vector spaces}

Let $V^{\bullet}$ be a cosimplicial vector space.   We denote by
$C^{\bullet}(V)$ the associated complex with component $C^n(V)=V^n$
and the differential $\partial^n : C^n(V)\to C^{n+1}(V)$ defined by
$\partial^n = \sum_i (-1)^i
\partial^n_i$, where $\partial^n_i$ is the map induced by the
$i^{\text{th}}$ face map $[n]\to[n+1]$. We denote cohomology of this
complex by $H^{\bullet}(V)$.

The complex $C^{\bullet}(V)$ contains a normalized subcomplex
$\overline{C}^{\bullet}(V)$. Here $\overline{C}^{n}(V)=\{V \in V^n \
|  \ s^n_i v=0\}$, where $s^n_i:[n] \to [n-1]$ is the
$i^{\text{th}}$ degeneracy map. Recall that the inclusion
$\overline{C}^{\bullet}(V) \to C^{\bullet}(V)$ is a
quasiisomorphism.

Starting from a cosimplicial vector space $V^{\bullet}$ one can
construct a new cosimplicial vector space $\widehat{V}^{\bullet}$ as
follows. For every $\lambda:[n] \to \Delta$ set
$\widehat{V}^{\lambda} =V^{\lambda(n)}$. Suppose given another
simplex $\mu : [m]\to\Delta$ and morphism $\phi : [m]\to[n]$ such
that $\mu = \lambda\circ\phi$ i.e. $\phi$ is a morphism of simplices
$\mu\to\lambda$. The morphism $(0n)$ factors uniquely into
$0\to\phi(0)\to\phi(m)\to n$, which, under $\lambda$, gives the
factorization of $\lambda(0n):\lambda(0)\to\lambda(n)$ (in $\Delta$)
into
\begin{equation}
\begin{CD}
\lambda(0) @>f>> \mu(0) @>g>> \mu(m) @>h>> \lambda(n) \ ,
\end{CD}
\end{equation}
where $g=\mu(0m)$. The map $\mu(m)\to\lambda(n)$ induces the map
\begin{equation}\label{cosim}
\phi_* :  \widehat{V}^{\mu}\to \widehat{V}^{\lambda}
\end{equation}
Set now   $\widehat{V}^n = \prod_{[n]\xrightarrow{\lambda}\Delta}
\widehat{V}^{\lambda}$.  The maps \eqref{cosim} endow
$\widehat{V}^{\bullet}$ with the structure of a cosimplicial vector
space.   We then have the following well-known result:
\begin{lemma}\label{subdivision}
\begin{equation}
H^{\bullet}(V)\cong H^{\bullet}(\widehat{V})
\end{equation}
\end{lemma}
\begin{proof}
We construct morphisms of complexes inducing the isomorphisms in
cohomology.  We will use the following notations. If $f \in
\widehat{V}^{n}$ and $\lambda: [n] \to \Delta$ we will denote by
$f(\lambda) \in \widehat{V}^{\lambda}$ its component in
$\widehat{V}^{\lambda}$. For $\lambda:[n] \to \Delta$ we denote by
$\lambda|_{[jl]}: [l-j] \to \Delta$ its truncation:
$\lambda|_{[jl]}(i)=\lambda (i+j)$, $\lambda|_{[jl]}(i, k)=\lambda
((i+j) \ (k+j))$. For $\lambda_1: [n_1]\to \Delta$ and $\lambda_2: [n_2]
\to \Delta$ with $\lambda_1(n_1)=\lambda_2(0)$ define their
concatenation $\Lambda=\lambda_1 \ast \lambda_2 :[n_1+n_2]\to
\Delta$ by the following formulas.
\begin{equation*}
\Lambda(i)= \begin{cases}\lambda_1(i) \text{ if } i \le n_1\\
\lambda_2(i-n_1)  \text{ if } i \ge n_1\end{cases}
\end{equation*}
\begin{equation*}
\Lambda(ik)=\begin{cases}\lambda_1(ik) \text{ if } i, k \le n_1\\
\lambda_2((i-n_1)\ (k-n_1))  \text{ if } i, k \ge n_1\\
\lambda_2(0 \ (k-n_1)) \circ \lambda (i n_1) \text{ if } i \le n_1
\le k
\end{cases}
\end{equation*}
This operation is associative. Finally we will identify in our
notations $\lambda:[1] \to \Delta$ with the morphism $\lambda(01)$
in $\Delta$

The morphism $C^{\bullet}(V) \to C^{\bullet}(\widehat{V})$ is
constructed as follows. Let $\lambda: [n] \to \Delta $ be a simplex
in $\Delta$ and define $\lambda_k$ by  $\lambda(k) = [\lambda_k]$,
$k=0, 1,\ldots, n$. Let $\Upsilon(\lambda): [n] \to \lambda(n)$ be a
morphism in $\Delta$ defined by
\begin{equation}
(\Upsilon(\lambda))(k)=\lambda(kn)(\lambda_k)
\end{equation}
Then define the map $\iota:V^{\bullet} \to \widehat{V}^{\bullet}$ by
the formula
\begin{equation*}
(\iota(v))(\lambda)=\Upsilon(\lambda)_*v \text{ for } v \in V^n
\end{equation*}
This is a map of cosimplicial vector spaces, and therefore it
induces a morphism of complexes.

The morphism $\pi: C^{\bullet}(\widehat{V}) \to C^{\bullet}(V)  $ is
defined by the formula
\begin{equation*}
\pi(f) = (-1)^{\frac{n(n+1)}{2}}\sum \limits_{0\le i_k \le k+1}
(-1)^{i_0+\dots+ i_{n-1}}f(
\partial_{i_0}^0\ast
\partial_{i_1}^1\ast \ldots \ast \partial_{i_{n-1}}^{n-1}) \text{
for } f \in \widehat{V}^{n}
\end{equation*}
when $n>0$, and $\pi(f)$ is $V^0$ component of $f$ if $n=0$.

The morphism $\iota \circ \pi$ is homotopic to $\id$ with the
homotopy   $h: C^{\bullet}(\widehat{V}) \to
C^{\bullet-1}(\widehat{V})$ given by the formula
\begin{equation*}
hf(\lambda) = \sum \limits_{j=0}^{n-1}  \ \sum \limits_{0\le i_k \le
k+1}(-1)^{i_0+\dots+ i_{j-1}} f(\partial_{i_0}^0\ast \dots \ast
\partial_{i_{j-1}}^{j-1}\ast\Upsilon(\lambda|_{[0j]})\ast
\lambda|_{[j\ (n-1)]})
\end{equation*}
$\text{ for } f \in \widehat{V}^{n}$ when $n>0$, and $h(f)=0$ if
$n=0$.

The composition $\pi \circ \iota: C^{\bullet}(V) \to C^{\bullet}(V)$
preserves the normalized subcomplex $\overline{C}^{\bullet}(V)$ and
acts by identity on it. Therefore $\pi \circ \iota$ induces the
identity map on cohomology. It follows that $\pi$ and $\iota$ are
quasiisomorphisms inverse to each other.
\end{proof}

\subsection{Covers}

A cover (open cover) of a space $X$ is a collection $\mathcal{U}$ of
open subsets of $X$ such that $\bigcup_{U\in\mathcal{U}}U = X$.

\subsubsection{The nerve of a cover}\label{nerve of cover}
Let $N_0\mathcal{U} = \coprod_{U\in\mathcal{U}}U$. There is a
canonical augmentation map
\begin{equation*}
\epsilon_0: N_0\mathcal{U}\xrightarrow{\coprod_{U\in\mathcal{U}}(U\hookrightarrow
X)} X
\end{equation*}
Let
\begin{equation*}
N_p\mathcal{U} = N_0\mathcal{U}\times_X\cdots\times_X N_0\mathcal{U}
\ ,
\end{equation*}
the $(p+1)$-fold fiber product.

The assignment $N\mathcal{U}: \Delta\ni [p]\mapsto N_p\mathcal{U}$
extends to a simplicial space called \emph{the nerve of the cover
$\mathcal{U}$}. The effect of the face map $\partial^n_i$
(respectively, the degeneracy map $s^n_i$) will be denoted by $d^i =
d_n^i$ (respectively, $\varsigma_i = \varsigma^n_i$) and is given by
the projection \emph{along} the $i^{\text{th}}$ factor
(respectively, the diagonal embedding on the $i^{\text{th}}$
factor). Therefore for every morphism $f:[p]\to [q]$ in $\Delta$ we
have a morphism $N_q\cu \to N_p\cu$ which we denote by $f^*$. We
will denote by $f_*$ the operation $(f^*)^*$ of pull-back along
$f^*$; if $\cf$ is a sheaf on $N_p\cu$  then $f_* \cf$ is a sheaf on
$N_q \cu$.

For $0\leq i\leq n+1$ we denote by $\pr_i^n : N_n\cu\to N_0\cu$ the
projection onto the $i^{\text{th}}$ factor. For $0\leq j\leq m$,
$0\leq i_j\leq n$  the map  $\pr_{i_0}\times\cdots\times\pr_{i_m}:
N_n \cu \to (N_0\cu)^m$ can be factored uniquely as a composition of
a map $ N_n \cu \to N_m\cu$ and the canonical imbedding $N_m\cu \to
(N_0\cu)^m$.
 We denote this map $ N_n \cu \to N_m\cu$ by
$\pr_{i_0\ldots i_m}^n$.

The augmentation map $\epsilon_0$ extends to a morphism  $\epsilon : N\mathcal{U}\to
X$ where the latter is regarded as a constant simplicial space. Its component
of degree $n$ $\epsilon_n: N_n \cu \to X$  is given by the formula
$\epsilon_n = \epsilon_0 \circ \pr_i^n$. Here $0\leq i\leq n+1$ is arbitrary.


\subsubsection{\v Cech complex}
Let $\cf$ be a sheaf of abelian groups on $X$. One defines a
cosimplicial group $\check{C}^{\bullet} (\cu, \cf) = \Gamma
(N_{\bullet}\cu; \epsilon^* \cf)$, with the cosimplicial structure
induced by the simplicial structure of $N\cu$. The associated
complex is the \v Cech complex of the cover $\cu$ with coefficients
in  $\cf$. The differential $\check{\partial}$ in  this complex is
given by $\sum (-1)^i (d^i)^*$.

\subsubsection{Refinement}
Suppose that $\mathcal{U}$ and $\mathcal{V}$ are two covers of $X$.
A morphism of covers $\rho : \mathcal{U}\to\mathcal{V}$ is a map of
sets $\rho : \mathcal{U}\to\mathcal{V}$ with the property
$U\subseteq \rho(U)$ for all $U\in\mathcal{U}$.

A morphism $\rho : \mathcal{U}\to\mathcal{V}$ induces the map $N\rho
: N\mathcal{U}\to N\mathcal{V}$ of simplicial spaces which commutes
with respective augmentations to $X$. The map $N_0\rho$ is determined by the commutativity of
\begin{equation*}
\begin{CD}
U @>>> N_0\mathcal{U} \\
@VVV @VV{N_0\rho}V \\
\rho(U) @>>> N_0\mathcal{V}
\end{CD}
\end{equation*}
It is clear that the map $N_0\rho$ commutes with the respective augmentations (i.e. is a map of spaces over $X$) and, consequently induces maps $N_n\rho = (N_0\rho)^{\times_X n+1}$ which commute with all structure maps.

\subsubsection{The category of covers}
Let $\Cov(X)_0$ denote the set of open covers of X. For
$\mathcal{U},\mathcal{V}\in\Cov(X)_0$ we denote by
$\Cov(X)_1(\mathcal{U},\mathcal{V})$ the set of morphisms
$\mathcal{U}\to\mathcal{V}$. Let $\Cov(X)$ denote the category
with objects $\Cov(X)_0$ and morphisms $\Cov(X)_1$. The
construction of \ref{nerve of cover} is a functor
\[
N : \Cov(X)\to\Top^{\Delta^{op}}/X \ .
\]

\section{Deligne $2$-groupoid and its cosimplicial analogues}\label{dgla}
In this section we begin by recalling the definition of Deligne
$2$-groupoid and its relation with the deformation theory. We then
describe the cosimplicial analogues of Deligne $2$-groupoids and
establish some of their properties.
\subsection{Deligne $2$-groupoid}
In this subsection we review the construction of Deligne
$2$-groupoid of a nilpotent differential graded algebra (DGLA). We
follow \cite{G, G1} and references therein.

Suppose that $\mathfrak{g}$ is a nilpotent DGLA such that
$\mathfrak{g}^i = 0$ for $i< -1$.

A Maurer-Cartan element of $\mathfrak{g}$  is an element $\gamma \in
\mathfrak{g}^1$ satisfying
\begin{equation} \label{eq:MC}
d\gamma + \frac{1}{2}[\gamma,\gamma]=0.
\end{equation}
We denote by $\MC^2(\mathfrak{g})_0$ the set of Maurer-Cartan
elements of $\mathfrak{g}$.

 The unipotent group $\exp
\mathfrak{g}^0$ acts on the set of Maurer-Cartan elements of
$\mathfrak{g}$ by the gauge equivalences. This action is given by
the formula
\begin{equation*}
(\exp X) \cdot \gamma= \gamma- \sum_{i=0}^{\infty} \frac{(\ad
X)^i}{(i+1)!}(dX +[\gamma, X])
\end{equation*}
If $\exp X$  is a gauge equivalence between two Maurer-Cartan
elements $\gamma_1 $ and $\gamma_2= (\exp X)\cdot \gamma_1$ then
\begin{equation} \label{eq:MC equivalence}
d+\ad \gamma_2=  \Ad \exp X\,(d+\ad \gamma_1).
\end{equation}
We denote by $\MC^2(\mathfrak{g})_1(\gamma_1, \gamma_2)$ the set of
gauge equivalences between $\gamma_1$, $\gamma_2$. The composition
\begin{equation*}
\MC^2(\mathfrak{g})_1(\gamma_2,
\gamma_3)\times\MC^2(\mathfrak{g})_1(\gamma_1,
\gamma_2)\to\MC^2(\mathfrak{g})_1(\gamma_1, \gamma_3)
\end{equation*}
is given by  the product in the group $\exp \mathfrak{g}^0$.

If $\gamma \in \MC^2(\mathfrak{g})_0$  we can define a Lie bracket
$[\cdot, \cdot]_{\gamma}$ on $\mathfrak{g}^{-1}$ by
\begin{equation}\label{eq:mu-bracket}
[a,\,b]_{\gamma}=[a,\, d b+[\gamma, \,b]].
\end{equation}
With this bracket $\mathfrak{g}^{-1}$ becomes a nilpotent Lie
algebra. We denote by $\exp_{\gamma} \mathfrak{g}^{-1}$ the
corresponding unipotent group, and by $\exp_{\gamma}$ the
corresponding exponential map $\mathfrak{g}^{-1} \to \exp_{\gamma}
\mathfrak{g}^{-1}$. If $\gamma_1$, $\gamma_2$ are two Maurer-Cartan
elements, then the group $\exp_{\gamma} \mathfrak{g}^{-1}$ acts on
$\MC^2(\mathfrak{g})_1(\gamma_1, \gamma_2)$. Let $\exp_{\gamma} t
\in \exp_{\gamma} \mathfrak{g}^{-1}$ and let $\exp X \in
\MC^2(\mathfrak{g})_1(\gamma_1, \gamma_2)$. Then
\begin{equation*}
(\exp_{\gamma} t) \cdot (\exp X) = {\operatorname
{exp}}(dt+[\gamma,t])\, {\operatorname {exp}}X \in \exp
\mathfrak{g}^0
\end{equation*}
Such an element $\exp_{\gamma} t$ is called a $2$-morphism between
$\exp X$ and $(\exp t) \cdot (\exp X)$. We denote by
$\MC^2(\mathfrak{g})_2(\exp X,\exp Y)$ the set of $2$-morphisms
between $\exp X$ and $\exp Y$. This set is endowed with a vertical
composition given by the product in the group $\exp_{\gamma}
\mathfrak{g}^{-1}$.

 Let $\gamma_1$, $\gamma_2$, $\gamma_3 \in
\MC^2(\mathfrak{g})_0$. Let $\exp X_{12}$, $\exp Y_{12}\in
\MC^2(\mathfrak{g})_1(\gamma_1, \gamma_2)$ and $\exp X_{23}$, $\exp
Y_{23}\in \MC^2(\mathfrak{g})_1(\gamma_2, \gamma_3)$. Then one
defines the horizontal composition
\begin{multline*}
\otimes: \MC^2(\mathfrak{g})_2(\exp X_{23}, \exp Y_{23}) \times
\MC^2(\mathfrak{g})_2(\exp X_{12}, \exp Y_{12}) \to\\
\MC^2(\mathfrak{g})_2(\exp X_{23}\exp X_{12}, \exp X_{23}\exp
Y_{12})
\end{multline*}
as follows. Let $\exp_{\gamma_2} t_{12} \in
\MC^2(\mathfrak{g})_2(\exp X_{12}, \exp Y_{12})$, $\exp_{\gamma_3}
t_{23} \in \MC^2(\mathfrak{g})_2(\exp X_{23}, \exp Y_{23})$. Then
\begin{equation*}
\exp_{\gamma_{3}} t_{23} \otimes \exp_{\gamma_{2}} t_{12}
=\exp_{\gamma_{3}} t_{23}\exp_{\gamma_3}( e^{\ad X_{23}}(t_{12}) )
\end{equation*}

To summarize, the data described above forms a $2$-groupoid which we
denote by $\MC^2(\mathfrak{g})$ as follows:
\begin{enumerate}
\item
the set of objects is $\MC^2(\mathfrak{g})_0$

\item the groupoid of morphisms
$\MC^2(\mathfrak{g})(\gamma_1,\gamma_2)$,
$\gamma_i\in\MC^2(\mathfrak{g})_0$ consists of:
\begin{itemize}
\item objects i.e. $1$-morphisms in $\MC^2(\mathfrak{g})$ are given
by $\MC^2(\mathfrak{g})_1(\gamma_1,\gamma_2) $  -- the gauge
transformations between $\gamma_1$ and $\gamma_2$.
\item morphisms between
$\exp X$, $\exp Y \in \MC^2(\mathfrak{g})_1(\gamma_1,\gamma_2)$ are
given by $\MC^2(\mathfrak{g})_2(\exp X,\exp Y)$.
\end{itemize}
\end{enumerate}

 A morphism of nilpotent DGLA $ \phi :
\mathfrak{g} \to \mathfrak{h}$ induces a functor $\phi:
\MC^2(\mathfrak{g}) \to \MC^2(\mathfrak{g})$.

We have the following important result (\cite{GM}, \cite{G} and
references therein).
\begin{thm}\label{thm: quism invariance of mc}
Suppose that $\phi : \mathfrak{g}\to\mathfrak{h}$ is a
quasi-isomorphism of DGLA and  let $\mathfrak{m}$ be a nilpotent
commutative ring. Then the induced map $\phi :
\MC^2(\mathfrak{g}\otimes \mathfrak{m})\to\MC^2(\mathfrak{h}\otimes
\mathfrak{m})$ is an equivalence of $2$-groupoids.
\end{thm}

\subsection{Deformations and Deligne $2$-groupoid}\label{ddgla}
Let $k$ be an algebraically closed field of characteristic zero.

\subsubsection{Hochschild cochains}
Suppose that $A$ is a $k$-vector space. The $k$-vector space $C^n(A)$ of Hochschild cochains of degree $n\ge 0$ is defined by
\begin{equation*}
C^n(A) := \Hom_{k}(A^{\otimes n}, A) \ .
\end{equation*}
The graded vector space $\mathfrak{g}(A) := C^{\bullet}(A)[1]$ has a canonical structure of a graded Lie algebra under the Gerstenhaber bracket denoted by $[\ ,\ ]$ below. Namely, $C^{\bullet}(A)[1]$ is canonically isomorphic to the (graded) Lie algebra of derivations of the free associative co-algebra generated by $A[1]$.

Suppose in addition that $A$ is equipped with a bilinear operation $\mu : A\otimes A\to A$, i.e. $\mu\in C^2(A)= \mathfrak{g}^1(A)$. The condition $[\mu,\mu] = 0$ is equivalent to the associativity of $\mu$.

Suppose that $A$ is an associative $k$-algebra with the product $\mu$. For $a\in \mathfrak{g}(A)$ let $\delta(a) = [\mu, a]$. Thus, $\delta$ is a derivation of the graded Lie algebra $\mathfrak{g}(A)$. The associativity of $\mu$ implies that $\delta^2 = 0$, i.e. $\delta$ defines a differential on $\mathfrak{g}(A)$ called the Hochschild differential.

For a unital algebra the subspace of \emph{normalized cochains} $\overline{C}^{n}(A) \subset C^n(A)$ is defined by
\[
\overline{C}^{n}(A) :=
\Hom_{k}((A/k\cdot 1)^{\otimes n}, A) \ .
\]
The subspace $\overline{C}^{\bullet}(A)[1]$ is closed under the Gerstenhaber bracket
and the action of the Hochschild differential and the inclusion
$\overline{C}^{\bullet}(A)[1]\hookrightarrow C^\bullet(A)[1]$ is a quasi-isomorphism
of DGLA.

Suppose in addition that $R$ is a commutative Artin $k$-algebra with
the nilpotent maximal ideal $\mathfrak{m}_R$ The DGLA
$\mathfrak{g}(A)\otimes_k\mathfrak{m}_R$ is nilpotent and satisfies
$\mathfrak{g}^i(A)\otimes_k\mathfrak{m}_R = 0$ for $i < -1$.
Therefore, the Deligne $2$-groupoid
$\MC^2(\mathfrak{g}(A)\otimes_k\mathfrak{m}_R)$ is defined.
Moreover, it is clear that the assignment $R\mapsto \MC^2(
\mathfrak{g}(A)\otimes_k\mathfrak{m}_R)$ extends to a functor on the
category of commutative Artin algebras.

\subsubsection{Star products}
Suppose that $A$ is an associative unital $k$-algebra. Let $m$
denote the product on $A$.

Let $R$ be a commutative Artin $k$-algebra with maximal ideal
$\mathfrak{m}_R$. There is a canonical isomorphism
$R/\mathfrak{m}_R\cong k$.

A \emph{($R$-)star product} on $A$ is an associative $R$-bilinear product on
$A\otimes_k R$ such that the canonical isomorphism of $k$-vector spaces
$(A\otimes_k R)\otimes_R k\cong A$ is an isomorphism of algebras. Thus, a star
product is an $R$-deformation of $A$.

The $2$-category of $R$-star products on $A$, denoted $\Def(A)(R)$,
is defined as the subcategory of the $2$-category $\Alg^2_R$ of
$R$-algebras (see \ref{subsubsection: algebroids}) with
\begin{itemize}
\item Objects: $R$-star products on $A$,

\item $1$-morphisms $\phi : m_1 \to m_2$ between the star products
$\mu_i$ those $R$-algebra homomorphisms $\phi : (A\otimes_k R,
m_1)\to (A\otimes_k R, m_2)$ which reduce to the identity map modulo
$\mathfrak{m}_R$, i.e. $\phi\otimes_R k = \id_A$

\item $2$-morphisms $b : \phi \to \psi$, where $\phi, \psi : m_1
\to m_2$ are two $1$-morphisms, are elements $b \in 1+A \otimes_k
\mathfrak{m}_R \subset A \otimes_k R$ such that $m_2(\phi(a),b) =
m_2(b,\psi(a))$ for all $a\in A\otimes_k R$.
\end{itemize}
It follows easily from the above definition and the nilpotency of
$\mathfrak{m}_R$ that $\Def(A)(R)$ is a $2$-groupoid.

Note that $\Def(A)(R)$ is non-empty: it contains the trivial
deformation, i.e. the star product, still denoted $m$, which is the
$R$-bilinear extension of the product on $A$.

It is clear that the assignment $R\mapsto \Def(A)(R)$ extends to a
functor on the category of commutative Artin $k$-algebras.

\subsubsection{Star products and the Deligne $2$-groupoid}
We continue in notations introduced above. In particular, we are
considering an associative unital $k$-algebra $A$. The product $m\in
C^2(A)$ determines a cochain, still denoted $m\in
\mathfrak{g}^1(A)\otimes_k R$, hence the Hochschild differential
$\delta = [m,\ ]$ in $\mathfrak{g}(A)\otimes_k R$ for any
commutative Artin $k$-algebra $R$.

Suppose that $m'$ is an $R$-star product on $A$. Since $\mu(m') :=
m' - m = 0 \mod\mathfrak{m}_R$ we have
$\mu(m')\in\mathfrak{g}^1(A)\otimes_k\mathfrak{m}_R$. Moreover, the
associativity of $m'$ implies that $\mu(m')$ satisfies the
Maurer-Cartan equation, i.e.
$\mu(m')\in\MC^2(\mathfrak{g}(A)\otimes_k \mathfrak{m}_R)_0$.

It is easy to see that the assignment $m'\mapsto\mu(m')$ extends to
a functor
\begin{equation}\label{functor Def to MC}
\Def(A)(R)\to \MC^2(\mathfrak{g}(A)\otimes_k \mathfrak{m}_R) \ .
\end{equation}

The following proposition is well-known (cf. \cite{Ge, G, G1}).

\begin{prop}\label{Def equivalent to MC}
The functor \eqref{functor Def to MC} is an isomorphism of
$2$-groupoids.
\end{prop}

\subsubsection{Star products on sheaves of algebras}
\label{subsubsection: star product on sheaves}
The above
considerations generalize to sheaves of algebras in a
straightforward way.

Suppose that $\mathcal A$ is a sheaf of $k$-algebras on a space $X$.
Let $m : \mathcal{A}\otimes_k\mathcal{A}\to \mathcal{A}$ denote the
product.

An $R$-star product on $\mathcal A$ is a structure of a sheaf of an
associative algebras on $\mathcal{A}\otimes_k R$ which reduces to
$\mu$ modulo the maximal ideal $\mathfrak{m}_R$. The $2$-category
(groupoid) of $R$ star products on $\mathcal A$, denoted
$\Def(\mathcal{A})(R)$ is defined just as in the case of algebras;
we leave the details to the reader.

The sheaf of Hochschild cochains of degree $n$ is defined by
\begin{equation*}
C^n({\mathcal A}) := \shHom({\mathcal A}^{\otimes n},{\mathcal A}).
\end{equation*}
We have the sheaf of DGLA $\mathfrak{g}(A) :=
C^\bullet(\mathcal{A})[1]$, hence the nilpotent DGLA
$\Gamma(X;\mathfrak{g}(\mathcal{A})\otimes_k \mathfrak{m}_R)$ for
every commutative Artin $k$-algebra $R$ concentrated in degrees
$\geq -1$. Therefore, the $2$-groupoid
$\MC^2(\Gamma(X;\mathfrak{g}(\mathcal{A})\otimes_k \mathfrak{m}_R)$
is defined.

The canonical functor $\Def(\mathcal{A})(R)\to
\MC^2(\Gamma(X;\mathfrak{g}(\mathcal{A})\otimes_k \mathfrak{m}_R)$
defined just as in the case of algebras is an isomorphism of
$2$-groupoids.

\subsection{$\mathfrak{G}$-stacks}
Suppose that $\mathfrak{G}$: $[n]\to \mathfrak{G}^n$  is a
cosimplicial DGLA. We assume that each $\mathfrak{G}^n$ is a
nilpotent DGLA. We denote its component of degree $i$ by
$\mathfrak{G}^{n, i}$ and assume that $\mathfrak{G}^{n, i}=0$ for
$i<-1$.

\begin{definition}
A $\mathfrak G$-stack is a triple $\gamma =
(\gamma^0,\gamma^1,\gamma^2)$, where
\begin{itemize}
\item
$\gamma^0\in\MC^2(\mathfrak{G}^0)_0$,

\item $\gamma^1\in
\MC^2(\mathfrak{G}^1)_1(\partial^0_0\gamma^0,\partial^0_1\gamma^0)$,

satisfying the condition
\begin{equation*}
s_0^1 \gamma^1 =\id
\end{equation*}

\item $\gamma^2\in
\MC^2(\mathfrak{G}^2)_2(\partial^1_2(\gamma^1)\circ\partial^1_0(\gamma^1),
\partial^1_1(\gamma^1))$
\end{itemize}
satisfying the conditions
\begin{equation}\label{gamma2}\begin{split}
\partial_2^2 \gamma_2 \circ (\id \otimes \partial_0^2 \gamma_2) &=
\partial_1^2 \gamma_2 \circ (\partial_3^2 \gamma_2 \otimes \id)\\
s_0^2 \gamma^2 &= s_1^2 \gamma^2 =\id
\end{split}\end{equation}
\end{definition}
Let $\Stack(\mathfrak{G})_0$ denote the set of
$\mathfrak{G}$-stacks.

\begin{definition}
For $\gamma_1,\gamma_2\in\Stack(\mathfrak{G})_0$ a $1$-morphism
$\gimel : \gamma_1\to\gamma_2$ is a pair $\gimel =
(\gimel^1,\gimel^2)$, where
$\gimel^1\in\MC^2(\mathfrak{G}^0)_1(\gamma_1^0,\gamma_2^0)$,
$\gimel^2\in\MC^2(\mathfrak{G}^1)_2(\gamma_2^1\circ\partial^0_0(\gimel^1),
\partial^0_1(\gimel^1)\circ\gamma_1^1)$, satisfying
\begin{equation}\label{1-morphism}
\begin{split}
(\id \otimes \gamma_1^2) \circ (\partial_2^1 \gimel^2 \otimes \id)
\circ (\id \otimes \partial_0^1 \gimel^2) &=
\partial_1^1\gimel^2 \circ(\gamma_2^2 \otimes \id)\\
s_0^1 \gimel^2 =\id
\end{split}
\end{equation}
\end{definition}
Let $\Stack(\mathfrak{G})_1(\gamma_1,\gamma_2)$ denote the set of
$1$-morphisms $\gamma_1\to\gamma_2$

Composition of $1$-morphisms $\gimel: \gamma_1 \to \gamma_2$ and
$\daleth: \gamma_2 \to \gamma_3$ is given by $(\daleth^1 \circ
\gimel^1, (\gimel^2\otimes \id) \circ (\id \otimes \daleth^2))$.

\begin{definition}
For $\gimel_1,\gimel_2\in\Stack(\mathfrak{G})_1(\gamma_1,\gamma_2)$
a $2$-morphism $\phi:\gimel_1\to\gimel_2$ is a $2$-morphism
$\phi\in\MC^2(\mathfrak{G}^0)_2(\gimel_1^1,\gimel_2^1)$ which
satisfies
\begin{equation}\label{2-morphism}
\gimel_2^2\circ (\id \otimes \partial_0^0 \phi  )  = (
\partial_1^0 \phi \otimes \id ) \circ \gimel_1^2
\end{equation}
\end{definition}
Let $\Stack(\mathfrak{G})_2(\gimel_1,\gimel_2)$ denote the set of
$2$-morphisms.

Compositions of $2$-morphisms are given by the compositions in
$\MC^2(\mathfrak{G}^0)_2$.

For $\gamma_1,\gamma_2\in\Stack(\mathfrak{G})_0$, we have the
groupoid $\Stack(\mathfrak{G})(\gamma_1,\gamma_2)$ with the set of
objects  $\Stack(\mathfrak{G})_1(\gamma_1,\gamma_2)$ and the set of
morphisms $\Stack(\mathfrak{G})_2(\gimel_1,\gimel_2)$ under vertical
composition.


Every morphism $\theta$ of cosimplicial DGLA induces in an obvious
manner a functor $\theta_*:
\Stack(\mathfrak{G}\otimes\mathfrak{m})\to\Stack(\mathfrak{H}\otimes\mathfrak{m})$

We have the following cosimplicial analogue of the Theorem \ref{thm:
quism invariance of mc}:
\begin{thm}\label{quism invariance of stack}
Suppose that $\theta : \mathfrak{G}\to\mathfrak{H}$ is a
quasi-isomorphism of cosimplicial DGLA and $\mathfrak{m}$ is a
commutative nilpotent ring. Then the induced map $\theta_* :
\Stack(\mathfrak{G}\otimes\mathfrak{m})\to\Stack(\mathfrak{H}\otimes\mathfrak{m})$
is an equivalence.
\end{thm}
\begin{proof}
The proof can be obtained by applying the Theorem \ref{thm: quism
invariance of mc} repeatedly.

Let $\gamma_1$, $\gamma_2 \in
\Stack(\mathfrak{G}\otimes\mathfrak{m})_0$, and let $\gimel_1$,
$\gimel_2$ be two $1$-morphisms between $\gamma_1$ and $\gamma_2$.
Note that by the Theorem \ref{thm: quism invariance of mc}
$\theta_*: \MC^2(\mathfrak{G}^0\otimes\mathfrak{m})_2(\gimel_1^1,
\gimel_2^1) \to
\MC^2(\mathfrak{H}^0\otimes\mathfrak{m})_2(\theta_*\gimel_1^1,
\theta_*\gimel_2^1)$ is a bijection. Injectivity of the map $
\theta_*: \Stack(\mathfrak{G}\otimes\mathfrak{m})_2(\gimel_1,
\gimel_2) \to
\Stack(\mathfrak{H}\otimes\mathfrak{m})_2(\theta_*\gimel_1,
\theta_*\gimel_2)$ follows immediately. For the surjectivity, notice
that an element of
$\Stack(\mathfrak{H}\otimes\mathfrak{m})_2(\theta_*\gimel_1,
\theta_*\gimel_2)$ is necessarily given by $\theta_* \phi$ for some
$\phi \in \MC^2 (\mathfrak{G}^0\otimes\mathfrak{m})(\gimel_1^1,
\gimel_2^1)_2$ and the following identity is satisfied in
$\MC^2(\mathfrak{H}^1\otimes\mathfrak{m})(\theta_*(\gamma_2^1 \circ
\partial_0^0\gimel_1^1),\theta_*(\partial_1^0\gimel_2^1\circ \gamma_1^1))_2$:
$\theta_*(\gimel_2^2\circ (\partial_0^1 \phi \otimes \id)) =
\theta_*((\id \otimes
\partial_1^1 \phi) \circ \gimel_1^2)$.  Since $\theta_*:
\MC^2(\mathfrak{G}^1\otimes\mathfrak{m})(\gamma_2^1 \circ
\partial_0^0\gimel_1^1, \partial_1^0\gimel_2^1\circ \gamma_1^1)_2 \to
\MC^2(\mathfrak{H}^1\otimes\mathfrak{m})(\theta_*(\gamma_2^1 \circ
\partial_0^0\gimel_1^1),\theta_*(\partial_1^0\gimel_2^1\circ \gamma_1^1))_2$ is bijective,
and in particular injective, $  \gimel_2^2\circ (\partial_0^1 \phi
\otimes \id) =   (\id \otimes
\partial_1^1 \phi) \circ \gimel_1^2$, and $\phi$ defines an element
in $\Stack(\mathfrak{G}\otimes\mathfrak{m})_2(\gimel_1, \gimel_2)$.

Next, let  $\gamma_1$, $\gamma_2 \in
\Stack(\mathfrak{G}\otimes\mathfrak{m})_0$, and let $\gimel$ be a
$1$-morphisms between $\theta_*\gamma_1$ and $\theta_*\gamma_2$. We
show that there exists $\daleth \in
\Stack(\mathfrak{G}\otimes\mathfrak{m})_1(\gamma_1, \gamma_2)$ such
that $\theta_* \daleth$ is isomorphic to $\gimel$. Indeed, by
Theorem \ref{thm: quism invariance of mc}, there exists $\daleth^1
\in \MC^2(\mathfrak{G^0}\otimes\mathfrak{m})_1(\gamma_1, \gamma_2)$
such that $\MC^2(\mathfrak{H^0}\otimes\mathfrak{m})_2(\theta_*
\daleth^1, \gimel) \ne \emptyset$. Let $\phi \in
\MC^2(\mathfrak{H^0}\otimes\mathfrak{m})_2(\theta_* \daleth^1,
\gimel)$. Define $\psi \in
\MC^2(\mathfrak{H}^1\otimes\mathfrak{m})_2(\theta_*(\gamma_2^1 \circ
\partial_0^0 \daleth^1), \theta_*(\partial_1^0 \daleth^1 \circ \gamma_1^1))$
 by $\psi=
(\partial_1^0 \phi \otimes \id)^{-1} \circ \gimel^2 \circ (\id
\otimes\partial_0^0 \phi)$. It is easy to verify that the following
identities holds:
\begin{equation*}
\begin{split}
(\id \otimes \theta_*\gamma_1^2) \circ (\partial_2^1 \psi \otimes
\id) \circ (\id \otimes \partial_0^1 \psi) &=
\partial_1^1\psi \circ(\theta_*\gamma_2^2 \otimes \id)
\\
s_0^1\psi= \id
\end{split}
\end{equation*}
By bijectivity of $\theta_*$ on $\MC^2$ there exists a unique
$\daleth^2 \in \MC^2(\mathfrak{H}^1\otimes\mathfrak{m})_2(
\gamma_2^1 \circ
\partial_0^0 \daleth^1,  \partial_1^0 \daleth^1 \circ \gamma_1^1)$
such that $\theta_* \daleth^2 = \psi$. Moreover, as before,
injectivity of $\theta_*$ implies that the conditions
\eqref{1-morphism} are satisfied. Therefore $\daleth=(\daleth^1,
\daleth^2)$ defines a $1$-morphism $\gamma_1 \to \gamma_2$ and
$\phi$ is a $2$-morphism $\theta_*\daleth \to \gimel$.

Now, let $\gamma \in \Stack(\mathfrak{H}\otimes\mathfrak{m})_0$. We
construct $\lambda \in \Stack(\mathfrak{G}\otimes\mathfrak{m})_0$
such that
$\Stack(\mathfrak{H}\otimes\mathfrak{m})_1(\theta_*\lambda, \gamma)
\ne \varnothing$. Indeed, by the Theorem \ref{thm: quism invariance
of mc} there exists $\lambda^0 \in
\MC^2(\mathfrak{G}^0\otimes\mathfrak{m})_0$ such that
$\MC^2(\mathfrak{H}^0\otimes\mathfrak{m})_1(\theta_* \lambda^0,
\gamma^0)\ne \varnothing$. Let $\gimel^1 \in
\MC^2(\mathfrak{H}^0\otimes\mathfrak{m})_1(\theta_* \lambda,
\gamma)$. Applying Theorem \ref{thm: quism invariance of mc} again
we obtain that there exists $\mu \in
\MC^2(\mathfrak{G}^1\otimes\mathfrak{m})_1(\partial_0^0 \lambda^0,
\partial_1^0\lambda^0)$ such that there exists  $\phi \in \MC^2(\mathfrak{H}^1\otimes\mathfrak{m})_2(\gamma^1 \circ
\partial_0^0 \gimel^1, \partial_1^0 \gimel^1\circ
\theta_*\mu)$. Then $s_0^1 \phi$ is then a $2$-morphism $\gimel^1
\to \gimel^1 \circ (s_0^1 \mu)$, which induces a $2$-morphism $\psi:
 (s_0^1 \mu)^{-1}\to \id$.

 Set now $ \lambda^1= \mu \circ (\partial_0^0(s_0^1
\mu))^{-1} \in
\MC^2(\mathfrak{G}^1\otimes\mathfrak{m})_1(\partial_0^0 \lambda^0,
\partial_1^0\lambda^0)$, $\gimel^2 = \phi \otimes (\partial_0^0
\psi)^{-1}\in \MC^2(\mathfrak{H}^1\otimes\mathfrak{m})_2(\gamma^1
\circ
\partial_0^0 \gimel^1, \partial_1^0 \gimel^1\circ
\theta_*\lambda^1) $. It is easy to see that  $s_0^1 \lambda^1 =
\id$, $s_0^1 \gimel^2 =\id$.

 We then conclude that there exists a unique
$\lambda^2$ such that
\begin{equation*}
(\id \otimes \theta_*\lambda^2) \circ (\partial_2^1 \gimel^2 \otimes
\id) \circ (\id \otimes \partial_0^1 \gimel^2) =
\partial_1^1\gimel^2 \circ(\gamma^2 \otimes \id).
\end{equation*}
Such a $\lambda^2$ necessarily satisfies the  conditions
\begin{equation*}
\begin{split}
\partial_2^2 \lambda^2 \circ (\id \otimes \partial_0^2 \lambda^2) &=
\partial_1^2 \lambda^2 \circ (\partial_3^2 \lambda^2 \otimes \id).\\
s_0^2\lambda^2 &=s_1^2\lambda^2=\id
\end{split}
\end{equation*}
Therefore $\lambda=(\lambda^0, \lambda^1, \lambda^2) \in
\Stack(\mathfrak{G}\otimes\mathfrak{m})$, and $\gimel= (\gimel^1,
\gimel^2)$ defines a $1$-morphism $\theta_* \lambda \to \gamma$.
\end{proof}

\subsection{Acyclicity and strictness}
\begin{definition}
A $\mathfrak G$-stack $(\gamma^0,\gamma^1,\gamma^2)$ is called
\emph{strict} if $\partial^0_0\gamma^0 =
\partial^0_1\gamma^0$, $\gamma_1=\id$ and $\gamma_2=\id$.
\end{definition}
Let $\Stack_{Str}(\mathfrak{G})_0$ denote the subset of strict
$\mathfrak{G}$-stacks.

\begin{lemma}
$\Stack_{Str}(\mathfrak{G})_0 =
\MC^2(\ker(\mathfrak{G}^0\rightrightarrows\mathfrak{G}^1))_0$
\end{lemma}

\begin{definition}
For $\gamma_1,\gamma_2\in\Stack_{Str}(\mathfrak{G})_0$ a
$1$-morphism
$\gimel=(\gimel^1,\gimel^2)\in\Stack(\mathfrak{G})_1(\gamma_1,\gamma_2)$
is called \emph{strict} if
$\partial^0_0(\gimel^1)=\partial^0_1(\gimel^1)$ and $\gimel^2=\id$.
\end{definition}
For $\gamma_1,\gamma_2\in\Stack_{Str}(\mathfrak{G})_0$ we denote by
$\Stack_{Str}(\mathfrak{G})_1(\gamma_1,\gamma_2)$ the subset of
strict morphisms.

\begin{lemma}
For $\gamma_1,\gamma_2\in\Stack_{Str}(\mathfrak{G})_0$
$\Stack_{Str}(\mathfrak{G})_1(\gamma_1,\gamma_2) =
\MC^2(\ker(\mathfrak{G}^0\rightrightarrows\mathfrak{G}^1))_1$
\end{lemma}

For $\gamma_1,\gamma_2\in\Stack_{Str}(\mathfrak{G})_0$ let
$\Stack_{Str}(\mathfrak{G})(\gamma_1,\gamma_2)$ denote the full
subcategory of $\Stack(\mathfrak{G})(\gamma_1,\gamma_2)$ with
objects $\Stack_{Str}(\mathfrak{G})_1(\gamma_1,\gamma_2)$.

Thus, we have the $2$-groupoids $\Stack(\mathfrak{G})$ and
$\Stack_{Str}(\mathfrak{G})$ and an embedding of the latter into the
former which is fully faithful on the respective groupoids of
$1$-morphisms.

\begin{lemma}\label{lemma:strict stacks are mc ker}
$\Stack_{Str}(\mathfrak{G}) =
\MC^2(\ker(\mathfrak{G}^0\rightrightarrows\mathfrak{G}^1))$
\end{lemma}

Suppose that $\mathfrak{G}$ is a cosimplicial DGLA. For each $n$ and
$i$ we have the vector space $\mathfrak{G}^{n,i}$, namely the degree
$i$ component of $\mathfrak{G}^n$. The assignment
$n\mapsto\mathfrak{G}^{n,i}$ is a cosimplicial vector space
$\mathfrak{G}^{\bullet,i}$.

We will be considering the following acyclicity condition on the
cosimplicial DGLA $\mathfrak{G}$:
\begin{equation}\label{acyclicity condition}
\text{for all $i\in\mathbb{Z}$, $H^p(\mathfrak{G}^{\bullet,i})=0$
for $p\neq 0$}
\end{equation}

\begin{thm}\label{acyclic vs strictness}
Suppose that $\mathfrak{G}$ is a cosimplicial DGLA which satisfies
the condition \eqref{acyclicity condition}, and $\mathfrak{m}$ a
commutative nilpotent ring. Then, the functor
$\iota:\Stack_{Str}(\mathfrak{G}\otimes
\mathfrak{m})\to\Stack(\mathfrak{G}\otimes \mathfrak{m}))$ is an
equivalence.
\end{thm}
\begin{proof}
As we noted before, it is immediate from the definitions that if
$\gamma_1$, $\gamma_2 \in \Stack_{Str}(\mathfrak{G}\otimes
\mathfrak{m})_0$, and $\gimel_1, \gimel_2 \in
\Stack_{Str}(\mathfrak{G}\otimes \mathfrak{m})_1(\gamma_1,
\gamma_2)$, then $\iota: \Stack_{Str}(\mathfrak{G}\otimes
\mathfrak{m})_2(\gimel_1, \gimel_2)\to
\Stack_{Str}(\mathfrak{G}\otimes \mathfrak{m})_2(\gimel_1,
\gimel_2)$ is a bijection.

 Suppose now that
$\gamma_1,\gamma_2\in
\Stack_{Str}(\mathfrak{G}\otimes\mathfrak{m})_0$, $\gimel=(\gimel^1,
\gimel^2)
\in\Stack(\mathfrak{G}\otimes\mathfrak{m})_1(\gamma_1,\gamma_2)$. We
show that then there exists
$\daleth\in\Stack_{Str}(\mathfrak{G})_1(\gamma_1,\gamma_2)$ and a
$2$-morphism $\phi:\gimel\to\daleth$. Let $\gimel^2 =
\exp_{\partial_0^0 \gamma_2^0} g$, $g\in
(\mathfrak{G}^1\otimes\mathfrak{m}) $. Then $\partial_0^1 g -
\partial_1^1 g+\partial_2^1 g=0 \mod \mathfrak{m}^2$. By the
acyclicity condition there exists $a \in (\mathfrak{G}^0\otimes
\mathfrak{m})$ such that $\partial_0^0 a -\partial_1^0 a=g \mod
\mathfrak{m}^2$. Set $\phi_1 =\exp_{\gamma_2^0} a$. Define then
$\daleth_1 =(\phi_1\cdot \gimel^1, (\partial_1^0\phi_1 \otimes \id)
\circ \gimel^2 \circ(\id \otimes
\partial_0^0\phi_1)^{-1})
\in\Stack(\mathfrak{G}\otimes\mathfrak{m})_1(\gamma_1,\gamma_2)$.
Note that $\phi_1$ defines a $2$-morphism $\gimel\to \daleth_1$.
Note also  that $\daleth_1^2 \in \exp_{\partial_0^0
\gamma_2^0}(\mathfrak{G} \otimes \mathfrak{m}^2)$. Proceeding
inductively one constructs a sequence $\daleth_k
\in\Stack(\mathfrak{G}\otimes\mathfrak{m})_1(\gamma_1,\gamma_2)$
such that $\daleth_k^2 \in \exp(\mathfrak{G} \otimes
\mathfrak{m}^{k+1})$ and $2$-morphisms $\phi_k: \gimel\to
\daleth_k$, $\phi_{k+1}=\phi_k \mod \mathfrak{m}^k$. Since
$\mathfrak{m}$ is nilpotent, for $k$ large enough $\daleth_k \in
\Stack_{Str}(\mathfrak{G})_1(\gamma_1,\gamma_2)$.

Assume now that $\gamma \in
\Stack(\mathfrak{G}\otimes\mathfrak{m})_0$. We will construct
$\lambda \in \Stack_{Str}(\mathfrak{G}\otimes\mathfrak{m})_0$ and a
$1$-morphism $\gimel: \gamma \to \lambda$.

We begin by constructing $\mu \in
\Stack(\mathfrak{G}\otimes\mathfrak{m})_0$ such that $\mu^2=\id$ and
a $1$-morphism $\daleth: \gamma \to \lambda$. We have: $\gamma^2 =
\exp_{\partial_1^1 \partial_1^0 \gamma^0} c$, $c \in
(\mathfrak{G}^2\otimes\mathfrak{m})$. In view of the equations
\eqref{gamma2} $c$ satisfies the following identities:
\begin{equation} \begin{split}
   \partial_0^2 c -\partial_1^2 c +\partial_2^2 c -\partial_3^2 c=
 0& \mod \mathfrak{m}^2\\
s_0^2 c = s_1^2 c =0&
\end{split}\end{equation}

By the acyclicity of the normalized complex we can find $b \in
\mathfrak{G}^1 \otimes \mathfrak{m}$ such that $\partial_0^1 b
-\partial_1^1 b+\partial_2^1 b = c\mod \mathfrak{m}^2$, $s^1_0b=0$.
Let $\phi_1 = \exp_{\partial_1^0\gamma^0}(-b)$. Define then
$\mu_1=(\mu_1^0, \mu_1^1, \mu_1^2) $ where $\mu_1^0 =\gamma^0$,
$\mu_1^1 =\phi_1 \cdot \gamma^1$, and $\mu_1^2$ is such that
\begin{equation*}
(\id \otimes \gamma^2) \circ (\partial_2^1 \phi_1 \otimes \id) \circ
(\id \otimes \partial_0^1 \phi_1) =
\partial_1^1\phi_1 \circ(\mu_1^2 \otimes \id)
\end{equation*}
Note that $\mu_1^2 = \id \mod \mathfrak{m}^2$ and $(\id, \phi_1)$ is
 a $1$-morphism $\gamma \to \mu_1$. As before we can construct a
 sequence $\mu_k$ such that $\mu_k^2 = \id \mod \mathfrak{m}^{k+1}$,
 and $1$-morphisms $(\id, \phi_k): \gamma \to \mu_k$,
 $\phi_{k+1}=\phi_k \mod \mathfrak{m}^k$. As before we conclude that
this gives the desired construction of $\mu$. The rest of the proof,
i.e. the construction of $\lambda$ is completely analogous.
\end{proof}
\begin{cor}\label{cor: stacks are mc ker for acyclic}
Suppose that $\mathfrak{G}$ is a cosimplicial DGLA which satisfies
the condition \eqref{acyclicity condition}. Then there is a
canonical equivalence:
\[\Stack(\mathfrak{G}\otimes \mathfrak{m})\cong
\MC^2(\ker(\mathfrak{G}^0\rightrightarrows\mathfrak{G}^1)\otimes
\mathfrak{m})\]
\end{cor}
\begin{proof}
Combine Lemma \ref{lemma:strict stacks are mc ker} with Theorem
\ref{acyclic vs strictness}.
\end{proof}

\section{Algebroid stacks}\label{algstack}
In this section we review the notions of algebroid stack and twisted
form. We also define the notion of descent datum and relate it with
algebroid stacks.

\subsection{Algebroids and algebroid stacks}
\subsubsection{Algebroids}\label{subsubsection: algebroids}
For a category $\mathcal{C}$ we denote by $i\mathcal{C}$ the
subcategory of isomorphisms in $\mathcal{C}$; equivalently,
$i\mathcal{C}$ is the maximal subgroupoid in $\mathcal{C}$.

Suppose that $R$ is a commutative $k$-algebra.

\begin{definition}
An \emph{$R$-algebroid} is a \emph{nonempty} $R$-linear category
$\mathcal{C}$ such that the groupoid $i\mathcal{C}$ is connected
\end{definition}

Let $\Algd_R$ denote the $2$-category of $R$-algebroids
(full $2$-subcategory of the $2$-category of $R$-linear categories).

Suppose that $A$ is an $R$-algebra. The $R$-linear category with one
object and morphisms $A$ is an $R$-algebroid denoted $A^+$.

Suppose that $\mathcal{C}$ is an $R$-algebroid and $L$ is an object
of $\mathcal{C}$. Let $A = \End_\mathcal{C}(L)$. The functor
$A^+\to\mathcal{C}$ which sends the unique object of $A^+$ to $L$ is
an equivalence.

Let $\Alg^2_R$ denote the $2$-category of with
\begin{itemize}
\item objects $R$-algebras
\item $1$-morphisms homomorphism of $R$-algebras
\item $2$-morphisms $\phi\to\psi$, where $\phi,\psi : A\to B$ are two $1$-morphisms are elements $b\in B$ such that $\phi(a)\cdot b = b\cdot\psi(a)$ for all $a\in A$.
\end{itemize}
It is clear that the $1$- and the $2$- morphisms in $\Alg^2_R$ as
defined above induce $1$- and $2$-morphisms of the corresponding
algebroids under the assignment $A \mapsto A^+$. The structure of a
$2$-category on $\Alg^2_R$ (i.e. composition of $1$- and $2$-
morphisms) is determined by the requirement that the assignment
$A\mapsto A^+$ extends to an embedding $(\bullet)^+ : \Alg^2_R\to
\Algd_R$.

Suppose that $R\to S$ is a morphism of commutative $k$-algebras. The
assignment $A\to A\otimes_R S$ extends to a functor
$(\bullet)\otimes_R S : \Alg^2_R\to \Alg^2_S$.

\subsubsection{Algebroid stacks}
We refer the reader to \cite{SGA1} and \cite{Vist} for basic
definitions. We will use the notion of fibered category
interchangeably with that of a pseudo-functor. A \emph{prestack}
$\mathcal C$ on a space $X$ is a category fibered over the category
of open subsets of $X$, equivalently, a pseudo-functor $U\mapsto
\mathcal{C}(U)$, satisfying the following additional requirement.
For an open subset $U$ of $X$ and two objects $A,B\in\mathcal{C}(U)$
we have the presheaf $\shHom_\mathcal{C}(A,B)$ on $U$ defined by
$U\supseteq V\mapsto \Hom_{\mathcal{C}(V)}(A\vert_V, B\vert_B)$. The
fibered category $\mathcal{C}$ is a prestack if for any $U$,
$A,B\in\mathcal{C}(U)$, the presheaf $\shHom_\mathcal{C}(A,B)$ is a
sheaf. A prestack is a \emph{stack} if, in addition, it satisfies
the condition of effective descent for objects. For a prestack
$\mathcal{C}$ we denote the associated stack by
$\widetilde{\mathcal{C}}$.

\begin{definition}
A stack in $R$-linear categories $\mathcal{C}$ on $X$ is an
\emph{$R$-algebroid stack} if it is locally nonempty and locally
connected, i.e. satisfies
\begin{enumerate}
\item any point $x\in X$ has a neighborhood $U$ such that
$\mathcal{C }(U)$ is nonempty;

\item for any $U\subseteq X$, $x\in U$, $A, B\in\mathcal{C}(U)$
there exits a neighborhood $V\subseteq U$ of $x$ and an isomorphism
$A\vert_V\cong B\vert_V$.
\end{enumerate}
\end{definition}

\begin{remark}
Equivalently, the stack associated to the substack of isomorphisms
$\widetilde{i\mathcal{C}}$ is a gerbe.
\end{remark}

\begin{example}
Suppose that $\mathcal{A}$ is a sheaf of $R$-algebras on $X$. The
assignment $X\supseteq U\mapsto\mathcal{A}(U)^+$ extends in an
obvious way to a prestack in $R$-algebroids denoted $\mathcal{A}^+$.
The associated stack $\widetilde{\mathcal{A}^+}$ is canonically
equivalent to the stack of locally free $\mathcal{A}^{op}$-modules
of rank one. The canonical morphism
$\mathcal{A}^+\to\widetilde{\mathcal{A}^+}$ sends the unique
(locally defined) object of $\mathcal{A}^+$ to the free module of
rank one.
\end{example}

$1$-morphisms and $2$-morphisms of $R$-algebroid stacks are those of
stacks in $R$-linear categories. We denote the $2$-category of
$R$-algebroid stacks by $\AlgStack_R(X)$.

\subsection{Descent data}
\subsubsection{Convolution data}
\begin{definition}\label{definition: convolution datum}
An \emph{$R$-linear convolution datum} is a triple
$(\mathcal{U},\mathcal{A}_{01},\mathcal{A}_{012})$ consisting
of:
\begin{itemize}
\item a cover $\mathcal{U}\in\Cov(X)$

\item a sheaf $\mathcal{A}_{01}$ of $R$-modules $\mathcal{A}_{01}$ on $N_1\mathcal{U}$

\item a morphism
\begin{equation}\label{convolution}
\mathcal{A}_{012} : (\pr_{01}^2)^*\mathcal{A}_{01}\otimes_R
(\pr_{12}^2)^*\mathcal{A}_{01} \to (\pr_{02}^2)^*\mathcal{A}_{01}
\end{equation}
of $R$-modules
\end{itemize}
subject to the associativity condition expressed by the
commutativity of the diagram
\begin{equation*}
\begin{CD}
(\pr_{01}^3)^*\mathcal{A}_{01}\otimes_R(\pr_{12}^3)^*\mathcal{A}_{01}\otimes_R
(\pr_{23}^3)^*\mathcal{A}_{01}
@>{(\pr_{012}^3)^*(\mathcal{A}_{012})\otimes\id}>>
(\pr_{02}^3)^*\mathcal{A}_{01}\otimes_R (\pr_{23}^3)^*\mathcal{A}_{01} \\
@V{\id\otimes(\pr_{123}^3)^*(\mathcal{A}_{012})}VV @VV{(\pr_{023}^3)^*(\mathcal{A}_{012})}V \\
(\pr_{01}^3)^*\mathcal{A}_{01}\otimes_R
(\pr_{13}^3)^*\mathcal{A}_{01}
@>{(\pr_{013}^3)^*(\mathcal{A}_{012})}>>
(\pr_{03}^3)^*\mathcal{A}_{01}
\end{CD}
\end{equation*}
\end{definition}

For a convolution datum $(\mathcal{U}, \mathcal{A}_{01},
\mathcal{A}_{012})$ we denote by $\underline{\mathcal{A}}$ the
pair $(\mathcal{A}_{01}, \mathcal{A}_{012})$ and  abbreviate
the convolution datum by
$(\mathcal{U},\underline{\mathcal{A}})$.

For a convolution datum $(\cu,\uca)$ let
\begin{itemize}
\item $\mathcal{A} := (\pr_{00}^0)^*\mathcal{A}_{01}$; $\mathcal{A}$
is a sheaf of $R$-modules on $N_0\mathcal{U}$

\item $\mathcal{A}_i^p := (\pr^p_i)^*\mathcal{A}$; thus for every $p$ we
get sheaves $\ca_i^p$, $0 \le i \le p$ on $N_p\cu$.
\end{itemize}

The identities $\pr_{01}^0\circ\pr^0_{000} =
\pr_{12}^0\circ\pr_{000}^0= \pr_{02}\circ \pr^0_{000} = \pr_{00}^0$
imply that the pull-back of $\mathcal{A}_{012}$ to $N_0\mathcal{U}$
by $\pr^0_{000}$ gives the pairing
\begin{equation}\label{pairing A A}
(\pr^0_{000})^*(\mathcal{A}_{012}) :
\mathcal{A}\otimes_R\mathcal{A}\to\mathcal{A} \ .
\end{equation}
The associativity condition implies that the pairing \eqref{pairing
A A} endows $\mathcal{A}$ with a structure of a sheaf of associative
$R$-algebras on $N_0\mathcal{U}$.

The sheaf $\mathcal{A}_i^p$ is endowed with the associative
$R$-algebra structure induced by that on $\mathcal{A}$. We denote by
$\mathcal{A}_{ii}^p$ the $\mathcal{A}_i^p\otimes_R
(\mathcal{A}_i^p)^{op}$-module $\mathcal{A}_i^p$, with the module
structure given by the left and right multiplication.

The identities $\pr^2_{01}\circ\pr^1_{001} =
\pr^0_{00}\circ\pr^1_0$, $\pr^2_{12}\circ\pr^1_{001} =
\pr^2_{02}\circ\pr^1_{001} = \id$ imply that the pull-back of
$\mathcal{A}_{012}$ to $N_1\mathcal{U}$ by $\pr^1_{001}$ gives the
pairing
\begin{equation}\label{pairing A M}
(\pr^1_{001})^*(\mathcal{A}_{012}) :
\mathcal{A}_0^1\otimes_R\mathcal{A}_{01}\to \mathcal{A}_{01}
\end{equation}
The associativity condition implies that the pairing \eqref{pairing
A M} endows $\mathcal{A}_{01}$ with a structure of a
$\mathcal{A}_0^1$-module. Similarly, the pull-back of
$\mathcal{A}_{012}$ to $N_1\mathcal{U}$ by $\pr^1_{011}$ endows
$\mathcal{A}_{01}$ with a structure of a
$(\mathcal{A}_1^1)^{op}$-module. Together, the two module structures
define a structure of a $\mathcal{A}_0^1\otimes_R
(\mathcal{A}_1^1)^{op}$-module on $\mathcal{A}_{01}$.

The map \eqref{convolution} factors through the map
\begin{equation}\label{factor convolution}
(\pr_{01}^2)^*\mathcal{A}_{01}\otimes_{\ca_1^2}
(\pr^2_{12})^*\mathcal{A}_{01} \to (\pr_{02}^2)^*\mathcal{A}_{01}
\end{equation}

\begin{definition}\label{definition: unit}
A \emph{unit} for a convolution datum $\underline{\mathcal{A}}$ is a
morphism of $R$-modules
\begin{equation*}
 \vac : R\to\mathcal{A}
\end{equation*}
such that the compositions
\begin{equation*}
\mathcal{A}_{01} \xrightarrow{\vac\otimes\id}
\mathcal{A}_0^1\otimes_R\mathcal{A}_{01}
\xrightarrow{(\pr^1_{001})^*(\mathcal{A}_{012})} \mathcal{A}_{01}
\end{equation*}
and
\begin{equation*}
\mathcal{A}_{01} \xrightarrow{\id\otimes\vac}
\mathcal{A}_{01}\otimes_R \mathcal{A}_1^1
\xrightarrow{(\pr^1_{011})^*(\mathcal{A}_{012})} \mathcal{A}_{01}
\end{equation*}
are equal to the respective identity morphisms.
\end{definition}

\subsubsection{Descent data}\label{definition: descent datum}
\begin{definition}
A \emph{descent datum} on $X$ is an $R$-linear convolution
datum $(\mathcal{U},\underline{\mathcal{A}})$ on $X$ with a
unit which satisfies the following additional conditions:
\begin{enumerate}
\item $\mathcal{A}_{01}$ is locally free of rank one as a $\mathcal{A}_0^1$-module
and as a $(\mathcal{A}_1^1)^{op}$-module;

\item the map \eqref{factor convolution} is an isomorphism.
\end{enumerate}
\end{definition}

\subsubsection{$1$-morphisms}\label{1-morphisms in desc}
Suppose given convolution data $(\cu,\underline{\mathcal{A}})$ and
$(\cu,\underline{\mathcal{B}})$ as in Definition \ref{definition:
convolution datum}.

\begin{definition}
A \emph{$1$-morphism of convolution data}
\begin{equation*}
\underline{\phi} : (\cu,\underline{\mathcal{A}})\to
(\cu,\underline{\mathcal{B}})
\end{equation*}
is a morphism of $R$-modules $\phi_{01} :
\mathcal{A}_{01}\to\mathcal{B}_{01}$ such that the diagram
\begin{equation}\label{compatibility with convolution}
\begin{CD}
(\pr^2_{01})^*\mathcal{A}_{01}\otimes_R (\pr^2_{12})^*\mathcal{A}_{01} @>{\mathcal{A}_{012}}>> (\pr^2_{02})^*\mathcal{A}_{01} \\
@V{(\pr^2_{01})^*(\phi_{01})\otimes (\pr^2_{12})^*(\phi_{01})}VV @VV{(\pr^2_{02})^*(\phi_{01})}V \\
(\pr^2_{01})^*\mathcal{B}_{01}\otimes_R
(\pr^2_{12})^*\mathcal{B}_{01} @>{\mathcal{B}_{012}}>>
(\pr^2_{02})^*\mathcal{B}_{01}
\end{CD}
\end{equation}
is commutative.

\end{definition}

The $1$-morphism $\underline{\phi}$ induces a morphism of
$R$-algebras $\phi := (\pr_{00}^0)^*(\phi_{01}) :
\mathcal{A}\to \mathcal{B}$ on $N_0\mathcal{U}$ as well as
morphisms $\phi_i^p := (\pr_i^p)^*(\phi)
:\mathcal{A}_i^p\to\mathcal{B}_i^p$ on $N_p \cu$. The morphism
$\phi_{01}$ is compatible with the morphism of algebras
$\phi_0\otimes\phi_1^{op} : \mathcal{A}_0^1\otimes_R
(\mathcal{A}_1^1)^{op}\to \mathcal{B}_0^1\otimes_R
(\mathcal{B}_1^1)^{op}$ and the respective module structures.

A \emph{$1$-morphism of descent data} is a $1$-morphism of the
underlying convolution data which preserves respective units.

\subsubsection{$2$-morphisms}\label{2-morphisms in desc}
Suppose that we are given descent data
$(\cu,\underline{\mathcal{A}})$ and
$(\cu,\underline{\mathcal{B}})$ as in \ref{definition: descent
datum} and two $1$-morphisms
\begin{equation*}
\underline{\phi}, \underline{\psi} :
(\cu,\underline{\mathcal{A}})\to (\cu,\underline{\mathcal{B}})
\end{equation*}
as in \ref{1-morphisms in desc}.

A $2$-morphism
\begin{equation*}
\underline{b} : \underline{\phi}\to \underline{\psi}
\end{equation*}
is a section $b\in\Gamma(N_0\mathcal{U};\mathcal{B})$ such that the
diagram
\begin{equation*}
\begin{CD}
\mathcal{A}_{01} @>{b\otimes\phi_{01}}>> \mathcal{B}_0^1\otimes_R\mathcal{B}_{01} \\
@V{\psi_{01}\otimes b}VV @VVV \\
\mathcal{B}_{01}\otimes_R\mathcal{B}_1^1 @>>> \mathcal{B}_{01}
\end{CD}
\end{equation*}
is commutative.

\subsubsection{The $2$-category of descent data}
Fix a cover $\cu$ of $X$.

Suppose that we are  given descent data $(\cu,\uca)$,
$(\cu,\underline{\mathcal{B}})$,
$(\cu,\underline{\mathcal{C}})$ and $1$-morphisms
$\underline{\phi} : (\cu,\uca)\to
(\cu,\underline{\mathcal{B}})$ and $\underline{\psi} :
(\cu,\underline{\mathcal{B}})\to
(\cu,\underline{\mathcal{C}})$. The map
$\psi_{01}\circ\phi_{01}: \ca_{01}\to\mathcal{C}_{01}$ is a
$1$-morphism of descent data
$\underline{\psi}\circ\underline{\phi} : (\cu,\uca)\to
(\cu,\underline{\mathcal{C}})$, the \emph{composition of
$\underline{\phi}$ and $\underline{\psi}$}.

Suppose that $\underline{\phi}^{(i)} : (\cu,\uca)\to
(\cu,\underline{\mathcal{B}})$, $i=1,2,3$, are $1$-morphisms and
$\underline{b}^{(j)}:
\underline{\phi}^{(j)}\to\underline{\phi}^{(j+1)}$, $j=1,2$, are
$2$-morphisms. The section $b^{(2)}\cdot
b^{(1)}\in\Gamma(N_0\mathcal{U};\mathcal{B})$ defines a
$2$-morphism, denoted $\underline{b}^{(2)}\underline{b}^{(1)} :
\underline{\phi}^{(1)}\to\underline{\phi}^{(3)}$, the \emph{vertical
composition of $\underline{b}^{(1)}$ and $\underline{b}^{(2)}$}.

Suppose that $\underline{\phi}^{(i)} : (\cu,\uca)\to
(\cu,\underline{\mathcal{B}})$, $\underline{\psi}^{(i)} :
(\cu,\underline{\mathcal{B}})\to (\cu,\underline{\mathcal{C}})$,
$i=1,2$, are $1$-morphisms and $\underline{b}:
\underline{\phi}^{(1)}\to\underline{\phi}^{(2)}$, $\underline{c}:
\underline{\psi}^{(1)}\to\underline{\psi}^{(2)}$ are $2$-morphisms.
The section $c\cdot
\psi^{(1)}(b)\in\Gamma(N_0\mathcal{U};\mathcal{C})$ defines a
$2$-morphism, denoted $\underline{c}\otimes\underline{b}$, the
\emph{horizontal composition of $\underline{b}$ and
$\underline{c}$}.

We leave it to the reader to check that with the compositions
defined above descent data, $1$-morphisms and $2$-morphisms form a
$2$-category, denoted $\Desc_R(\cu)$.

\subsubsection{Fibered category of descent data}
Suppose that $\rho : \mathcal{V}\to \mathcal{U}$ is a morphism of
covers and $(\cu,\uca)$ is a descent datum. Let $\ca^\rho_{01} =
(N_1\rho)^*\ca_{01}$, $\ca^\rho_{012} = (N_2\rho)^*(\ca_{012})$.
Then, $(\mathcal{V},\uca^\rho)$ is a descent datum. The assignment
$(\cu,\uca)\mapsto(\mathcal{V},\uca^\rho)$ extends to a functor,
denoted $\rho^* : \Desc_R(\cu)\to \Desc_R(\mathcal{V})$.

The assignment $\Cov(X)^{op}\ni\cu\to\Desc_R(\cu)$, $\rho\to\rho^*$
is (pseudo-)functor. Let $\Desc_R(X)$ denote the corresponding
$2$-category fibered in $R$-linear $2$-categories over $\Cov(X)$
with object pairs $(\cu,\uca)$ with $\cu\in\Cov(X)$ and
$(\cu,\uca)\in\Desc_R(\cu)$; a morphism $(\cu',\uca')\to(\cu,\uca)$
in $\Desc_R(X)$ is a pair $(\rho,\underline{\phi})$, where $\rho :
\cu'\to\cu$ is a morphism in $\Cov(X)$ and $\underline{\phi} :
(\cu',\uca')\to\rho^*(\cu,\uca)=(\cu',\uca^\rho)$.

\subsection{Trivializations}
\subsubsection{Definition of a trivialization}
\begin{definition}
A \emph{trivialization} of an algebroid stack $\mathcal{C}$ on $X$
is an object in $\mathcal{C}(X)$.
\end{definition}

Suppose that $\mathcal{C}$ is an algebroid stack on $X$ and
$L\in\mathcal{C}(X)$ is a trivialization. The object $L$ determines
a morphism $\shEnd_{\mathcal{C}}(L)^+\to\mathcal{C}$.

\begin{lemma}
The induced morphism
$\widetilde{\shEnd_{\mathcal{C}}(L)^+}\to\mathcal{C}$ is an
equivalence.
\end{lemma}

\begin{remark}
Suppose that $\mathcal{C}$ is an $R$-algebroid stack on $X$. Then,
there exists a cover $\mathcal{U}$ of $X$ such that the stack
$\epsilon_0^*\mathcal{C}$ admits a trivialization.
\end{remark}

\subsubsection{The $2$-category of trivializations}
Let $\Triv_R(X)$ denote the $2$-category with
\begin{itemize}
\item objects the triples $(\mathcal{C},\mathcal{U},L)$ where $\mathcal{C}$
is an $R$-algebroid stack on $X$, $\mathcal{U}$ is an open cover of
$X$ such that $\epsilon_0^*\mathcal{C}(N_0\mathcal{U})$ is nonempty
and $L$ is a trivialization of $\epsilon_0^*\mathcal{C}$.

\item $1$-morphisms $(\mathcal{C}',\mathcal{U}',L)\to (\mathcal{C},\mathcal{U},L)$
are pairs $(F,\rho)$ where $\rho : \cu'\to\cu$ is a morphism of
covers, $F: \mathcal{C}'\to\mathcal{C}$ is a functor such that
$(N_0\rho)^*F(L') = L$

\item $2$-morphisms $(F,\rho)\to (G,\rho)$, where $(F,\rho),(G,\rho):(\mathcal{C}',\mathcal{U}',L)\to (\mathcal{C},\mathcal{U},L)$, are the morphisms of functors $F\to G$.
\end{itemize}
The assignment $(\mathcal{C},\mathcal{U},L)\mapsto\mathcal{C}$
extends in an obvious way to a functor
$\Triv_R(X)\to\AlgStack_R(X)$.

The assignment $(\mathcal{C},\mathcal{U},L)\mapsto\mathcal{U}$ extends to a functor $\Triv_R(X)\to\Cov(X)$ making $\Triv_R(X)$ a fibered $2$-category over $\Cov(X)$. For $\cu\in\Cov(X)$ we denote the fiber over $\cu$ by $\Triv_R(X)(U)$.

\subsubsection{Algebroid stacks from descent
data}\label{subsubsection: descent data to stacks} Consider
$(\mathcal{U},\uca)\in\Desc_R(\cu)$.

The sheaf of algebras $\ca$ on $N_0\mathcal{U}$ gives rise to the
algebroid stack $\widetilde{\ca^+}$. The sheaf $\mathcal{A}_{01}$
defines an equivalence
\begin{equation*}
\phi_{01} := (\bullet)\otimes_{\ca_0^1}\ca_{01} :
(\pr^1_0)^*\widetilde{\ca^+} \to (\pr^1_1)^*\widetilde{\ca^+} \ .
\end{equation*}
The convolution map $\ca_{012}$ defines an isomorphism of functors
\begin{equation*}
\phi_{012} :
(\pr_{01}^2)^*(\phi_{01})\circ(\pr_{12}^2)^*(\phi_{01})\to(\pr_{02}^2)^*(\phi_{01})
\ .
\end{equation*}
We leave it to the reader to verify that the triple
$(\widetilde{\ca^+},\phi_{01}, \phi_{012})$ constitutes a descent
datum for an algebroid stack on $X$ which we denote by
$\stack(\cu,\uca)$.

By construction there is a canonical equivalence $\widetilde{\uca^+}\to\epsilon_0^*\stack(\cu,\uca)$ which endows $\epsilon_0^*\stack(\cu,\uca)$ with a canonical trivialization $\one$.

The assignment $(\cu,\uca)\mapsto (\stack(\cu,\uca),\cu,\one)$ extends to a cartesian functor
\[
\stack : \Desc_R(X)\to\Triv_R(X) \ .
\]

\subsubsection{Descent data from trivializations}\label{subsubsection: descent data for algebroid stacks}
Consider $(\mathcal{C},\mathcal{U},L)\in\Triv_R(X)$. Since
$\epsilon_0\circ\pr^1_0 = \epsilon_0\circ\pr^1_1 = \epsilon_1$ we
have canonical identifications
$(\pr^1_0)^*\epsilon_0^*\mathcal{C}\cong
(\pr_1^1)^*\epsilon_0^*\mathcal{C} \cong \epsilon_1^*\mathcal{C}$.
The object $L\in\epsilon_0^*\mathcal{C}(N_0\mathcal{U})$ gives rise
to the objects $(\pr_0^1)^*L$ and $(\pr_1^1)^*L$ in
$\epsilon_1^*\mathcal{C}(N_0\mathcal{U})$. Let $\mathcal{A}_{01} =
\shHom_{\epsilon_1^*\mathcal{C}}((\pr^1_1)^*L,(\pr^1_0)^*L)$. Thus,
$\mathcal{A}_{01}$ is a sheaf of $R$-modules on $N_1\mathcal{U}$.

The object $L\in\epsilon_0^*\mathcal{C}(N_0\mathcal{U})$ gives rise
to the objects $(\pr^2_0)^*L$, $(\pr^2_1)^*L$ and $(\pr_2^2)^*L$ in
$\epsilon_2^*\mathcal{C}(N_2\mathcal{U})$. There are canonical
isomorphisms
$(\pr_{ij}^2)^*\mathcal{A}_{01}\cong\shHom_{\epsilon_2^*\mathcal{C}}((\pr^2_i)^*L,(\pr^2_j)^*L)$.
The composition of morphisms
\begin{equation*}
\shHom_{\epsilon_2^*\mathcal{C}}((\pr^2_1)^*L,(\pr^2_0)^*L)\otimes_R
\shHom_{\epsilon_2^*\mathcal{C}}((\pr^2_2)^*L,(\pr^2_1)^*L) \to
\shHom_{\epsilon_2^*\mathcal{C}}((\pr^2_2)^*L,(\pr^2_0)^*L)
\end{equation*}
gives rise to the map
\begin{equation*}
\mathcal{A}_{012} : (\pr_{01}^2)^*\mathcal{A}_{01}\otimes_R
(\pr_{12}^2)^*\mathcal{A}_{01} \to (\pr_{02}^2)^*\mathcal{A}_{02}
\end{equation*}

Since $\pr^1_i\circ\pr^0_{00}=\id$ there is a canonical isomorphism
$\mathcal{A} := (\pr^0_{00})^*\mathcal{A}_{01}\cong\shEnd(L)$ which
supplies $\mathcal{A}$ with the unit section $\vac :
R\xrightarrow{1}\shEnd(L)\to \mathcal{A}$.

The pair $(\mathcal{U},\underline{\mathcal{A}})$, together with the section
$\vac$ is a decent datum which we denote $\desc(\mathcal{C},\mathcal{U},L)$.

The assignment $(\mathcal{U},\underline{\mathcal{A}})\mapsto \desc(\mathcal{C},\mathcal{U},L)$ extends to a cartesian functor
\[
\desc : \Triv_R(X)\to\Desc_R(X) \ .
\]

\begin{lemma}
The functors $\stack$ and $\desc$ are mutually quasi-inverse equivalences.
\end{lemma}

\subsection{Base change}
For an $R$-linear category $\mathcal{C}$ and homomorphism of
algebras $R\to S$ we denote by $\mathcal{C}\otimes_R S$ the category
with the same objects as $\mathcal{C}$ and morphisms defined by
$\Hom_{\mathcal{C}\otimes_R S}(A,B) = \Hom_\mathcal{C}(A,B)\otimes_R
S$.

For an $R$-algebra $A$ the categories $(A\otimes_R S)^+$ and
$A^+\otimes_R S$ are canonically isomorphic.

For a prestack $\mathcal{C}$ in $R$-linear categories we denote by
$\mathcal{C}\otimes_R S$ the prestack associated to the fibered
category $U\mapsto\mathcal{C}(U)\otimes_R S$.

For $U\subseteq X$, $A,B\in\mathcal{C}(U)$, there is an isomorphism
of sheaves $\shHom_{\mathcal{C}\otimes_R S}(A,B) =
\shHom_\mathcal{C}(A,B)\otimes_R S$.

\begin{lemma}
Suppose that $\ca$ is a sheaf of $R$-algebras and $\mathcal{C}$ is
an $R$-algebroid stack.
\begin{enumerate}
\item $(\widetilde{\ca^+}\otimes_R S)\widetilde{~~}$ is an algebroid stack equivalent to $\widetilde{(\ca\otimes_R S)^+}$ .

\item $\widetilde{\mathcal{C}\otimes_R S}$ is an algebroid stack.
\end{enumerate}
\end{lemma}
\begin{proof}
Suppose that $\ca$ is a sheaf of $R$-algebras. There is a canonical
isomorphism of prestacks $(\ca\otimes_R S)^+\cong \ca^+\otimes_R S$
which induces the canonical equivalence $\widetilde{(\ca\otimes_R
S)^+}\cong \widetilde{\ca^+\otimes_R S}$.

The canonical functor $\ca^+\to\widetilde{\ca^+}$ induces the
functor $\ca^+\otimes_R S \to\widetilde{\ca^+}\otimes_R S$, hence
the functor $\widetilde{\ca^+\otimes_R S} \to
(\widetilde{\ca^+}\otimes_R S)\widetilde{~~}$.

The map $\ca\to\ca\otimes_R S$ induces the functor
$\widetilde{\ca^+}\to \widetilde{(\ca\otimes_R S)^+}$ which
factors through the functor $\widetilde{\ca^+}\otimes_R S\to
\widetilde{(\ca\otimes_R S)^+}$. From this we obtain the
functor $(\widetilde{\ca^+}\otimes_R S)\widetilde{~~}\to
\widetilde{(\ca\otimes_R S)^+}$.

We leave it to the reader to check that the two constructions are
mutually inverse equivalences. It follows that
$(\widetilde{\ca^+}\otimes_R S)\widetilde{~~}$ is an algebroid stack
equivalent to $\widetilde{(\ca\otimes_R S)^+}$.

Suppose that $\mathcal{C}$ is an $R$-algebroid stack. Let $\cu$ be a
cover such that $\epsilon_0^*\mathcal{C}(N_0\cu)$ is nonempty. Let
$L$ be an object in $\epsilon_0^*\mathcal{C}(N_0\cu)$; put $\ca :=
\shEnd_{\epsilon_0^*\mathcal{C}}(L)$. The equivalence
$\widetilde{\ca^+}\to\epsilon_0^*\mathcal{C}$ induces the
equivalence $(\widetilde{\ca^+}\otimes_R S)\widetilde{~~}\to
(\epsilon_0^*\mathcal{C}\otimes_R S)\widetilde{~~}$. Since the
former is an algebroid stack so is the latter. There is a canonical
equivalence $(\epsilon_0^*\mathcal{C}\otimes_R S)\widetilde{~~}\cong
\epsilon_0^*(\widetilde{\mathcal{C}\otimes_R S})$; since the former
is an algebroid stack so is the latter. Since the property of being
an algebroid stack is local, the stack
$\widetilde{\mathcal{C}\otimes_R S}$ is an algebroid stack.
\end{proof}

\subsection{Twisted forms}\label{tworms}
 Suppose that $\mathcal{A}$ is a sheaf of
$R$-algebras on $X$. We will call an $R$-algebroid stack
locally equivalent to $\widetilde{\mathcal{A}^{op +}}$ a
\emph{twisted form of $\mathcal{A}$}.

Suppose that $\mathcal{S}$ is twisted form of $\mathcal{O}_X$. Then,
the substack $i\mathcal{S}$ is an $\mathcal{O}_X^\times$-gerbe. The
assignment $\mathcal{S}\mapsto i\mathcal{S}$ extends to an
equivalence between the $2$-groupoid of twisted forms of
$\mathcal{O}_X$ (a subcategory of $\AlgStack_\mathbb{C}(X)$) and the
$2$-groupoid of $\mathcal{O}_X^\times$-gerbes.

Let $\mathcal{S}$ be a twisted form of $\mathcal{O}_X$. Then
for any $U\subseteq X$, $A\in \mathcal{S}(U)$ the canonical map
$\mathcal{O}_U\to\shEnd_\mathcal{S}(A)$ is an isomorphism.
Consequently, if $\mathcal{U}$ is a cover of $X$ and $L$ is a
trivialization of $\epsilon_0^*\mathcal{S}$, then there is a
canonical isomorphism of sheaves of algebras
$\mathcal{O}_{N_0\cu}\to\shEnd_{\epsilon_0^*\mathcal{S}}(L)$.

Conversely, suppose that $(\cu,\uca)$ is a $\mathbb{C}$-descent
datum. If the sheaf of algebras $\ca$ is isomorphic to
$\mathcal{O}_{N_0\cu}$ then such an isomorphism is unique since the
latter has no non-trivial automorphisms. Thus, we may and will
identify $\ca$ with $\mathcal{O}_{N_0\cu}$. Hence, $\ca_{01}$ is a
line bundle on $N_1\cu$ and the convolution map $\ca_{012}$ is a
morphism of line bundles. The stack which corresponds to
$(\cu,\uca)$ (as in \ref{subsubsection: descent data to stacks}) is
a twisted form of $\mathcal{O}_X$.

Isomorphism classes of twisted forms of $\mathcal{O}_X$ are
classified by $H^2(X;\mathcal{O}^\times_X)$. We recall the
construction presently. Suppose that the twisted form $\mathcal{S}$
of $\mathcal{O}_X$ is represented by the descent datum
$(\mathcal{U},\underline{\mathcal{A}})$. Assume in addition that the
line bundle $\mathcal{A}_{01}$ on $N_1\mathcal{U}$ is trivialized.
Then we can consider $\mathcal{A}_{012}$ as an element  in
$\Gamma(N_2\mathcal{U};\mathcal{O}^\times)$. The associativity
condition implies that $\mathcal{A}_{012}$ is a cocycle in $\check
C^2(\mathcal{U}; \mathcal{O}^\times)$. The class of this cocycle in
$\check H^2(\mathcal{U}; \mathcal{O}^\times)$ does not depend on the
choice of trivializations of the line bundle $\mathcal{A}_{01}$ and
yields a class in $H^2(X;\mathcal{O}^\times_X)$.

We can write this class using the de Rham complex for jets. We refer
to the Section \ref{remind jet} for the notations and a brief
review. The composition $\mathcal{O}^\times\to
\mathcal{O}^\times/\mathbb{C}^\times\xrightarrow{\log}
\mathcal{O}/\mathbb{C}\xrightarrow{j^\infty} \DR(\mathcal{J/O})$
induces the map $H^2(X;\mathcal{O}^\times)\to
H^2(X;\DR(\mathcal{J/O}))\cong
H^2(\Gamma(X;\Omega^\bullet_X\otimes\mathcal{J}_X/\mathcal{O}_X),
\nabla^{can})$. Here the latter isomorphism follows from   the fact
that the sheaf $\Omega^\bullet_X\otimes\mathcal{J}_X/\mathcal{O}_X$
is soft.
We denote by $[\mathcal{S}]\in
H^2(\Gamma(X;\Omega^\bullet_X\otimes\mathcal{J}_X/\mathcal{O}_X),
\nabla^{can})$ the image of the class of $\mathcal{S}$ in
$H^2(X;\mathcal{O}^\times)$. In the Lemma \ref{defw} (see also Lemma
\ref{propw}) we will construct an explicit representative for
$[\mathcal{S}]$.

\section{DGLA of local cochains on matrix algebras}\label{dglaloc}
In this section we define matrix algebras from a descent datum and
use them to construct a cosimplicial DGLA of local cochains. We also
establish the acyclicity of this cosimplicial DGLA.

\subsection{Definition of matrix algebras}\label{subsection: matrix algebras}
\subsubsection{Matrix entries}
Suppose that $(\cu,\uca)$ is an $R$-descent datum. Let
$\mathcal{A}_{10} := \tau^*\mathcal{A}_{01}$, where $\tau
=\pr^1_{10}: N_1\mathcal{U}\to N_1\mathcal{U}$ is the transposition
of the factors. The pairings $(\pr^1_{100})^*(\mathcal{A}_{012}):
\ca_{10} \otimes_R \ca^1_0 \to \ca_{10} $ and
$(\pr^1_{110})^*(\mathcal{A}_{012}):\ca^1_1 \otimes_R\ca_{10}\to
\ca_{10}$ of sheaves on $N_1\cu$ endow $\mathcal{A}_{10}$ with a
structure of a
$\mathcal{A}^1_1\otimes(\mathcal{A}^1_0)^{op}$-module.

The identities $\pr^2_{01}\circ\pr^1_{010} = \id$,
$\pr^2_{12}\circ\pr^1_{010}=\tau$ and $\pr^2_{02}\circ\pr^1_{010} =
\pr^0_{00}\circ\pr^1_0$ imply that the pull-back of $\ca_{012}$ by
$\pr^1_{010}$ gives the pairing
\begin{equation}\label{pairing 01 10}
(\pr^1_{010})^*(\mathcal{A}_{012}) :
\mathcal{A}_{01}\otimes_R\mathcal{A}_{10}\to \mathcal{A}^1_{00}
\end{equation}
which is a morphism of
$\mathcal{A}^1_0\otimes_K(\mathcal{A}^1_0)^{op}$-modules. Similarly,
we have the pairing
\begin{equation}\label{pairing 10 01}
(\pr^1_{101})^*(\mathcal{A}_{012})
:\mathcal{A}_{10}\otimes_R\mathcal{A}_{01}\to \mathcal{A}^1_{11}
\end{equation}

The pairings \eqref{pairing A A}, \eqref{pairing A M},
\eqref{pairing 01 10} and \eqref{pairing 10 01},
$\mathcal{A}_{ij}\otimes_R\mathcal{A}_{jk}\to\mathcal{A}_{ik}$ are
morphisms of $\mathcal{A}^1_i\otimes_R (\ca^1_k)^{op}$-modules
which, as a consequence of associativity, factor through maps
\begin{equation}
\mathcal{A}_{ij}\otimes_{\mathcal{A}^1_j}\mathcal{A}_{jk}\to
\mathcal{A}_{ik}
\end{equation}
induced by $\ca_{ijk} = (\pr^1_{ijk})^*(\ca_{012})$; here $i, j,
k=0,1$. Define now for every $p \ge 0$ the sheaves $\ca_{ij}^p$ , $0
\le i, j \le p$, on $N_p\cu$ by $\ca_{ij}^p=(\pr^p_{ij})^*
\ca_{01}$. Define also $\ca_{ijk}^p = (\pr^p_{ijk})^*(\ca_{012})$.
We immediately obtain for every $p$ the morphisms
\begin{equation}\label{map Aijk}
\mathcal{A}^p_{ijk}:
\mathcal{A}^p_{ij}\otimes_{\mathcal{A}^p_j}\mathcal{A}^p_{jk}\to
\mathcal{A}^p_{ik}
\end{equation}

\subsubsection{Matrix algebras}
Let $\Mat(\underline{\mathcal{A}})^0 = \mathcal{A}$; thus,
$\Mat(\underline{\mathcal{A}})^0$ is a sheaf of algebras on
$N_0\mathcal{U}$. For $p= 1,2,\ldots$ let
$\Mat(\underline{\mathcal{A}})^p$ denote the sheaf on
$N_p\mathcal{U}$ defined by
\begin{equation*}
\Mat(\underline{\mathcal{A}})^p =
\bigoplus_{i,j=0}^p\mathcal{A}_{ij}^p
\end{equation*}
The maps \eqref{map Aijk} define the pairing
\begin{equation*}
\Mat(\underline{\mathcal{A}})^p\otimes\Mat(\underline{\mathcal{A}})^p
\to \Mat(\underline{\mathcal{A}})^p
\end{equation*}
which endows the sheaf $\Mat(\underline{\mathcal{A}})^p$ with a
structure of an associative algebra by virtue of the associativity
condition. The unit section $1$ is given by $1=\sum_{i=0}^p 1_{ii}$,
where $1_{ii}$ is the image of the unit section of
$\mathcal{A}^p_{ii}$.

\subsubsection{Combinatorial restriction} The algebras
$\Mat(\underline{\mathcal{A}})^p$, $p=0,1,\ldots$, do not form a
cosimplicial sheaf of algebras on $N\mathcal{U}$ in the usual sense.
They are, however, related by \emph{combinatorial restriction} which
we describe presently.

For a morphism $f:[p] \to [q]$ in $\Delta$ define a sheaf on
$N_q\cu$ by
\begin{equation*}
f^{\sharp}\Mat(\underline{\mathcal{A}})^q = \bigoplus_{i,j=0}^p
\mathcal{A}^q_{f(i)f(j)}\ .
\end{equation*}
Note that $f^{\sharp}\Mat(\underline{\mathcal{A}})^q$ inherits a
structure of an algebra.

Recall from the Section \ref{nerve of cover} that the morphism $f$
induces the map $f^*: N_q \cu \to N_p \cu$ and that $f_*$ denotes
the pull-back along $f^*$.  We will also use $f_*$ to denote the
canonical isomorphism of algebras
\begin{equation}\label{combogr}
f_* : f_*\Mat(\underline{\mathcal{A}})^p \to
f^{\sharp}\Mat(\underline{\mathcal{A}})^q
\end{equation}
induced by the isomorphisms
$f_*\mathcal{A}^p_{ij}\cong\mathcal{A}^q_{f(i)f(j)}$.

\subsubsection{Refinement}
Suppose that $\rho:\mathcal{V}\to\mathcal{U}$ is a morphism of
covers. For $p=0,1,\ldots$ there is a natural isomorphism
\begin{equation*}
(N_p\rho)^*\Mat(\underline{\mathcal{A}})^p\cong
\Mat(\underline{\mathcal{A}^\rho})^p
\end{equation*}
of sheaves of algebras on $N_p\mathcal{V}$. The above isomorphisms
are obviously compatible with combinatorial restriction.

\subsection{Local cochains on matrix algebras}
\subsubsection{Local cochains}\label{subsubsection: local cochains}
A \emph{sheaf of matrix algebras} is a sheaf of algebras ${\mathcal B}$ together with a decomposition
\[
{\mathcal B}=\bigoplus_{i,j=0}^p{\mathcal B}_{ij}
\]
as a sheaf of vector spaces which satisfies ${\mathcal B}_{ij}\cdot{\mathcal B}_{jk}={\mathcal B}_{ik}$.

To a matrix algebra $\mathcal B$ one can associate to the DGLA
of \emph{local cochains} defined as follows. For $n=0$ let $C^0({\mathcal B})^{loc}= \bigoplus \mathcal{B}_{ii}
\subset \mathcal{B} = C^0({\mathcal B})$. For $n>0$ let $C^n({\mathcal
B})^{loc}$ denote the subsheaf of $C^n({\mathcal B})$ whose stalks
consist of multilinear maps $D$ such that for any collection of
$s_{i_kj_k}\in{\mathcal B}_{i_kj_k}$
\begin{enumerate}
\item $D(s_{i_1j_1}\otimes\cdots\otimes s_{i_nj_n}) = 0$ unless
$j_k=i_{k+1}$ for all $k = 1,\dots ,n-1$

\item $D(s_{i_0i_1}\otimes s_{i_1i_2}\otimes\cdots\otimes
s_{i_{n-1}i_n})\in{\mathcal B}_{i_0i_n}$
\end{enumerate}
For $I = (i_0,\ldots,i_n)\in [p]^{\times n+1}$ let
\begin{equation*}
C^I(\mathcal{B})^{loc} := \shHom_k(\otimes_{j=0}^{n-1}
\mathcal{B}_{i_j i_{j+1}},\mathcal{B}_{i_0 i_n}) \ .
\end{equation*}
The restriction maps along the embeddings $\otimes_{j=0}^{n-1}
\mathcal{B}_{i_j i_{j+1}}\hookrightarrow\mathcal{B}^{\otimes n}$
induce an isomorphism
$C^n(\mathcal{B})^{loc}\to\oplus_{I\in[p]^{\times n+1}}
C^I(\mathcal{B})^{loc}$.

The sheaf $C^\bullet({\mathcal B})^{loc}[1]$ is a subDGLA of
$C^\bullet({\mathcal B})[1]$ and the inclusion map is a
quasi-isomorphism.

For a matrix algebra $\mathcal B$ on $X$ we denote by
$\Def(\mathcal{B})^{loc}(R)$ the subgroupoid of
$\Def(\mathcal{B})(R)$ with objects $R$-star products which respect
the decomposition given by $(\mathcal{B}\otimes_k R)_{ij} =
\mathcal{B}_{ij}\otimes_k R$ and $1$- and $2$-morphisms defined
accordingly. The composition
\[
\Def(\mathcal{B})^{loc}(R)\to \Def(\mathcal{B})(R) \to
\MC^2(\Gamma(X; C^\bullet(\mathcal{B})[1])\otimes_k \mathfrak{m}_R)
\]
takes values in $\MC^2(\Gamma(X;
C^\bullet(\mathcal{B})^{loc}[1])\otimes_k \mathfrak{m}_R)$ and
establishes an isomorphism of $2$-groupoids
$\Def(\mathcal{B})^{loc}(R)\cong \MC^2(\Gamma(X;
C^\bullet(\mathcal{B})^{loc}[1])\otimes_k \mathfrak{m}_R)$.

\subsubsection{Combinatorial restriction of local
cochains}\label{subsubsection: comb restr}
Suppose given a matrix
algebra ${\mathcal B} = \bigoplus_{i,j=0}^q {\mathcal B}_{ij}$ is a
sheaf of matrix $k$-algebras.

The DGLA $C^\bullet({\mathcal B})^{loc}[1]$ has additional variance
not exhibited by $C^\bullet({\mathcal B})[1]$. Namely, for $f:
[p]\to [q]$ -- a morphism in $\Delta$ -- there is a natural map of
DGLA
\begin{equation}\label{combinatorial restiction of cochains}
f^\sharp : C^\bullet({\mathcal B})^{loc}[1]\to
C^\bullet(f^\sharp{\mathcal B})^{loc}[1]
\end{equation}
defined as follows. Let $f^{ij}_\sharp:
(f^\sharp\mathcal{B})_{ij}\to\mathcal{B}_{f(i)f(j)}$ denote the
tautological isomorphism. For each multi-index $I =
(i_0,\ldots,i_n)\in [p]^{\times n+1}$ let
\begin{equation*}
f^I_\sharp := \otimes_{j=0}^{n-1}f^{i_j i_{j+1}}_\sharp :
\otimes_{j=0}^{n-1} (f^\sharp\mathcal{B})_{i_j
i_{j+1}}\to\otimes_{i=0}^{n-1}\mathcal{B}_{f(i_j) f(i_{j+1})} \ .
\end{equation*}
Let $f^n_\sharp := \oplus_{I\in\SIGMA_p^{\times n+1}}f^I_\sharp$.
The map \eqref{combinatorial restiction of cochains} is defined as
restriction along $f^n_\sharp$.

\begin{lemma}
The map \eqref{combinatorial restiction of cochains} is a morphism
of DGLA
\begin{equation*}
f^\sharp : C^\bullet(\mathcal{B})^{loc}[1]\to
C^\bullet(f^\sharp\mathcal{B})^{loc}[1] \ .
\end{equation*}
\end{lemma}

If follows that combinatorial restriction of local cochains induces
the functor
\[
\MC^2(f^\sharp) : \MC^2(\Gamma(X;
C^\bullet(\mathcal{B})^{loc}[1])\otimes_k \mathfrak{m}_R) \to
\MC^2(\Gamma(X; C^\bullet(f^\sharp\mathcal{B})^{loc}[1])\otimes_k
\mathfrak{m}_R) \ .
\]

Combinatorial restriction with respect to $f$ induces the functor
\[
f^\sharp : \Def(\mathcal{B})^{loc}(R) \to
\Def(f^\sharp\mathcal{B})^{loc}(R) \ .
\]
It is clear that the diagram
\[
\begin{CD}
\Def(\mathcal{B})^{loc}(R) @>{f^\sharp}>>
\Def(f^\sharp\mathcal{B})^{loc}(R) \\
@VVV @VVV \\
\MC^2(\Gamma(X; C^\bullet(\mathcal{B})^{loc}[1])\otimes_k
\mathfrak{m}_R) @>{\MC^2(f^\sharp)}>> \MC^2(\Gamma(X;
C^\bullet(f^\sharp\mathcal{B})^{loc}[1])\otimes_k \mathfrak{m}_R)
\end{CD}
\]

\subsubsection{Cosimplicial DGLA from descent datum }\label{subsubsection: cosimplicial
DGLA from descent datum} Suppose that
$(\mathcal{U},\underline{\mathcal{A}})$ is a descent datum for a
twisted sheaf of algebras as in \ref{definition: descent datum}.
Then, for each $p=0,1,\ldots$  we have the matrix algebra
$\Mat(\underline{\mathcal{A}})^p$ as defined in \ref{subsection:
matrix algebras}, and therefore the DGLA of local cochains
$C^\bullet(\Mat(\underline{\mathcal{A}})^p)^{loc}[1]$ defined in
\ref{subsubsection: local cochains}. For each morphism $f:[p]\to[q]$
there is a morphism of DGLA
\begin{equation*}
f^\sharp : C^\bullet(\Mat(\underline{\mathcal{A}})^q)^{loc}[1]\to
C^\bullet(f^\sharp\Mat(\underline{\mathcal{A}})^q)^{loc}[1]
\end{equation*}
and an isomorphism of DGLA
\begin{equation*}
C^\bullet(f^\sharp\Mat(\underline{\mathcal{A}})^q)^{loc}[1] \cong
f_*C^\bullet(\Mat(\underline{\mathcal{A}})^p)^{loc}[1]
\end{equation*}
induced by the isomorphism $f_* : f_*\Mat(\underline{\mathcal{A}})^p
\to f^{\sharp}\Mat(\underline{\mathcal{A}})^q$ from the equation
\eqref{combogr}. These induce the morphisms of the DGLA of global
sections
\begin{equation*}
\xymatrix@C=-1cm@R=1cm{
\Gamma(N_q\mathcal{U};C^\bullet(\Mat(\underline{\mathcal{A}})^q)^{loc}[1])
\ar[rd]^{f^\sharp} & &
\Gamma(N_p\mathcal{U};C^\bullet(\Mat(\underline{\mathcal{A}})^p)^{loc}[1])
\ar[ld]_{f_*} \\
 &\Gamma(N_q\mathcal{U};f_*C^\bullet(\Mat(
\underline{\mathcal{A}})^p)^{loc}[1]) & }
\end{equation*}

For $\lambda : [n]\to\Delta$ let
\begin{equation*}
\mathfrak{G}(\underline{\mathcal{A}})^\lambda =
\Gamma(N_{\lambda(n)}\mathcal{U}; \lambda(0n)_*
C^\bullet(\Mat(\underline{\mathcal{A}})^{\lambda(0)})^{loc}[1])
\end{equation*}

Suppose given another simplex $\mu : [m]\to\Delta$ and morphism
$\phi : [m]\to[n]$ such that $\mu = \lambda\circ\phi$ (i.e. $\phi$
is a morphism of simplices $\mu\to\lambda$). The morphism $(0n)$
factors uniquely into $0\to\phi(0)\to\phi(m)\to n$, which, under
$\lambda$, gives the factorization of
$\lambda(0n):\lambda(0)\to\lambda(n)$ (in $\Delta$) into
\begin{equation}\label{defh}
\begin{CD}
\lambda(0) @>f>> \mu(0) @>g>> \mu(m) @>h>> \lambda(n) \ ,
\end{CD}
\end{equation}
where $g=\mu(0m)$. The map
\begin{equation*}
\phi_*: \mathfrak{G}(\underline{\mathcal{A}})^\mu\to
\mathfrak{G}(\underline{\mathcal{A}})^\lambda
\end{equation*}
is the composition
\begin{multline*}
\Gamma(N_{\mu(m)}\mathcal{U}; g_*
C^\bullet(\Mat(\underline{\mathcal{A}})^{\mu(0)})^{loc}[1])
\xrightarrow{h_*} \\
\Gamma(N_{\lambda(n)}\mathcal{U}; h_*g_*
C^\bullet(\Mat(\underline{\mathcal{A}})^{\mu(0)})^{loc}[1])
\xrightarrow{f^\sharp} \\
\Gamma(N_{\lambda(n)}\mathcal{U}; h_*g_*f_*
C^\bullet(\Mat(\underline{\mathcal{A}})^{\lambda(0)})^{loc}[1])
\end{multline*}

Suppose given yet another simplex, $\nu : [l]\to\Delta$ and a
morphism of simplices $\psi : \nu\to\mu$, i.e. a morphism $\psi :
[l]\to[m]$ such that $\nu = \mu\circ\psi$. Then, the composition
$\phi_*\circ\psi_* : \mathfrak{G}(\underline{\mathcal{A}})^\nu\to
\mathfrak{G}(\underline{\mathcal{A}})^\lambda$ coincides with the
map $(\phi\circ\psi)_*$.

For $n=0,1,2,\ldots$ let
\begin{equation}\label{def complex cosimp matrix alg}
\mathfrak{G}(\uca)^n = \prod_{[n]\xrightarrow{\lambda}\Delta}
\mathfrak{G}(\uca)^\lambda
\end{equation}
A morphism $\phi:[m]\to[n]$ in $\Delta$ induces the map of DGLA
$\phi_*: \mathfrak{G}(\uca)^m\to\mathfrak{G}(\uca)^n$. The
assignment $\Delta\ni [n]\mapsto \mathfrak{G}(\uca)^n$,
$\phi\mapsto\phi_*$ defines the cosimplicial DGLA denoted
$\mathfrak{G}(\uca)$.

\subsection{Acyclicity}
\begin{thm}\label{acyclicity of matrix dgla}
The cosimplicial DGLA $\mathfrak{G}(\underline{\mathcal{A}})$ is
acyclic, i.e. it satisfies the condition \eqref{acyclicity
condition}.
\end{thm}

The rest of the section is devoted to the proof of the Theorem
\ref{acyclicity of matrix dgla}. We fix a degree of Hochschild
cochains $k$.

For $\lambda : [n]\to\Delta$ let $\mathfrak{c}^\lambda =
\Gamma(N_{\lambda(n)}\mathcal{U};
\lambda(0n)_*C^k(\Mat(\uca)^{\lambda(0)})^{loc})$. For a morphism
$\phi: \mu\to\lambda$ we have the map $\phi_* :
\mathfrak{c}^\mu\to\mathfrak{c}^\lambda$ defined as in
\ref{subsubsection: cosimplicial DGLA from descent datum}.

Let $(\mathfrak{C}^\bullet,\partial)$ denote the corresponding
cochain complex whose definition we recall below. For $n=0,1,\ldots$
let $\mathfrak{C}^n = \prod \limits_{[n]\xrightarrow{\lambda}\Delta}
\mathfrak{c}^\lambda$. The differential $\partial^n :
\mathfrak{C}^n\to\mathfrak{C}^{n+1}$ is defined by the formula
$\partial^n = \sum_{i=0}^{n+1} (-1)^i(\partial^n_i)_*$.

\subsubsection{Decomposition of local cochains}
As was noted in \ref{subsubsection: local cochains}, for $n,q =
0,1,\ldots$ there is a direct sum decomposition
\begin{equation}\label{I-decomposition}
C^k(\Mat(\underline{\mathcal{A}})^q)^{loc} =
\bigoplus_{I\in[q]^{k+1}} C^I(\Mat(\underline{\mathcal{A}})^q)^{loc}
.
\end{equation}

In what follows we will interpret a multi-index
$I=(i_0,\ldots,i_n)\in[q]^{k+1}$ as a map $I :
\{0,\ldots,k\}\to[q]$. For $I$ as above let $s(I) = \vert\im(I)\vert
- 1$. The map $I$ factors uniquely into the composition
\begin{equation*}
\{0,\ldots,k\}\xrightarrow{I^\prime}[s(I)]\xrightarrow{m(I)}[q]
\end{equation*}
where the second map is a morphism in $\Delta$ (i.e. is order
preserving). Then, the isomorphisms
$m(I)_*\mathcal{A}_{I^\prime(i)I^\prime(j)}\cong\mathcal{A}_{I(i)I(j)}$
induce the isomorphism
\begin{equation*}
 m(I)_*C^{I^\prime}(\Mat(\underline{\mathcal{A}})^{s(I)})^{loc} \to C^I(\Mat(\underline{\mathcal{A}})^q)^{loc}
\end{equation*}
Therefore, the decomposition \eqref{I-decomposition} may be
rewritten as follows:
\begin{equation}\label{improved I-decomposition}
C^k(\Mat(\underline{\mathcal{A}})^q)^{loc} = \bigoplus_e \bigoplus_I
e_*C^I(\Mat(\underline{\mathcal{A}})^p)^{loc}
\end{equation}
where the summation is over \emph{injective} (monotone) maps
$e:[s(e)]\to[q]$ and \emph{surjective} maps
$I:\{0,\ldots,k\}\to[s(e)]$. Note that, for $e$, $I$ as above, there
is an isomorphism
\begin{equation*}
e_*C^I(\Mat(\underline{\mathcal{A}})^p)^{loc}\cong C^{e\circ
I}(e^\sharp \Mat(\underline{\mathcal{A}})^{s(e)})^{loc}
\end{equation*}

\subsubsection{Filtrations}\label{subsubsection: filtrations}
Let $F^\bullet C^k(\Mat(\underline{\mathcal{A}})^q)^{loc}$ denote
the \emph{decreasing} filtration defined by
\begin{equation*}
F^s C^k(\Mat(\underline{\mathcal{A}})^q)^{loc} = \bigoplus_{e:
s(e)\ge s} \bigoplus_I e_*C^I(\Mat(\underline{\mathcal{A}})^p)^{loc}
\end{equation*}
Note that $F^s C^k(\Mat(\underline{\mathcal{A}})^q)^{loc} = 0$ for
$s>k$ and $Gr^s C^k(\Mat(\underline{\mathcal{A}})^q)^{loc} = 0$ for
$s<0$.

The filtration $F^\bullet
C^k(\Mat(\underline{\mathcal{A}})^{\lambda(0)})^{loc}$ induces the
filtration $F^\bullet\mathfrak{c}^\lambda$, hence the filtration
$F^\bullet\mathfrak{C}^\bullet$ with $F^s\mathfrak{c}^\lambda =
\Gamma(N_{\lambda(n)}\mathcal{U}; \lambda(0n)_*F^s
C^k(\Mat(\underline{A})^{\lambda(0)})^{loc})$, $F^s\mathfrak{C}^n =
\prod
\limits_{[n]\xrightarrow{\lambda}\Delta}F^s\mathfrak{c}^\lambda$.

The following result then is an easy consequence of the definitions:
\begin{lemma}
For $\phi : \mu\to\lambda$ the induced map $\phi_*:
\mathfrak{c}^\mu\to \mathfrak{c}^\lambda$ preserves filtration.
\end{lemma}

\begin{cor}
The differential $\partial^n$ preserves the filtration.
\end{cor}

\begin{prop}\label{acyclicity of Gr}
For any $s$ the complex $Gr^s(\mathfrak{C}^\bullet,\partial)$ is
acyclic in non-zero degrees.
\end{prop}
\begin{proof}
Use the following notation: for a simplex $\mu : [m]\to\Delta$ and
an arrow $e : [s]\to[\mu(0)]$ the simplex $\mu^e : [m+1]\to\Delta$
is defined by $\mu^e(0)= [s]$, $\mu^e(01)=e$, and
$\mu^e(i)=\mu(i-1)$, $\mu^e(i,j) = \mu(i-1,j-1)$ for $i>0$.

For $D=(D_\lambda)\in F^s\mathfrak{C}^n$,
$D_\lambda\in\mathfrak{c}^\lambda$ and $\mu : [n-1]\to\Delta$ let
\begin{equation*}
h^n(D)_\mu = \sum_{I,e}e_*D_{\mu^e}^I
\end{equation*}
where $e: [s]\to\mu(0)$ is an injective map,
$I:\{0,\ldots,k\}\to[s]$ is a surjection, and $D_{\mu^e}^I$ is the
$I$-component of $D_{\mu^e}$ in the sense of the decomposition
\eqref{I-decomposition}. Put $h^n(D) =
(h^n(D)_\mu)\in\mathfrak{C}^{n-1}$. It is clear that $h^n(D)\in
F^s\mathfrak{C}^{n-1}$. Thus the assignment $D\mapsto h^n(D)$
defines a filtered map $h^n : \mathfrak{C}^n\to\mathfrak{C}^{n-1}$
for all $n\geq 1$. Hence, for any $s$ we have the map $Gr^s h^n :
Gr^s\mathfrak{C}^n\to Gr^s\mathfrak{C}^{n-1}$.

Next, we calculate the effect of the face maps $\partial^{n-1}_i :
[n-1]\to[n]$. First of all, note that, for $1\leq i\leq n-1$,
$(\partial^{n-1}_i)_* = \id$. The effect of $(\partial^{n-1}_0)_*$
is the map on global sections induced by
\begin{equation*}
\lambda(01)^\sharp :
\lambda(1n)_*C^k(\Mat(\underline{\mathcal{A}})^{\lambda(1)})^{loc}
\to \lambda(0n)_*C^k(\Mat(\underline{A})^{\lambda(0)})^{loc}
\end{equation*}
and $(\partial^{n-1}_n)_* = \lambda(n-1,n)_*$. Thus, for
$D=(D_\lambda)\in F^s\mathfrak{C}^n$, we have:
\begin{equation*}
((\partial^{n-1}_i)_* h^n(D))_\lambda =
(\partial^{n-1}_i)_*\sum_{I,e}
e_*D_{(\lambda\circ\partial^{n-1}_i)^e}^I = \left\lbrace
\begin{array}{ll}
 \lambda(01)^\sharp\sum_{I,e} e_*D_{(\lambda\circ\partial^{n-1}_0)^e}^I & \text{if $i = 0$} \\
\sum_{I,e} e_*D_{(\lambda\circ\partial^{n-1}_i)^e}^I  & \text{if $1\leq i\leq n-1$} \\
\lambda(n-1,n)_*\sum_{I,e} e_*D_{(\lambda\circ\partial^{n-1}_n)^e}^I
& \text{if $i=n$}      \end{array}\right.
\end{equation*}

On the other hand,
\begin{equation*}
h^{n+1}((\partial^{n}_j)_*D)_\lambda = \sum_{I,e}
e_*((\partial^{n}_j)_*D)^I_{\lambda^e} = \left\lbrace
\begin{array}{ll}
\sum_{I,e} e_*e^\sharp D^I_\lambda & \text{if $j = 0$} \\
\sum_{I,e} e_* D^I_{(\lambda\circ\partial^{n-1}_0)^{\lambda(01)\circ e}} & \text{if $j = 1$} \\
\sum_{I,e} e_*D_{(\lambda\circ\partial^{n-1}_{j-1})^e}^I & \text{if $2\leq j\leq n$} \\
\sum_{I,e}\lambda(n-1,n)_* e_*D_{(\lambda\circ\partial^{n-1}_n)^e}^I
& \text{if $j=n+1$}     \end{array}\right.
\end{equation*}

Note that $\sum_{I,e} e_*e^\sharp D^I_\lambda = D_\lambda \mod
F^{s+1}\mathfrak{c}^\lambda$, i.e. the map
$h^{n+1}\circ\partial^n_0$ induces the identity map on
$Gr^s\mathfrak{c}^\lambda$ for all $\lambda$, hence the identity map
on $Gr^s\mathfrak{C}^n$.

The identities $(\partial^{n-1}_i h^n(D))_\lambda =
h^{n+1}((\partial^{n}_{i+1})_*D)_\lambda \mod
F^s\mathfrak{c}^\lambda$ hold for $0\leq i\leq n$.

For $n=0,1,\ldots$ let $\partial^n = \sum_i (-1)\partial^n_i$. The
above identities show that, for $n>0$, $h^{n+1}\circ\partial^n +
\partial^{n-1}\circ h^n = \id$.
\end{proof}

 \begin{proof}[Proof of Theorem \ref{acyclicity of matrix dgla}]
 Consider for a fixed  $k$ the complex of Hochschild degree $k$ cochains
 $\mathfrak{C}^\bullet$.
It is shown  in Proposition \ref{acyclicity of Gr} that this complex
admits a finite filtration $F^\bullet\mathfrak{C}^\bullet$ such that
$Gr\mathfrak{C}^\bullet$ is acyclic in positive degrees. Therefore
$\mathfrak{C}^\bullet$is also acyclic in positive degree. As this
holds for every $k$, the condition \eqref{acyclicity condition} is
satisfied.
\end{proof}

\section{Deformations of algebroid stacks}\label{defalgstack}
In this section we define a $2$-groupoid  of  deformations of an
algebroid stack. We also define $2$-groupoids  of  deformations and
star products of a descent datum and relate it with the $2$-groupoid
$\mathfrak{G}$-stacks of an appropriate cosimplicial DGLA.

\subsection{Deformations of linear stacks}
\begin{definition}
Let   $\mathcal{B}$ be a prestack on $X$ in $R$-linear
categories. We say that  $\mathcal{B}$ is  \emph{flat} if for
any $U\subseteq X$, $A,B\in\mathcal{B}(U)$ the sheaf
$\shHom_\mathcal{B}(A,B)$ is flat (as a sheaf of $R$-modules).
\end{definition}
Suppose now that $\mathcal{C}$ is a stack in $k$-linear
categories on $X$ and $R$ is a commutative Artin $k$-algebra.
We denote by $\Def(\mathcal{C})(R)$ the $2$-category with
\begin{itemize}
\item objects: pairs $(\mathcal{B}, \varpi)$, where $\mathcal{B}$ is a stack in
$R$-linear categories flat over $R$ and $\varpi :
\widetilde{\mathcal{B}\otimes_R k} \to \mathcal{C}$ is an
equivalence of stacks in $k$-linear categories

\item $1$-morphisms: a $1$-morphism $(\mathcal{B}^{(1)}, \varpi^{(1)})\to
(\mathcal{B}^{(2)}, \varpi^{(2)})$ is a pair $(F,\theta)$ where $F :
\mathcal{B}^{(1)}\to \mathcal{B}^{(2)}$ is a $R$-linear functor and
$\theta : \varpi^{(2)}\circ (F\otimes_R k) \to \varpi^{(1)}$ is an
isomorphism of functors

\item $2$-morphisms: a $2$-morphism $(F',\theta')\to
(F'',\theta'')$ is a morphism of $R$-linear functors $\kappa : F'\to
F''$ such that
$\theta''\circ(\id_{\varpi^{(2)}}\otimes(\kappa\otimes_R k)) =
\theta'$
\end{itemize}

The $2$-category $\Def(\mathcal{C})(R)$ is a $2$-groupoid.

\begin{lemma}\label{lemma: def of algd is algd}
Suppose that $\mathcal{B}$ is a flat $R$-linear stack on $X$ such
that $\widetilde{\mathcal{B}\otimes_R k}$ is an algebroid stack.
Then, $\mathcal{B}$ is an algebroid stack.
\end{lemma}

\begin{proof}
Let $x\in X$. Suppose that for any neighborhood $U$ of $x$ the
category $\mathcal{B}(U)$ is empty. Then, the same is true
about $\widetilde{\mathcal{B}\otimes_R k}(U)$ which contradicts
the assumption that $\widetilde{\mathcal{B}\otimes_R k}$ is an
algebroid stack. Therefore, $\mathcal{B}$ is locally nonempty.

Suppose that $U$ is an open subset and $A,B$ are two objects in
$\mathcal{B}(U)$. Let $\overline{A}$ and $\overline{B}$ be
their respective images in $\widetilde{\mathcal{B}\otimes_R
k}(U)$. We have: $\Hom_{\widetilde{\mathcal{B}\otimes_R
k}(U)}(\overline{A},\overline{B}) =
\Gamma(U;\shHom_{\mathcal{B}}(A,B)\otimes_R k)$. Replacing $U$
by a subset we may assume that there is an isomorphism
$\overline{\phi} : \overline{A}\to\overline{B}$.

The short exact sequence
\[
0 \to \mathfrak{m}_R \to R \to k \to 0
\]
gives rise to the sequence
\[
0 \to \shHom_{\mathcal{B}}(A,B)\otimes_R\mathfrak{m}_R \to \shHom_{\mathcal{B}}(A,B) \to \shHom_{\mathcal{B}}(A,B)\otimes_R k \to 0
\]
of sheaves on $U$ which is exact due to flatness of
$\shHom_{\mathcal{B}}(A,B)$. The surjectivity of the map
$\shHom_{\mathcal{B}}(A,B) \to
\shHom_{\mathcal{B}}(A,B)\otimes_R k$ implies that for any
$x\in U$ there exists a neighborhood $x\in V\subseteq U$ and
$\phi\in\Gamma(V;\shHom_{\mathcal{B}}(A\vert_V,B\vert_V)) =
\Hom_{\mathcal{B}(V)}(A\vert_V,B\vert_V)$ such that
$\overline{\phi}\vert_V$ is the image of $\phi$. Since
$\overline{\phi}$ is an isomorphism and $\mathfrak{m}_R$ is
nilpotent it follows that $\phi$ is an isomorphism.
\end{proof}

\subsection{Deformations of descent data}\label{defdescdata}
Suppose that $(\cu,\uca)$ is an $k$-descent datum and $R$ is a
commutative Artin $k$-algebra.

We denote by $\Def'(\cu,\uca)(R)$ the $2$-category with
\begin{itemize}
\item objects: \emph{$R$-deformations of $(\cu,\uca)$}; such a
gadget is a flat $R$-descent datum $(\cu,\underline{\mathcal{B}})$
together with an isomorphism of $k$-descent data $\pi :
(\cu,\underline{\mathcal{B}}\otimes_R k))\to (\cu,\uca)$

\item $1$-morphisms: a $1$-morphism of deformations $(\cu,
\underline{\mathcal{B}}^{(1)}, \pi^{(1)})\to (\cu,
\underline{\mathcal{B}}^{(2)}, \pi^{(2)})$ is a pair

$(\underline{\phi},\underline{a})$, where $\underline{\phi}: (\cu,
\underline{\mathcal{B}}^{(1)})\to (\cu,
\underline{\mathcal{B}}^{(2)})$ is a $1$-morphism of $R$-descent
data and $\underline{a} : \pi^{(2)}\circ (\underline{\phi}\otimes_R
k)\to \pi^{(1)}$ is $2$-isomorphism

\item $2$-morphisms: a $2$-morphism $(\underline{\phi}',
\underline{a}')\to (\underline{\phi}'', \underline{a}'')$ is a
$2$-morphism $\underline{b} : \underline{\phi}'\to
\underline{\phi}''$ such that
$\underline{a}''\circ(\id_{\pi^{(2)}}\otimes(\underline{b}\otimes_R
k)) = \underline{a}'$
\end{itemize}

Suppose that $(\underline{\phi},\underline{a})$ is a $1$-morphism.
It is immediate from the definition above that the morphism of
$k$-descent data $\underline{\phi}\otimes_R k$ is an isomorphism.
Since $R$ is an extension of $k$ by a nilpotent ideal the morphism
$\underline{\phi}$ is an isomorphism. Similarly, any $2$-morphism is
an isomorphism, i.e. $\Def'(\cu,\uca)(R)$ is a $2$-groupoid.

The assignment $R\mapsto \Def'(\cu,\underline{\mathcal{A}})(R)$ is
fibered $2$-category in $2$-groupoids over the category of
commutative Artin $k$-algebras (\cite{Breen}).

\subsubsection{Star products}
A \emph{($R$-)star product} on $(\cu,\uca)$ is a deformation
$(\cu,\underline{\mathcal{B}},\pi)$ of $(\cu,\uca)$ such that
$\mathcal{B}_{01} = \mathcal{A}_{01}\otimes_k R$ and $\pi :
\mathcal{B}_{01}\otimes_R k\to \mathcal{A}_{01}$ is the canonical
isomorphism. In other words, a star product is a structure of an
$R$-descent datum on $(\cu,\underline{\mathcal{A}}\otimes_k R)$ such
that the natural map
\begin{equation*}
(\cu,\underline{\mathcal{A}}\otimes_k R)\to
(\cu,\underline{\mathcal{A}})
\end{equation*}
is a morphism of such.

We denote by $\Def(\cu,\uca)(R)$ the full $2$-subcategory of
$\Def'(\cu,\uca)(R)$ whose objects are star products.

The assignment $R\mapsto \Def(\cu,\underline{\mathcal{A}})(R)$ is
fibered $2$-category in $2$-groupoids over the category of
commutative Artin $k$-algebras (\cite{Breen}) and the inclusions
$\Def(\cu,\uca)(R)\to \Def'(\cu,\uca)(R)$ extend to a morphism of
fibered $2$-categories.

\begin{prop}\label{star products are equivalent to def}
Suppose that $(\cu,\uca)$ is a $\mathbb{C}$-linear descent datum
with $\mathcal{A}=\mathcal{O}_{N_0\cu}$. Then, the embedding
$\Def(\cu,\uca)(R)\to \Def'(\cu,\uca)(R)$ is an equivalence.
\end{prop}

\subsubsection{Deformations and  $\mathfrak{G}$-stacks}
Suppose that $(\cu,\uca)$ is a $k$-descent datum and
$(\cu,\underline{\mathcal{B}})$ is an $R$-star product on
$(\cu,\uca)$. Then, for every $p=0,1,\ldots$ the matrix algebra
$\Mat(\underline{\mathcal{B}})^p$ is a flat $R$-deformation of the
matrix algebra $\Mat(\uca)^p$. The identification
$\mathcal{B}_{01}=\ca_{01}\otimes_k R$ gives rise to the
identification $\Mat(\underline{\mathcal{B}})^p =
\Mat(\uca)^p\otimes_k R$ of the underlying sheaves of $R$-modules.
Using this identification we obtain the Maurer-Cartan element
$\mu^p\in \Gamma(N_p\cu; C^2(\Mat(\uca)^p)^{loc}\otimes_k
\mathfrak{m}_R)$. Moreover, the equation \eqref{combogr} implies
that for a morphism $f: [p]\to [q]$ in $\Delta$ we have $f_*\mu^p =
f^{\sharp} \mu^q$. Therefore the collection $\mu^p$ defines an
element in $\Stack_{str} (\mathfrak{G}(\uca)\otimes_k
\mathfrak{m}_R)_0$. The considerations of \ref{subsubsection: local
cochains} and \ref{subsubsection: comb restr} imply that this
construction extends to an isomorphism of $2$-groupoids
\begin{equation}\label{def isom StackStr}
\Def(\cu,\uca)(R)\to \Stack_{str} (\mathfrak{G}(\uca)\otimes_k
\mathfrak{m}_R)
\end{equation}

Combining \eqref{def isom StackStr} with the embedding
\begin{equation}\label{StackStr into Stack}
\Stack_{str} (\mathfrak{G}(\uca)\otimes_k \mathfrak{m}_R) \to \Stack
(\mathfrak{G}(\uca)\otimes_k \mathfrak{m}_R)
\end{equation}
we obtain the functor
\begin{equation}\label{def to stack R}
\Def(\cu,\uca)(R)\to \Stack (\mathfrak{G}(\uca)\otimes_k
\mathfrak{m}_R)
\end{equation}
The naturality properties of \eqref{def to stack R} with respect to
base change imply that \eqref{def to stack R} extends to morphism of
functors on the category of commutative Artin algebras.

Combining this with the results of Theorems \ref{acyclicity of
matrix dgla} and \ref{acyclic vs strictness} implies the following:
\begin{prop}\label{allst}
The functor \eqref{def to stack R} is an
equivalence.
\end{prop}
\begin{proof}
By Theorem \ref{acyclicity of matrix dgla} the DGLA
$\mathfrak{G}(\uca)\otimes_k \mathfrak{m}_R$ satisfies the
assumptions of Theorem \ref{acyclic vs strictness}. The latter says
that the inclusion \eqref{StackStr into Stack} is an equivalence.
Since \eqref{def isom StackStr} is an isomorphism, the composition
\eqref{def to stack R} is an equivalence as claimed.
\end{proof}

\section{Jets}\label{jets}
In this section we use constructions involving the infinite jets to
simplify the cosimplicial DGLA governing the deformations of a
descent datum.
\subsection{Infinite jets of a vector bundle}\label{remind jet}
Let $M$ be a smooth manifold, and ${\mathcal E}$ a locally-free
${\mathcal O}_{M}$-module of finite rank.

Let $\pi_i:M\times M\to M$, $i=1,2$ denote the projection on the
$i^{\text{th}}$. Denote by  $\Delta_M : M\to M\times M$  the
diagonal embedding and let $\Delta_M^*: {\mathcal O}_{M\times
M}\to{\mathcal O}_M$ be the induced map. Let  ${\mathcal
I}_{M}:=\ker(\Delta_{M}^*)$.

Let
\[
{\mathcal J}^k({\mathcal E}):=(\pi_1)_*\left({\mathcal O}_{M\times
M}/{\mathcal I}_{M}^{k+1}\otimes_{\pi_2^{-1}{\mathcal
O}_{M}}\pi_2^{-1}{\mathcal E}\right) \ ,
\]
${\mathcal J}^k_{M}:= {\mathcal J}^k({\mathcal O}_{M})$. It is clear
from the above definition that ${\mathcal J}^k_{M}$ is, in a natural
way, a sheaf of commutative algebras and ${\mathcal J}^k({\mathcal
E})$ is a sheaf of ${\mathcal J}^k_{M}$-modules. If moreover
$\mathcal{E}$ is a sheaf of algebras, ${\mathcal J}^k({\mathcal E})$
will canonically be a sheaf of algebras as well. We regard
${\mathcal J}^k({\mathcal E})$ as ${\mathcal O}_{M}$-modules via the
pull-back map $\pi_1^*:{\mathcal O}_{M}\to(\pi_1)_*{\mathcal
O}_{{M\times M}}$

For $0\leq k\leq l$ the inclusion ${\mathcal
I}_{M}^{l+1}\to{\mathcal I}_{M}^{k+1}$ induces the surjective map
${\mathcal J}^l({\mathcal E})\to{\mathcal J}^k({\mathcal E})$. The
sheaves ${\mathcal J}^k({\mathcal E})$, $k=0,1,\ldots$ together with
the maps just defined    form an inverse system. Define ${\mathcal
J}({\mathcal E}):=\underset{\longleftarrow}{\lim}{\mathcal
J}^k({\mathcal E})$. Thus, ${\mathcal J}({\mathcal E})$ carries a
natural topology.

We denote by $p_{\mathcal{E}}: \cj(\mathcal{E}) \to \mathcal{E}$ the
canonical projection. In the case when $\mathcal{E} =\co_M$ we
denote by $p$ the corresponding projection $p:\cj_M \to
\mathcal{O}_M$. By $j^k: {\mathcal E}\to{\mathcal J}^k({\mathcal
E})$ we denote the map $e\mapsto 1\otimes e$, and
$j^\infty:=\underset{\longleftarrow}{\lim}j^k$. In the case
$\mathcal{E} = \co_M$ we also have the canonical  embedding $\co_M
\to \cj_M$ given by $f \mapsto f \cdot j^{\infty}(1)$.

Let
\[
d_1 : {\mathcal O}_{{M\times M}}\otimes_{\pi_2^{*}{\mathcal
O}_{M}}\pi_2^{-1}{\mathcal E}  \to
\pi_1^{-1}\Omega^1_{M}\otimes_{\pi_1^{-1}{\mathcal O}_{M}}{\mathcal
O}_{{M\times M}}\otimes_{\pi_2^{-1}{\mathcal
O}_{M}}\pi_2^{-1}{\mathcal E}
\]
denote the exterior derivative along the first factor. It satisfies
\[
d_1({\mathcal I}_{M}^{k+1}\otimes_{\pi_2^{-1}{\mathcal
O}_{M}}\pi_2^{-1}{\mathcal E})\subset
\pi_1^{-1}\Omega^1_M\otimes_{\pi_1^{-1}{\mathcal O}_{M}}{\mathcal
I}_{M}^k\otimes_{\pi_2^{-1}{\mathcal O}_{M}}\pi_2^{-1}{\mathcal E}
\]
for each $k$ and, therefore, induces the map
\[
d_1^{(k)} : {\mathcal J}^k({\mathcal E})\to\Omega^1_{M/{\mathcal
P}}\otimes_{{\mathcal O}_{M}}{\mathcal J}^{k-1}({\mathcal E})
\]
The maps $d_1^{(k)}$ for different values of $k$ are compatible with
the maps ${\mathcal J}^l({\mathcal E})\to{\mathcal J}^k({\mathcal
E})$ giving rise to the \emph{canonical flat connection}
\[
\nabla^{can}  : {\mathcal J}({\mathcal E})\to\Omega^1_{M}\otimes
{\mathcal J}({\mathcal E})
\]

Here and below   we use notation $(\bullet)\otimes {\mathcal J}
({\mathcal E})$ for
$\underset{\longleftarrow}{\lim}(\bullet)\otimes_{{\mathcal
O}_{M}}{\mathcal J}^k({\mathcal E})$.

Since $\nabla^{can}$ is flat we obtain the complex of sheaves
$\DR(\cj(\mathcal{E}))=(\Omega^{\bullet}_M \otimes \cj(\mathcal{E}),
\nabla^{can})$. When $\mathcal{E}=\co_M$ embedding $\co_M \to \cj_M$
induces embedding of de Rham complex $\DR(\co)=(\Omega^{\bullet}_M,
d)$ into $\DR(\cj)$. We denote the quotient by $\DR(\cj/\co)$. All
the complexes above are complexes of soft sheaves.
 We have the following:
\begin{prop}
The (hyper)cohomology $H^i(M, \DR(\cj(\mathcal{E}))) \cong
H^i(\Gamma(M; \Omega_M\otimes\cj(\mathcal{E})), \nabla^{can})$ is
$0$ if $i>0$. The map $j^{\infty}:\mathcal{E} \to \cj(\mathcal{E})$
induces the isomorphism between $\Gamma(\mathcal{E})$ and $H^0(M,
\DR(\cj(\mathcal{E}))) \cong H^0(\Gamma(M;
\Omega_M\otimes\cj(\mathcal{E})), \nabla^{can})$
\end{prop}

\subsection{Jets of line bundles}
Let, as before, $M$ be a smooth manifold, $\cj_M$ be the sheaf of
infinite jets of smooth functions on $M$ and $p: \cj_M \to \co_M$ be
the canonical projection. Set $\cj_{M,0}= \ker p$. Note that
$\cj_{M,0}$ is a sheaf of $\co_M$ modules and therefore is soft.

Suppose now that $\cl$ is a line bundle on $M$. Let $\shIsom_0(\cl
\otimes{\mathcal J}_M,{\mathcal J}(\cl))$ denote the sheaf of local
${\mathcal J}_M$-module isomorphisms $ \cl \otimes{\mathcal
J}_M\to{\mathcal J}(\cl)$ such that the following diagram is
commutative:
\begin{equation*}
\xymatrix{ & \cl \otimes \cj_M \ar[rr]\ar[dr]_{\id \otimes p} &&\cj(\cl)\ar[dl]^{p_{\cl}}\\
           &&\cl &
          }
\end{equation*}

It is easy to see that the canonical map $\cj_M\to
\shEnd_{\cj_M}(\cl \otimes \cj_M)$ is an isomorphism. For
$\phi\in\cj_{M,0}$ the exponential series $\exp(\phi)$ converges.
The composition
\[
\cj_{M,0}\xrightarrow{\exp}\cj_M\to \shEnd_{\cj_M}(\cl \otimes
\cj_M)
\]
defines an isomorphism of sheaves of groups
\[
\exp : \cj_{M,0}\to \shAut_0(\cl \otimes \cj_M) \ ,
\]
where $\shAut_0(\cl \otimes \cj_M)$ is the sheaf of groups of
(locally defined) $\cj_M$-linear automorphisms of $\cl \otimes
\cj_M$ making the diagram
\begin{equation*}
\xymatrix{ & \cl \otimes \cj_M \ar[rr]\ar[rd]_{\id \otimes p} &&\cl \otimes \cj_M \ar[ld]^{\id \otimes p}\\
           &&\cl &
          }
\end{equation*}
commutative.

\begin{lemma}\label{lemma:isom is a torsor}
The sheaf $\shIsom_0({\mathcal L}\otimes{\mathcal J}_M,{\mathcal
J}({\mathcal L}))$ is a torsor under the sheaf of groups
$\exp\cj_{M,0}$.
\end{lemma}
\begin{proof}
Since $\cl$ is locally trivial, both ${\mathcal J}({\mathcal L})$
and ${\mathcal L}\otimes{\mathcal J}_M$ are locally isomorphic to
${\mathcal J}_M$. Therefore the sheaf $\shIsom_0({\mathcal
L}\otimes{\mathcal J}_M,{\mathcal J}({\mathcal L}))$ is locally
non-empty, hence a torsor.
\end{proof}

\begin{cor}\label{trivial torsor}
The torsor $\shIsom_0({\mathcal L}\otimes{\mathcal J}_M,{\mathcal
J}({\mathcal L}))$ is trivial, i.e. $\Isom_0({\mathcal
L}\otimes{\mathcal J}_M,{\mathcal J}({\mathcal
L})):=\Gamma(M;\shIsom_0({\mathcal L}\otimes{\mathcal J}_M,{\mathcal
J}({\mathcal L})))\ne \varnothing$.
\end{cor}
\begin{proof}
Since the sheaf of groups $\cj_{M,0}$ is soft we have $H^1(M,
\cj_{M,0})= 0$ (\cite{DD}, Lemme 22, cf. also \cite{Brylinski},
Proposition 4.1.7). Therefore every $\cj_{M,0}$-torsor is trivial.
\end{proof}

\begin{cor}\label{aff}
The set $\Isom_0({\mathcal L}\otimes{\mathcal J}_M,{\mathcal
J}({\mathcal L}))$ is an affine space with the underlying vector
space $\Gamma(M;\cj_{M,0})$.
\end{cor}

Let $\cl_1$ and $\cl_2$ be two line bundles, and $f:\cl_1 \to \cl_2$
an isomorphism.
Then $f$ induces a map $\Isom_0({\mathcal L_2}\otimes{\mathcal
J}_M,{\mathcal J}({\mathcal L_2})) \to \Isom_0({\mathcal
L_1}\otimes{\mathcal J}_M,{\mathcal J}({\mathcal L_1}))$ which we
denote by $\Ad f$:
\begin{equation*}
\Ad f (\sigma) = (j^{\infty}(f))^{-1}\circ \sigma\circ(f\otimes
\id).
\end{equation*}

Let $\cl$ be a line bundle on $M$ and $f: N \to M$ is a smooth map.
Then there is a pull-back  map $f^*: \Isom_0({\mathcal
L}\otimes{\mathcal J}_M,{\mathcal J}({\mathcal L})) \to
\Isom_0(f^*{\mathcal L}\otimes{\mathcal J}_N,{\mathcal
J}(f^*{\mathcal L}))$.

If $\cl_1$, $\cl_2$ are two line bundles, and $\sigma_i\in
\Isom_0({\mathcal L_i}\otimes{\mathcal J}_M,{\mathcal J}({\mathcal
L_i}))$, $i=1, 2$. Then we denote by $\sigma_1 \otimes \sigma_2$ the
induced element of $\Isom_0(({\mathcal L_1}\otimes
\cl_2)\otimes{\mathcal J}_M,{\mathcal J}(\cl_1\otimes{\mathcal
L_2}))$.

For a line bundle $\cl$ let $\cl^*$ be its dual. Then given $\sigma
\in \Isom_0({\mathcal L}\otimes{\mathcal J}_M,{\mathcal J}({\mathcal
L}))$ there exists a unique a unique $\sigma^* \in \Isom_0({\mathcal
L}^*\otimes{\mathcal J}_M,{\mathcal J}({\mathcal L}^*))$ such that
$\sigma^* \otimes \sigma = \id$.

For any bundle $E$ $\cj(E)$ has a canonical flat connection which we
denote  $\nabla^{can}$ . A choice of $\sigma\in \Isom_0({\mathcal
L}\otimes{\mathcal J}_M,{\mathcal J}({\mathcal L}))$ induces the
flat connection $\sigma^{-1}\circ\nabla^{can}_{\mathcal
L}\circ\sigma$ on ${\mathcal L}\otimes{\mathcal J}_M$.

Let $\nabla$ be a connection on $\cl$ with the curvature $\theta$.
It gives rise to the connection
$\nabla\otimes\id+\id\otimes\nabla^{can}$ on ${\mathcal
L}\otimes{\mathcal J}_M$.

\begin{lemma}\label{lemma:connections compared}

\begin{enumerate}
\item
Choose $\sigma$, $\nabla$ as above. Then the difference
\begin{equation}\label{difference of connections}
F(\sigma, \nabla)=\sigma^{-1}\circ\nabla^{can}_{\mathcal
L}\circ\sigma -(\nabla\otimes\id + \id\otimes\nabla^{can})
\end{equation}
is an element of $\in \Omega^1 \otimes \shEnd_{\cj_M}(\cl \otimes
\cj_M)\cong \Omega^1 \otimes \cj_M$.

\item
Moreover, $F$ satisfies
\begin{equation} \label{formula F}
 \nabla^{can} F(\sigma, \nabla) +\theta=0
\end{equation}
\end{enumerate}
\end{lemma}
\begin{proof}
We leave the verification of the first claim to the reader. The
flatness of $=\sigma^{-1}\circ\nabla^{can}_{\mathcal L}\circ\sigma$
implies the second claim.
\end{proof}

The following properties of our construction are immediate
\begin{lemma}\label{propp}
\begin{enumerate}

\item Let $\cl_1$ and $\cl_2$ be two line bundles, and $f:\cl_1 \to \cl_2$
an isomorphism. Let $\nabla$ be a connection on $\cl_2$ and $\sigma
\in \Isom_0({\mathcal L_2}\otimes{\mathcal J}_M,{\mathcal
J}({\mathcal L_2}))$.  Then
\begin{equation*}
F(\sigma, \nabla)=F(\Ad f(\sigma), \Ad f(\nabla))
\end{equation*}
\item Let $\cl$ be a line bundle on $M$, $\nabla$  a connection on $\cl$ and $\sigma
\in \Isom_0({\mathcal L}\otimes{\mathcal J}_M,{\mathcal J}({\mathcal
L}))$. Let $f: N \to M$ be a smooth map. Then
\begin{equation*}
f^*F(\sigma, \nabla)=F(f^*\sigma,  f^*\nabla)
\end{equation*}

\item Let $\cl$ be a line bundle on $M$, $\nabla$  a connection on $\cl$ and $\sigma
\in \Isom_0({\mathcal L}\otimes{\mathcal J}_M,{\mathcal J}({\mathcal
L}))$. Let $ \phi \in \Gamma(M;\cj_{M,0})$. Then
\begin{equation*}
F(\phi \cdot \sigma, \nabla)=F(\sigma,  \nabla)+\nabla^{can}\phi
\end{equation*}

\item Let $\cl_1$, $\cl_2$ be two line bundles with connections $\nabla_1$ and $\nabla_2$ respectively,
and let $\sigma_i\in \Isom_0({\mathcal L_i}\otimes{\mathcal
J}_M,{\mathcal J}({\mathcal L_i}))$, $i=1, 2$.  Then
\begin{equation*}
F(\sigma_1\otimes \sigma_2, \nabla_1 \otimes \id +\id \otimes
\nabla_2)=F(\sigma_1, \nabla_1)+F(\sigma_2, \nabla_2)
\end{equation*}
\end{enumerate}

\end{lemma}

\subsection{DGLAs of infinite jets}
Suppose that $(\mathcal{U},\underline{\ca})$ is a descent datum
representing a twisted form of $\co_X$. Thus, we have the matrix
algebra $\Mat(\underline{\ca})$ and the cosimplicial DGLA
$\mathfrak{G}(\uca)$ of local $\mathbb{C}$-linear Hochschild
cochains.

The descent datum $(\mathcal{U},\underline{\ca})$ gives rise to the
descent datum $(\mathcal{U},\mathcal{J}(\underline{\ca}))$,
$\mathcal{J}(\underline{\ca}) = (\mathcal{J}(\mathcal{A}),
\mathcal{J}(\mathcal{A}_{01}), j^\infty(\mathcal{A}_{012}))$,
representing a twisted form of $\mathcal{J}_X$, hence to the matrix
algebra $\Mat(\cj(\underline{\ca}))$ and the corresponding
cosimplicial DGLA $\mathfrak{G}(\cj(\uca))$  of local
\emph{$\mathcal{O}$-linear continuous} Hochschild cochains.

The canonical flat connection $\nabla^{can}$ on
$\mathcal{J}(\mathcal{A})$ induces the flat connection, still
denoted $\nabla^{can}$ on $\Mat(\mathcal{J}(\uca))^p$ for each $p$;
the product on $\Mat(\mathcal{J}(\uca))^p$ is horizontal with
respect to $\nabla^{can}$. The flat connection $\nabla^{can}$
induces the flat connection, still denoted $\nabla^{can}$ on
$C^\bullet(\Mat(\mathcal{J}(\uca))^p)^{loc}[1]$ which acts by
derivations of the Gerstenhaber bracket and commutes with the
Hochschild differential $\delta$. Therefore we have the sheaf of
DGLA $\DR( C^\bullet(\Mat(\mathcal{J}(\uca))^p)^{loc}[1])$ with the
underlying sheaf of graded Lie algebras
$\Omega^\bullet_{N_p\cu}\otimes
C^\bullet(\Mat(\mathcal{J}(\uca))^p)^{loc}[1]$ and the differential
$\nabla^{can} + \delta$.

For $\lambda : [n]\to\Delta$ let
\begin{equation*}
\mathfrak{G}_\DR(\cj(\uca))^\lambda = \Gamma(N_{\lambda(n)}\cu;
\lambda(0n)_*
\DR(C^\bullet(\Mat(\mathcal{J}(\uca))^{\lambda(0)})^{loc}[1]))
\end{equation*}
be the DGLA of global sections. The ``inclusion of horizontal sections'' map induces
the morphism of DGLA
\begin{equation*}
j^\infty : \mathfrak{G}(\uca)^\lambda \to
 \mathfrak{G}_\DR(\cj(\uca))^\lambda
\end{equation*}

For $\phi : [m]\to[n]$ in $\Delta$, $\mu = \lambda\circ\phi$ there
is a morphism of DGLA $\phi_*:\mathfrak{G}_\DR(\cj(\uca))^\mu\to
\mathfrak{G}_\DR(\cj(\uca))^\lambda$ making the diagram
\begin{equation*}
\begin{CD}
\mathfrak{G}(\uca)^\lambda @>{j^\infty}>> \mathfrak{G}_\DR(\cj(\uca))^\lambda \\
@V{\phi_*}VV @VV{\phi_*}V \\\
\mathfrak{G}(\uca)^\mu @>{j^\infty}>>
\mathfrak{G}_\DR(\cj(\uca))^\mu
\end{CD}
\end{equation*}
commutative.

Let $\mathfrak{G}_\DR(\cj(\uca))^n =
\prod_{[n]\xrightarrow{\lambda}\Delta}
\mathfrak{G}_\DR(\cj(\uca))^\lambda$. The assignment $\Delta\ni
[n]\mapsto \mathfrak{G}_\DR(\cj(\uca))^n$, $\phi\mapsto\phi_*$
defines the cosimplicial DGLA $\mathfrak{G}_\DR(\cj(\uca))$.

\begin{prop}
The map $j^\infty : \mathfrak{G}(\uca)\to \mathfrak{G}(\cj(\uca))$
extends to a quasiisomorphism of cosimplicial DGLA.
\begin{equation*}
j^\infty : \mathfrak{G}(\uca) \to
 \mathfrak{G}_\DR(\cj(\uca))
\end{equation*}
\end{prop}

The  goal of this section is to construct a quasiisomorphism of the
latter DGLA with the simpler DGLA.

The canonical flat connection $\nabla^{can}$ on $\mathcal{J}_X$
induces a flat connection on $\overline
C^\bullet(\mathcal{J}_X)[1]$, the complex of
\emph{$\mathcal{O}$-linear continuous normalized} Hochschild
cochains, still denoted $\nabla^{can}$ which acts by derivations of
the Gerstenhaber bracket and commutes with the Hochschild
differential $\delta$. Therefore we have the sheaf of DGLA
$\DR(\overline C^\bullet(\mathcal{J}_X)[1])$ with the underlying
graded Lie algebra $\Omega^\bullet_X\otimes \overline
C^\bullet(\mathcal{J}_X)[1]$ and the differential $\delta +
\nabla^{can}$.

Recall that the Hochschild differential $\delta$ is zero on
$C^0(\mathcal{J}_X)$ due to commutativity of $\mathcal{J}_X$. It
follows that the action of the sheaf of abelian Lie algebras
$\mathcal{J}_X = C^0(\mathcal{J}_X)$ on $\overline
C^\bullet(\mathcal{J}_X)[1]$ via the restriction of the adjoint
action (by derivations of degree $-1$)  commutes with the Hochschild
differential $\delta$. Since the cochains we consider are
$\mathcal{O}_X$-linear, the subsheaf $\mathcal{O}_X =
\mathcal{O}_X\cdot j^{\infty}(1)\subset\mathcal{J}_X$ acts
trivially. Hence the action of $\mathcal{J}_X$ descends to an action
of the quotient $\mathcal{J}_X/\mathcal{O}_X$. This action  induces
the action of the abelian graded Lie algebra
$\Omega^\bullet_X\otimes\mathcal{J}_X/\mathcal{O}_X$  by derivations
on the graded Lie algebra $\Omega^\bullet_X\otimes \overline
C^\bullet(\mathcal{J}_X)[1]$. Moreover, the subsheaf
$(\Omega^\bullet_X\otimes\mathcal{J}_X/\mathcal{O}_X)^{cl} :=
\ker(\nabla^{can})$ acts by derivation which commute with the
differential $\delta + \nabla^{can}$. For $\omega\in
\Gamma(X;(\Omega^2_X\otimes\mathcal{J}_X/\mathcal{O}_X)^{cl})$ we
denote by $\DR(\overline C^\bullet(\mathcal{J}_X)[1])_\omega$ the
sheaf of DGLA with the underlying graded Lie algebra
$\Omega^\bullet_X\otimes \overline C^\bullet(\mathcal{J}_X)[1]$ and
the differential $\delta + \nabla^{can} + \iota_\omega$. Let
\begin{equation*}
\mathfrak{g}_\DR(\mathcal{J}_X)_\omega = \Gamma(X;\DR(\overline
C^\bullet(\mathcal{J}_X)[1])_\omega) \ ,
\end{equation*}
be the corresponding  DGLA of global sections.

Suppose now that $\mathcal{U}$ is a cover of $X$; let $\epsilon :
N\mathcal{U}\to X$ denote the canonical map. For $\lambda :
[n]\to\Delta$ let
\begin{equation*}
\mathfrak{G}_\DR(\mathcal{J})^\lambda_\omega =
\Gamma(N_{\lambda(n)}\mathcal{U}; \epsilon^*\DR(\overline
C^\bullet(\mathcal{J}_X)[1])_\omega)
\end{equation*}
For $\mu : [m]\to\Delta$ and a morphism $\phi : [m]\to [n]$ in $\Delta$ such that $\mu = \lambda\circ\phi$ the map $\mu(m)\to\lambda(n)$ induces the map
\begin{equation*}
\phi_* : \mathfrak{G}_\DR(\mathcal{J})^\mu_\omega \to
\mathfrak{G}_\DR(\mathcal{J})^\lambda_\omega
\end{equation*}
For $n = 0, 1, \ldots$ let $\mathfrak{G}_\DR(\mathcal{J})^n_\omega =
\prod_{[n]\xrightarrow{\lambda}\Delta}
\mathfrak{G}_\DR(\mathcal{J})^\lambda_\omega$. The assignment
$[n]\mapsto \mathfrak{G}_\DR(\mathcal{J})^n_\omega$ extends to a
cosimplicial DGLA $\mathfrak{G}_\DR(\mathcal{J})_\omega$. We will
also denote this DGLA by $\mathfrak{G}_\DR(\mathcal{J})_\omega(\cu)$
if we need to explicitly indicate the cover.

\begin{lemma}
The cosimplicial DGLA $\mathfrak{G}_\DR(\mathcal{J})_\omega$ is
acyclic, i.e. satisfies the condition \eqref{acyclicity condition}.
\end{lemma}
\begin{proof}
Consider the cosimplicial vector space $V^{\bullet}$ with $V^n
=\Gamma(N_{n}\mathcal{U}; \epsilon^*\DR(\overline
C^\bullet(\mathcal{J}_X)[1])_\omega)$ and the cosimplicial structure
induced by the simplicial structure of $N\cu$. The cohomology of the
complex $(V^{\bullet}, \partial)$ is the \v{C}ech cohomology of
$\cu$ with the coefficients in the soft sheaf of vector spaces
$\Omega^{\bullet}\otimes\overline C^\bullet(\mathcal{J}_X)[1]$ and,
therefore, vanishes in the positive degrees.
$\mathfrak{G}_\DR(\mathcal{J})_\omega$ as a cosimplicial vector
space can be identified with $\widehat{V}^{\bullet}$ in the
notations of Lemma \ref{subdivision}. Hence the result follows from
the Lemma \ref{subdivision}.
\end{proof}

We leave the proof of the following lemma to the reader.
\begin{lemma}\label{lemma: inclusion of kernel}
The map $\epsilon^* : \mathfrak{g}_\DR(\mathcal{J}_X)_\omega
\to\mathfrak{G}_\DR(\mathcal{J})^0_\omega$ induces an isomorphism of
DGLA
\begin{equation*}
\mathfrak{g}_\DR(\mathcal{J}_X)_\omega \cong
\ker(\mathfrak{G}_\DR(\mathcal{J})^0_\omega \rightrightarrows
\mathfrak{G}_\DR(\mathcal{J})_\omega^1)
\end{equation*}
where the two maps on the right are $(\partial^0_0)_*$ and
$(\partial^0_1)_*$.
\end{lemma}

Two previous lemmas together with the Corollary \ref{cor: stacks are
mc ker for acyclic} imply the following:
\begin{prop}\label{rem}
 Let $\omega\in
\Gamma(X;(\Omega^2_X\otimes\mathcal{J}_X/\mathcal{O}_X)^{cl})$ and
let $\mathfrak{m}$ be a commutative nilpotent ring. Then the map
$\epsilon^*$ induces equivalence of groupoids:
\begin{equation*}
\MC^2(\mathfrak{g}_\DR(\mathcal{J}_X)_\omega\otimes
\mathfrak{m})\cong
\Stack(\mathfrak{G}_\DR(\mathcal{J})_\omega\otimes \mathfrak{m})
\end{equation*}
\end{prop}

For $\beta\in\Gamma(X;\Omega^1\otimes\mathcal{J}_X/\mathcal{O}_X)$
there is a canonical isomorphism of cosimplicial DGLA
$\exp(\iota_\beta) : \mathfrak{G}_\DR(\mathcal{J})_{\omega +
\nabla^{can}\beta} \to \mathfrak{G}_\DR(\mathcal{J})_\omega$.
Therefore, $\mathfrak{G}_\DR(\mathcal{J})_\omega$ depends only on
the class of $\omega$ in $H^2_\DR(\mathcal{J}_X/\mathcal{O}_X)$.

The rest of this section is devoted to the proof of the following theorem.

\begin{thm}\label{theorem: quism jet to mat jet}
Suppose that $(\mathcal{U},\underline{\mathcal{A}})$ is a descent
datum representing a twisted form $\mathcal{S}$ of $\mathcal{O}_X$.
There exists a quasi-isomorphism of cosimplicial DGLA
$\mathfrak{G}_\DR(\mathcal{J})_{[\mathcal{S}]} \to
\mathfrak{G}_\DR(\cj(\uca))$.
\end{thm}

\subsection{Quasiisomorphism}
Suppose that $(\mathcal{U},\underline{\mathcal{A}})$ is a descent
datum for a twisted form of $\mathcal{O}_X$. Thus, $\mathcal{A}$ is
identified with $\mathcal{O}_{N_0\mathcal{U}}$ and
$\mathcal{A}_{01}$ is a line bundle on $N_1\mathcal{U}$.

\subsubsection{Multiplicative connections}
Let $\mathcal{C}^\mu(\mathcal{A}_{01})$ denote the set of
connections $\nabla$ on $\mathcal{A}_{01}$ which satisfy
\begin{enumerate}
\item $\Ad\mathcal{A}_{012}((\pr^2_{02})^*\nabla) =
(\pr_{01}^2)^*\nabla\otimes\id + \id\otimes(\pr^2_{12})^*\nabla$

\item $(\pr^0_{00})^*\nabla$ is the canonical flat connection on
$\mathcal{O}_{N_0\mathcal{U}}$.
\end{enumerate}

Let
$\Isom_0^\mu(\mathcal{A}_{01}\otimes\mathcal{J}_{N_1\mathcal{U}},
\mathcal{J}(\mathcal{A}_{01}))$ denote the subset of
$\Isom_0(\mathcal{A}_{01}\otimes\mathcal{J}_{N_1\mathcal{U}},
\mathcal{J}(\mathcal{A}_{01}))$ which consists of $\sigma$ which
satisfy
\begin{enumerate}
\item $\Ad\mathcal{A}_{012}((\pr_{02}^2)^*\sigma) =
(\pr_{01}^2)^*\sigma\otimes(\pr_{12}^2)^*\sigma$

\item $(\pr_{00}^0)^*\sigma = \id$
\end{enumerate}

Note that the vector space $\overline{Z}^1(\mathcal{U};\Omega^1)$ of
cocycles in the \emph{normalized} \v{C}ech complex of the cover
$\cu$ with coefficients in the sheaf of $1$-forms $\Omega^1$ acts on
the set $\mathcal{C}^\mu(\mathcal{A}_{01})$, with the action given
by
\begin{equation}\label{conn}
\alpha\cdot \nabla = \nabla+\alpha
\end{equation}
Here $\nabla \in \mathcal{C}^\mu(\mathcal{A}_{01})$, $\alpha \in
\overline{Z}^1(\mathcal{U};\Omega^1)  \subset \Omega^1(N_1\cu)$.

Similarly, the vector space
$\overline{Z}^1(\mathcal{U};\mathcal{J}_0)$ acts on the set
$\Isom_0^\mu(\mathcal{A}_{01}\otimes\mathcal{J}_{N_1\mathcal{U}},
\mathcal{J}(\mathcal{A}_{01}))$, with the action given as in
Corollary \ref{aff}.

Note that since the sheaves involved are soft, cocycles coincide
with coboundaries:
$\overline{Z}^1(\mathcal{U};\Omega^1)=\overline{B}^1(\mathcal{U};\Omega^1)$,
$\overline{Z}^1(\mathcal{U};\mathcal{J}_0)=\overline{B}^1(\mathcal{U};\mathcal{J}_0)$.

\begin{prop}\label{sigma}
The set $\mathcal{C}^\mu(\mathcal{A}_{01})$ (respectively,
$\Isom_0^\mu(\mathcal{A}_{01}\otimes\mathcal{J}_{N_1\mathcal{U}},
\mathcal{J}(\mathcal{A}_{01}))$) is an affine space with the
underlying vector space being $\overline{Z}^1(\mathcal{U};\Omega^1)$
(respectively, $\overline{Z}^1(\mathcal{U};\mathcal{J}_0)$).
\end{prop}

\begin{proof}
Proofs of the both statements  are completely analogous. Therefore
we explain the proof of the statement concerning
$\Isom_0^\mu(\mathcal{A}_{01}\otimes\mathcal{J}_{N_1\mathcal{U}},
\mathcal{J}(\mathcal{A}_{01}))$ only.

We show first that
$\Isom_0^\mu(\mathcal{A}_{01}\otimes\mathcal{J}_{N_1\mathcal{U}},
\mathcal{J}(\mathcal{A}_{01}))$ is nonempty.
 Choose an arbitrary $\sigma
\in \Isom_0( \ca_{01}\otimes{\mathcal J}_{N_1\cu},{\mathcal
J}(\ca_{ij}))$  such that $(\pr^0_{00})^*\sigma = \id$. Then, by
Corollary \ref{aff}, there exists $c\in  \Gamma(N_2\cu;\cj_0)$ such
that $ c \cdot (\Ad \ca_{012} ((d^1)^*\sigma_{02}))=
(d^0)^*\sigma_{01}\otimes (d^2)^*\sigma_{12}$. It is easy to see
that $ c\in \overline{Z}^2(\cu; \exp \cj_0)$. Since the sheaf $\exp
\cj_0$ is soft, corresponding \v{C}ech cohomology  is trivial.
Therefore, there exists $ \phi \in \overline{C}^1(\cu;  \cj_0)$ such
that $c = \check{\partial}  \phi $. Then, $\phi \cdot \sigma \in
\Isom_0^\mu(\mathcal{A}_{01}\otimes\mathcal{J}_{N_1\mathcal{U}},
\mathcal{J}(\mathcal{A}_{01}))$.

Suppose that  $\sigma, \sigma' \in
\Isom_0^\mu(\mathcal{A}_{01}\otimes\mathcal{J}_{N_1\mathcal{U}},
\mathcal{J}(\mathcal{A}_{01}))$.  By the Corollary \ref{aff}
$\sigma=\phi \cdot \sigma'$ for some uniquely defined $\phi \in
\Gamma(N_2\cu;\cj_0)$. It is easy to see that $\phi \in
\overline{Z}^1(\mathcal{U};\mathcal{J}_0)$.
\end{proof}

We assume from  now on that we have chosen $\sigma \in
\Isom_0^\mu(\mathcal{A}_{01}\otimes\mathcal{J}_{N_1\mathcal{U}},
\mathcal{J}(\mathcal{A}_{01}))$, $\nabla \in
\mathcal{C}^\mu(\mathcal{A}_{01})$. Such a choice defines
$\sigma_{ij}^p \in \Isom_0( \ca_{ij}\otimes{\mathcal
J}_{N_p\cu},{\mathcal J}(\ca_{ij}))$ for every $p$ and $0 \le i, j
\le p$ by $\sigma^p_{ij}=(\pr^p_{ij})^* \sigma$. This collection of
$\sigma^p_{ij}$ induces for every $p$ algebra isomorphism $\sigma^p:
\Mat (\uca)^p \otimes \cj_{N_p \cu} \to \Mat(\cj(\uca))^p$. The
following compatibility holds for these isomorphisms. Let $f: [p]\to
[q]$ be a morphism in $\Delta$. Then the following diagram commutes:
\begin{equation}\label{combs}
\begin{CD}
f_* (\Mat (\uca)^p \otimes \cj_{N_p\cu}) @>{f_*}>>
f^{\sharp}\Mat(\underline{\mathcal{A}}\otimes \cj_{N_q \cu})^q \\
@Vf_*(\sigma^p)VV @VV{f^\sharp(\sigma^p)}V \\
f_*\Mat (\cj(\uca))^p @>{f_*}>> f^{\sharp}\Mat(\cj(\uca))^q
\end{CD}
\end{equation}
Similarly define the connections
$\nabla^p_{ij}=(\pr_{ij}^p)^*\nabla$. For $p=0,1,\ldots$ set
$\nabla^p = \oplus_{i,j=0}^p\nabla_{ij}$; the connections $\nabla^p$
on $\Mat(\uca)^p$ satisfy
\begin{equation}\label{kak nazvat}
f_*\nabla^p =(\Ad f_*) (f^\sharp \nabla^q) \ .
\end{equation}
Note that $F(\sigma, \nabla)\in \Gamma(N_1\cu;
\Omega^1_{N_1\cu}\otimes \cj_{N_1\cu})$ is a cocycle of degree one
in $\check C^\bullet(\cu; \Omega^1_X\otimes\cj_X)$. Vanishing of the
corresponding \v Cech cohomology implies that  there exists $F^0 \in
\Gamma(N_0\cu; \Omega^1_{N_0\cu}\otimes \cj_{N_0\cu})$ such that
\begin{equation}\label{deff0}
(d^1)^*F^0 -(d^0)^*F^0=F (\sigma, \nabla).
\end{equation}

For $p=0,1,\ldots$, $0\leq i\leq p$, let $F^p_{ii} = (\pr_i^p)^*
F^0$; put $F^p_{ij} = 0$ for $i\neq j$. Let
$F^p\in\Gamma(N_p\cu;\Omega^1_{N_p\cu}\otimes\Mat(\uca)^p)\otimes
\cj_{N_p\cu}$ denote the diagonal matrix with components $F^p_{ij}$.
For $f : [p]\to [q]$ we have
\begin{equation}\label{consf}
f_* F^p = f^\sharp F^q \ .
\end{equation}
Then, we obtain the following
equality of connections on $\Mat(\uca)^p\otimes\cj_{N_p\cu}$:
\begin{equation}\label{nablasigma}
(\sigma^p)^{-1}\circ\nabla^{can}\circ\sigma^p =\nabla^p\otimes\id +
\id\otimes\nabla^{can}+\ad F^p
\end{equation}

The matrices $F^p$ also have the following property. Let
$\nabla^{can} F^p$ be the diagonal matrix with the entries
$(\nabla^{can} F^p)_{ii} =\nabla^{can} F^p_{ii}$. Denote by
$\overline{\nabla^{can}F^p}$ the image of $\nabla^{can}F^p$ under
the natural map
$\Gamma(N_p\cu;\Omega^2_{N_p\cu}\otimes\Mat(\uca)^p\otimes
\cj_{N_p\cu}) \to
\Gamma(N_p\cu;\Omega^2_{N_p\cu}\otimes\Mat(\uca)^p\otimes
(\cj_{N_p\cu}/\co_{N_p\cu}))$. Recall the canonical map $\epsilon_p:
N_p\cu \to X$. Then, we have the following:
\begin{lemma}\label{defw}
There exists a unique $\omega \in \Gamma(X;(\Omega^2_X\otimes\mathcal{J}_X/\mathcal{O}_X)^{cl})$
 such that
\begin{equation*} \overline{\nabla^{can}F^p}=-\epsilon_p^* \omega \otimes \id_p,
\end{equation*}
where $\id_p$ denotes the  $(p+1)\times(p+1)$ identity matrix.
\end{lemma}
\begin{proof}
Using the definition of $F^0$ and formula \eqref{formula F} we
obtain: $(d^1)^*\nabla^{can}F^0
-(d^0)^*\nabla^{can}F^0=\nabla^{can}F (\sigma, \nabla) \in \Omega^2
(N_1\cu)$. Therefore
$(d^1)^*\overline{\nabla^{can}F^0}-(d^0)^*\overline{\nabla^{can}F^0}=0$,
and there exists a unique $\omega \in \Gamma(X; \Omega^2_X \otimes
(\cj_X/\co_X))$ such that $\epsilon_0^* \omega
=\overline{\nabla^{can}F^0}$. Since $(\nabla^{can})^2=0$ it follows
that $\nabla^{can}\omega=0$. For any $p$ we have:
$(\overline{\nabla^{can}F^p})_{ii}
=\pr_i^*\overline{\nabla^{can}F^0}= \epsilon_p^* \omega$, and the
assertion of the Lemma follows.
\end{proof}

\begin{lemma}
The class of $\omega$ in
$H^2(\Gamma(X;\Omega^\bullet_X\otimes\mathcal{J}_X/\mathcal{O}_X),\nabla^{can})$
does not depend on the choices made in its construction.
\end{lemma}
\begin{proof}
The construction of $\omega$ depends on the choice of $\sigma \in
\Isom_0^\mu(\mathcal{A}_{01}\otimes\mathcal{J}_{N_1\mathcal{U}},
\mathcal{J}(\mathcal{A}_{01}))$, $\nabla \in
\mathcal{C}^\mu(\mathcal{A}_{01})$ and $F^0$ satisfying the equation
\eqref{deff0}. Assume that we make  different choices: $\sigma'
=(\check{\partial} \phi) \cdot \sigma$,
$\nabla'=(\check{\partial}\alpha) \cdot \nabla$ and $(F^0)'$
satisfying $\check{\partial}(F^0)'=F (\sigma', \nabla')$. Here,
$\phi\in \overline{C}^0(\mathcal{U};\mathcal{J}_0)$ and $\alpha \in
\overline{C}^0(\mathcal{U};\Omega^1)$. We have:  $F (\sigma',
\nabla')= F (\sigma, \nabla)- \check{\partial} \alpha
+\check{\partial}\nabla^{can} \phi$. It follows that
$\check{\partial}((F^0)'-F^0-\nabla^{can}\phi+\alpha)=0$. Therefore
$(F^0)'-F^0-\nabla^{can}\phi+\alpha =-\epsilon_0^* \beta$ for some
$\beta \in \Gamma(X; \Omega^1_X \otimes \cj_X)$. Hence if $\omega'$
is constructed using $\sigma'$, $\nabla'$, $(F^0)'$ then
$\omega'-\omega =\nabla^{can} \overline{\beta}$ where
$\overline{\beta}$ is the image of $\beta$ under the natural
projection $\Gamma(X; \Omega^1_X \otimes \cj_X) \to \Gamma(X;
\Omega^1_X \otimes (\cj_X/\co_X))$.
\end{proof}

Let $\rho: \cv \to \cu$ be a refinement of the cover $\cu$, and
$(\cv, \uca^{\rho})$ the corresponding descent datum. Choice of
$\sigma$, $\nabla$, $F^0$ on $\cu$ induces the corresponding choice
$(N\rho)^*\sigma$, $(N\rho)^*\nabla$, $(N\rho)^*(F^0)$ on $\cv$. Let
$\omega^{\rho}$ denotes the form constructed as  in Lemma \ref{defw}
using $(N\rho)^*\sigma$, $(N\rho)^*\nabla$, $(N\rho)^*F^0$. Then,
\begin{equation*}
\omega^{\rho}= (N\rho)^* \omega.
\end{equation*}

The following result now follows easily and we leave the details to
the reader.

\begin{prop}\label{propw} The class of $\omega$ in
$H^2(\Gamma(X;\Omega^\bullet_X\otimes\mathcal{J}_X/\mathcal{O}_X),\nabla^{can})$
coincides with the image $[\mathcal{S}]$ of the class of the gerbe.
\end{prop}

\subsection{Construction of the quasiisomorphism}
For $\lambda: [n] \to \Delta$ let
\begin{equation*}
\mathfrak{H}^\lambda := \Gamma(N_{\lambda(n)}\mathcal{U};\
\Omega^{\bullet}_{N_{\lambda(n)}\cu}\otimes \lambda(0n)_*
C^\bullet(\Mat(\underline{\mathcal{A}})^{\lambda(0)} \otimes \cj_{N_{\lambda(0)}\cu})^{loc}[1])
\end{equation*}
considered as a graded Lie algebra.
For $\phi:[m]\to[n]$, $\mu = \lambda\circ\phi$ there is a morphism of graded Lie algebras $\phi_* : \mathfrak{H}^\mu\to \mathfrak{H}^\lambda$.
For $n=0,1,\ldots$ let $\mathfrak{H}^n := \prod_{[n]\xrightarrow{\lambda}\Delta} \mathfrak{H}^\lambda$. The assignment $\Delta\ni [n]\mapsto \mathfrak{H}^n$, $\phi\mapsto \phi_*$ defines a cosimplicial graded Lie algebra $\mathfrak{H}$.

For each $\lambda: [n] \to \Delta$ the map
\begin{equation*}
\sigma^\lambda_* := \id \otimes \lambda(0n)_*(\sigma^{\lambda(0)}):
\mathfrak{H}^\lambda\to \mathfrak{G}_\DR(\cj(\uca))^\lambda
\end{equation*}
is an isomorphism of graded Lie algebras. It follows from \eqref{combs} that the
maps $\sigma^\lambda_*$ yield an isomorphism of cosimplicial graded Lie algebras
\begin{equation*}
\sigma_*: \mathfrak{H}\to \mathfrak{G}_{\DR}(\cj(\uca)).
\end{equation*}
Moreover, the equation
\eqref{nablasigma} shows that if we equip $\mathfrak{H}$ with the
differential given on $\Gamma(N_{\lambda(n)}\mathcal{U};\
\Omega^{\bullet}_{N_{\lambda(n)}\cu}\otimes \lambda(0n)_*
C^\bullet(\Mat(\underline{\mathcal{A}})^{\lambda(0)} \otimes
\cj_{N_{\lambda(0)}\cu})^{loc}[1])$ by
\begin{multline}\label{differ}
\delta+\lambda(0n)_*(\nabla^{\lambda(0)})\otimes\id +
\id\otimes\nabla^{can}+\ad
\lambda(0n)_*(F^{\lambda(0)})\\=\delta+\lambda(0n)^{\sharp}(\nabla^{\lambda(n)})\otimes\id
+ \id\otimes\nabla^{can}+\ad \lambda(0n)^{\sharp}(F^{\lambda(n)}))
\end{multline}
then $\sigma_*$ becomes an isomorphism of DGLA. Consider now an
automorphism $ \exp \iota_F$ of the cosimplicial graded Lie algebra
$\mathfrak{H}$ given on $\mathfrak{H}^{\lambda}$ by $\exp
\iota_{\lambda(0n)_* F_p}$. Note the fact that this morphism
preserves  the cosimplicial structure follows from the relation
\eqref{consf}.

The following result is proved by the direct calculation; see
\cite{bgnt2}, Lemma 16.
\begin{lemma}
\begin{multline}\label{conj}
\exp (\iota_{ F_p})\circ (\delta+\nabla^p\otimes\id +
\id\otimes\nabla^{can}+\ad F^p) \circ\exp (-\iota_{F_p}) =\\
\delta+\nabla^p\otimes\id + \id\otimes\nabla^{can} - \iota_{\nabla
F^p}.
\end{multline}
\end{lemma}

Therefore the morphism
\begin{equation}
\exp \iota_F :\mathfrak{H} \to \mathfrak{H}
\end{equation}
conjugates the differential given by the formula \eqref{differ} into
the differential which on $\mathfrak{H}^{\lambda}$ is given by
\begin{equation}\label{2differential}
\delta+\lambda(0n)_*(\nabla^{\lambda(0)})\otimes\id +
\id\otimes\nabla^{can}-\iota_{\lambda(0n)_*(\nabla^{can}F^{\lambda(0)})}
\end{equation}

Consider the map
\begin{equation}\label{cotrace map}
\cotr:   \overline C^\bullet({\mathcal J}_{N_p\cu})[1]\to
  C^\bullet(\Mat(\uca)^p\otimes{\mathcal J}_{N_p\cu})[1]
\end{equation}
defined as follows:
\begin{equation}\label{define cotrace}
\cotr(D)(a_1\otimes j_1, \dots, a_n \otimes j_n) = a_0\ldots a_n
D(j_1, \ldots, j_n).
\end{equation}

The map $\cotr$ is a quasiisomorphism of DGLAs (cf. \cite{Loday},
section 1.5.6; see also \cite{bgnt2} Proposition 14).

\begin{lemma}
For every $p$ the map
\begin{equation}\label{cotrace with forms}
\id\otimes\cotr: \Omega^\bullet_{N_p\cu}\otimes\overline
C^\bullet({\mathcal J}_{N_p\cu})[1]\to
\Omega^\bullet_{N_p\cu}\otimes
C^\bullet(\Mat(\uca)^p\otimes{\mathcal J}_{N_p\cu})[1] \ .
\end{equation}
is a quasiisomorphism of DGLA, where the source and the target are
equipped with the differentials
$\delta+\nabla^{can}+\iota_{\epsilon^* \omega}$  and
$\delta+\nabla^p\otimes\id + \id\otimes\nabla^{can} - \iota_{\nabla
F^p}$ respectively.
\end{lemma}
\begin{proof}
It is easy to see that $\id \otimes \cotr$ is a morphism of graded
Lie algebras, which satisfies  $ (\nabla^p\otimes\id +
\id\otimes\nabla^{can} )\circ (\id \otimes \cotr)= (\id \otimes
\cotr) \circ \nabla^{can} $ and $\delta \circ (\id \otimes \cotr)=
(\id \otimes \cotr) \circ \delta $. Since the domain of $(\id
\otimes \cotr)$ is the normalized complex, in view of the Lemma
\ref{defw} we also have $\iota_{\nabla F^p} \circ(\id \otimes
\cotr)=-(\id \otimes \cotr) \circ \iota_{\epsilon^*\omega}$. This
implies that $(\id \otimes \cotr)$ is a morphism of DGLA.

To see that this map is a quasiisomorphism, introduce filtration on
$\Omega^{\bullet}_{N_p\cu}$ by
$F_i\Omega^\bullet_{N_p\cu}=\Omega^{\geq -i}_{N_p\cu}$ and consider
the complexes $\overline C^\bullet({\mathcal J}_{N_p\cu})[1]$ and
$C^\bullet(\Mat(\uca))^p\otimes{\mathcal J}_{N_p\cu})[1]$ equipped
with the trivial filtration. The map \eqref{cotrace with forms} is a
morphism of filtered complexes with respect to the induced
filtrations on the source and the target. The differentials induced
on the associated graded complexes are $\delta$ (or, more precisely,
$\id\otimes\delta$) and the induced map of the associated graded
objects is $\id\otimes\cotr$ which is a quasi-isomorphism.
Therefore, the map \eqref{cotrace with forms} is a quasiisomorphism
as claimed.
\end{proof}
The map \eqref{cotrace with forms} therefore induces  for every
$\lambda:[n] \to \Delta$ a morphism $\id \otimes \cotr:
\mathfrak{G}_\DR(\mathcal{J})^\lambda_\omega \to
\mathfrak{H}^{\lambda}$. These morphisms are clearly compatible with
the cosimplicial structure and hence induce a quasiisomorphism of
cosimplicial DGLAs
\begin{equation*}
\id \otimes \cotr: \mathfrak{G}_\DR(\mathcal{J})_\omega \to
\mathfrak{H}
\end{equation*}
where the differential in the right hand side is given by
\eqref{2differential}.

We summarize our consideration in the following:
\begin{thm}\label{prop:quism of DGLA}
For a any choice of $\sigma \in
\Isom_0^\mu(\mathcal{A}_{01}\otimes\mathcal{J}_{N_1\mathcal{U}},
\mathcal{J}(\mathcal{A}_{01}))$, $\nabla \in
\mathcal{C}^\mu(\mathcal{A}_{01})$
  and $F^0$ as in \eqref{deff0} the
composition $\Phi_{\sigma,\nabla, F}:=
\sigma_*\circ\exp(\iota_{F})\circ(\id\otimes\cotr)$
\begin{equation}\label{map:the comparison map}
\Phi_{\sigma,\nabla, F}: \mathfrak{G}_\DR(\mathcal{J})_\omega \to
\mathfrak{G}_{\DR}(\cj(\uca))
\end{equation}
is a quasiisomorphism of cosimplicial DGLAs.
\end{thm}

Let $\rho: \cv \to \cu$ be a refinement of the cover $\cu$ and let
$(\cv, \uca^{\rho})$ be the induced descent datum. We will denote
the corresponding cosimplicial DGLAs by
$\mathfrak{G}_\DR(\mathcal{J})_\omega(\cu)$ and
$\mathfrak{G}_\DR(\mathcal{J})_\omega(\cv)$ respectively. Then the
map $N\rho: N\cv \to N\cu$ induces a morphism of cosimplicial DGLAs
\begin{equation}\label{nro}
(N\rho)^*: \mathfrak{G}_{\DR}(\cj(\uca)) \to
\mathfrak{G}_{\DR}(\cj(\uca^{\rho}))
\end{equation}and
\begin{equation}
(N\rho)^*: \mathfrak{G}_\DR(\mathcal{J})_\omega(\cu) \to
\mathfrak{G}_\DR(\mathcal{J})_\omega(\cv).
\end{equation}
Notice also that the choice that the choice of data $\sigma$,
$\nabla$, $F^0$ on $N \cu$ induces the corresponding data
$(N\rho)^*\sigma$, $(N\rho)^*\nabla$, $(N\rho)^*F^0$ on $N \cv$.
This data allows one to construct using the equation \eqref{map:the
comparison map} the map
\begin{equation}  \Phi_{(N\rho)^*\sigma,(N\rho)^*\nabla,
(N\rho)^*F}: \mathfrak{G}_\DR(\mathcal{J})_\omega \to
\mathfrak{G}_{\DR}(\cj(\uca^\rho))
\end{equation}

The following Proposition is an easy consequence of  the description
of the map $\Phi$, and we leave the proof to the reader.
\begin{prop}\label{prop: refinement commutes with construction}
The following diagram commutes:
\begin{equation}\hspace{-2cm}
\begin{CD}
\mathfrak{G}_\DR(\mathcal{J})_\omega(\cu) @>{(N\rho)^*}>>
\mathfrak{G}_\DR(\mathcal{J})_\omega(\cv) \\
@V\Phi_{\sigma,\nabla, F}VV
@VV{\Phi_{(N\rho)^*\sigma,(N\rho)^*\nabla,
(N\rho)^*F}}V \\
\mathfrak{G}_{\DR}(\cj(\uca)) @>{(N\rho)^*}>>
\mathfrak{G}_{\DR}(\cj(\uca^{\rho}))
\end{CD}
\end{equation}
\end{prop}

\section{Proof of the main theorem}\label{section: proof of main
theorem} In this section we prove the main result of this paper.
Recall the statement of the theorem from the introduction:
\begin{thm1}
Suppose that $X$ is a $C^\infty$ manifold and $\mathcal{S}$ is an
algebroid stack on $X$ which is a twisted form of $\mathcal{O}_X$.
Then, there is an equivalence of $2$-groupoid valued functors of
commutative Artin $\mathbb{C}$-algebras
\begin{equation*}
\Def_X(\mathcal{S})\cong
\MC^2(\mathfrak{g}_\DR(\mathcal{J}_X)_{[\mathcal{S}]}) \ .
\end{equation*}
\end{thm1}
\begin{proof}
 Suppose   $\mathcal{U}$ is a cover of $X$ such
that $\epsilon_0^*\mathcal{S}(N_0\cu)$ is nonempty. There is a
descent datum $(\cu,\uca)\in\Desc_\mathbb{C}(\cu)$ whose image under
the functor $\Desc_\mathbb{C}(\cu)\to\AlgStack_\mathbb{C}(X)$ is
equivalent to $\mathcal{S}$.

The proof proceeds as follows. Recall the $2$-groupoids
$\Def'(\cu,\uca)(R)$ and $\Def(\cu,\uca)(R)$ of deformations of and
star-products on the descent datum $(\cu,\uca)$ defined in the
Section \ref{defdescdata}. Note that the composition
$\Desc_R(\cu)\to \Triv_R(X)\to \AlgStack_R(X)$ induces functors
$\Def'(\cu,\uca)(R)\to\Def(\mathcal{S})(R)$ and
$\Def(\cu,\uca)(R)\to\Def(\mathcal{S})(R)$, the second one being the
composition of the first one with the equivalence
$\Def(\cu,\uca)(R)\to \Def'(\cu,\uca)(R)$.
 We are going to show that for a commutative Artin
$\mathbb{C}$-algebra $R$ there are equivalences
\begin{enumerate}
\item $\Def(\cu,\uca)(R)\cong
\MC^2(\mathfrak{g}_\DR(\mathcal{J}_X)_{[\mathcal{S}]})(R)$ and

\item the functor
$\Def(\cu,\uca)(R)\cong \Def(\mathcal{S})(R)$ induced by the functor
$\Def(\cu,\uca)(R)\to \Def(\mathcal{S})(R)$ above.
\end{enumerate}

Let $\mathcal{J}(\uca) = (\mathcal{J}(\ca_{01}),
j^\infty(\ca_{012}))$. Then, $(\cu, \mathcal{J}(\uca))$ is a descent
datum for a twisted form of $\mathcal{J}_X$. Let $R$ be a
commutative Artin $\mathbb{C}$-algebra with maximal ideal
$\mathfrak{m}_R$.

Then the first statement follows from the equivalences
\begin{eqnarray}
\Def(\cu,\uca)(R)  & \cong &
\Stack(\mathfrak{G}(\uca)\otimes\mathfrak{m}_R)
\label{to all stacks}\\
& \cong & \Stack(\mathfrak{G}_\DR(\cj(\uca)))\otimes\mathfrak{m}_R)
\label{to jets}\\
& \cong &
\Stack(\mathfrak{G}_\DR(\mathcal{J}_X)_{[\mathcal{S}]}\otimes\mathfrak{m}_R)
\label{effect of construction}\\
& \cong &
\MC^2(\mathfrak{g}_\DR(\mathcal{J}_X)_{[\mathcal{S}]}\otimes\mathfrak{m}_R)
\label{inclusion of kernel}
\end{eqnarray}

Here the equivalence \eqref{to all stacks} is the subject of the
Proposition \ref{allst}. The inclusion of horizontal sections is a
quasi-isomorphism and the induced map in \eqref{to jets} is an
equivalence by Theorem \ref{quism invariance of stack}. In the
Theorem \ref{prop:quism of DGLA} we have constructed a
quasi-isomorphism $\mathfrak{G}(\mathcal{J}_X)_{[\mathcal{S}]} \to
\mathfrak{G}(\cj(\uca))$; the induced map \eqref{effect of
construction} is an equivalence by another application of Theorem
\ref{quism invariance of stack}. Finally, the equivalence
\eqref{inclusion of kernel} is shown in the Proposition \ref{rem}

We now proceed with the proof of the second statement.
We begin by considering the behavior of $\Def'(\cu,\uca)(R)$ under the refinement.
Consider a
refinement $\rho : \mathcal{V} \to \mathcal{U}$ . Recall that by the
Proposition \ref{rem} the map $\epsilon^*$ induces equivalences
$\MC^2(\mathfrak{g}_\DR(\mathcal{J}_X)_\omega\otimes
\mathfrak{m})\to
\Stack(\mathfrak{G}_\DR(\mathcal{J})_\omega(\cu)\otimes
\mathfrak{m})$ and
$\MC^2(\mathfrak{g}_\DR(\mathcal{J}_X)_\omega\otimes
\mathfrak{m})\to
\Stack(\mathfrak{G}_\DR(\mathcal{J})_\omega(\cv)\otimes
\mathfrak{m})$. It is clear that the diagram
\begin{equation*}
\xymatrix@C=-1cm@R=1cm{
  &\MC^2(\mathfrak{g}_\DR(\mathcal{J}_X)_\omega\otimes
\mathfrak{m}) \ar[rd] \ar[ld]&
  \\
\Stack(\mathfrak{G}_\DR(\mathcal{J})_\omega(\cu)\otimes
\mathfrak{m})\ar[rr]^{(N \rho)^*} &
  &\Stack(\mathfrak{G}_\DR(\mathcal{J})_\omega(\cv)\otimes
\mathfrak{m}) }
\end{equation*}
commutes, and therefore $(N\rho)^*:
\Stack(\mathfrak{G}_\DR(\mathcal{J})_\omega(\cu)\otimes
\mathfrak{m})\to
  \Stack(\mathfrak{G}_\DR(\mathcal{J})_\omega(\cv)\otimes
\mathfrak{m})$ is an equivalence. Then the  Proposition \ref{prop:
refinement commutes with construction} together with the Theorem
\ref{quism invariance of stack} implies that $(N\rho)^*:
\Stack(\mathfrak{G}_{\DR}(\cj(\uca)))\otimes \mathfrak{m})\to
  \Stack(\mathfrak{G}_{\DR}(\cj(\uca^{\rho})))\otimes
\mathfrak{m})$ is an equivalence. It follows that the functor
$\rho^* : \Def(\cu,\uca)(R)\to \Def(\mathcal{V},\uca^\rho)(R)$ is an
equivalence. Note also that the diagram
\[
\begin{CD}
\Def(\cu,\uca)(R) @>{\rho^*}>> \Def(\mathcal{V},\uca^\rho)(R) \\
@VVV @VVV \\
\Def'(\cu,\uca)(R) @>{\rho^*}>> \Def'(\mathcal{V},\uca^\rho)(R)
\end{CD}
\]
is commutative with the top horizontal and both vertical maps being
equivalences. Hence it follows that the bottom horizontal map is an
equivalence.

Recall now that by the Proposition \ref{star products are equivalent to
def} the embedding $\Def(\cu,\uca)(R)\to \Def'(\cu,\uca)(R)$ is an
equivalence. Therefore  it is sufficient to show that the functor
$\Def'(\cu,\uca)(R)\to \Def(\mathcal{S})(R)$ is an equivalence.
Suppose that $\mathcal{C}$ is an $R$-deformation of $\mathcal{S}$.
It follows from Lemma \ref{lemma: def of algd is algd} that
$\mathcal{C}$ is an algebroid stack. Therefore, there exists a cover
$\mathcal{V}$ and an $R$-descent datum
$(\mathcal{V},\underline{\mathcal{B}})$ whose image under the
functor $\Desc_R(\mathcal{V})\to\AlgStack_R(X)$ is equivalent to
$\mathcal{C}$. Replacing $\mathcal{V}$ by a common refinement of
$\mathcal{U}$ and $\mathcal{V}$ if necessary we may assume that
there is a morphism of covers $\rho: \mathcal{V}\to\mathcal{U}$.
Clearly, $(\mathcal{V},\underline{\mathcal{B}})$ is a deformation of
$(\mathcal{V}, \uca^\rho)$. Since the functor $\rho^* :
\Def'(\cu,\uca)(R)\to \Def'(\mathcal{V},\uca^\rho)$ is an
equivalence there exists a deformation $(\cu,
\underline{\mathcal{B}}')$ such that
$\rho^*(\mathcal{U},\underline{\mathcal{B}}')$ is isomorphic to
$(\mathcal{V},\underline{\mathcal{B}})$. Let $\mathcal{C}'$ denote
the image of $(\cu, \underline{\mathcal{B}}')$ under the functor
$\Desc_R(\mathcal{V})\to\AlgStack_R(X)$. Since the images of $(\cu,
\underline{\mathcal{B}}')$ and
$\rho^*(\mathcal{U},\underline{\mathcal{B}}')$ in $\AlgStack_R(X)$
are equivalent it follows that $\mathcal{C}'$ is equivalent to
$\mathcal{C}$. This shows that the functor
$\Def'(\cu,\uca)(R)\to\Def(\mathcal{C})(R)$ is essentially
surjective.

Suppose now that $(\cu,\underline{\mathcal{B}}^{(i)})$, $i=1,2$, are
$R$-deformations of $(\cu,\uca)$. Let $\mathcal{C}^{(i)}$ denote the
image of $(\cu,\underline{\mathcal{B}}^{(i)})$ in
$\Def(\mathcal{S})(R)$. Suppose that $F :
\mathcal{C}^{(1)}\to\mathcal{C}^{(2)}$ is a $1$-morphism.

Let $L$ denote the image of the canonical trivialization under the composition
\[
\epsilon_0^*(F) : \widetilde{\mathcal{B}^{(1)+}}\cong
\epsilon_0^*\mathcal{C}^{(1)}\xrightarrow{F}
\epsilon_0^*\mathcal{C}^{(1)}
\cong\widetilde{\mathcal{B}^{(2)+}} \ .
\]
Thus, $L$ is a $\mathcal{B}^{(1)}\otimes_R
\mathcal{B}^{(2)op}$-module  such that the line bundle
$L\otimes_R\mathbb{C}$ is trivial. Therefore, $L$ admits a
non-vanishing global section. Moreover, there is an isomorphism $f :
\mathcal{B}^{(1)}_{01}\otimes_{(\mathcal{B}^{(1)})_1^1}(\pr^1_1)^*L
\to
(\pr^1_0)^*L\otimes_{(\mathcal{B}^{(2)})_0^1}\mathcal{B}^{(2)}_{01}$
of $(\mathcal{B}^{(1)})_0^1\otimes_R
((\mathcal{B}^{(2)})^1_1)^{op}$-modules.

A choice of a non-vanishing global section of $L$ gives rise to
isomorphisms $\mathcal{B}^{(1)}_{01}\cong
\mathcal{B}^{(1)}_{01}\otimes_{(\mathcal{B}^{(1)})_1^1}(\pr^1_1)^*L$
and $\mathcal{B}^{(2)}_{01}\cong
(\pr^1_0)^*L\otimes_{(\mathcal{B}^{(2)})_0^1}\mathcal{B}^{(2)}_{01}$.

The composition
\[
\mathcal{B}^{(1)}_{01}\cong
\mathcal{B}^{(1)}_{01}\otimes_{(\mathcal{B}^{(1)})_1^1}(\pr^1_1)^*L
\xrightarrow{f}
(\pr^1_0)^*L\otimes_{(\mathcal{B}^{(2)})_0^1}\mathcal{B}^{(2)}_{01}
\cong \mathcal{B}^{(2)}_{01}
\]
defines a $1$-morphism of deformations of $(\cu,\uca)$ such that the
induced $1$-morphism $\mathcal{C}^{(1)}\to\mathcal{C}^{(2)}$ is
isomorphic to $F$. This shows that the functor
$\Def'(\cu,\uca)(R)\to\Def(\mathcal{C})$ induces essentially
surjective functors on groupoids of morphisms. By similar arguments
left to the reader one shows that these are fully faithful.

This completes the proof of Theorem \ref{main theorem}.
\end{proof}

\end{document}